\documentclass[preprint,1p,leqno]{elsarticle}
\usepackage{amsmath}
\usepackage{amssymb}
\usepackage{mathrsfs}
\usepackage{color}
\usepackage{marvosym}

\usepackage[]{graphicx}
\usepackage{psfrag}
\usepackage{pdfpages}
\usepackage[dvipsnames]{xcolor}

\usepackage[makeroom]{cancel}

\usepackage{amsthm}
\newtheorem{remark}{Remark}

\usepackage{bigstrut}
\usepackage{multirow}
\usepackage{subfig}

\usepackage{hyperref}

\usepackage{geometry}
\newcommand{\bu}{\boldsymbol{u}}
\newcommand{\bw}{\boldsymbol{\omega}}
\newcommand{\be}{\boldsymbol{\epsilon}}
\newcommand{\eR}{\dfrac{1}{\mathrm{Re}}}
\newcommand{\Rn}{\mathrm{Re}}

\usepackage{blindtext}

\newcommand{\REVO}[1]{\textcolor{black}{#1}} 
\newcommand{\REVT}[1]{\textcolor{black}{#1}} 
\newcommand{\MOD}[1]{\textcolor{black}{#1}}  

\journal{Journal of Computational Physics}
\begin{document}
\begin{frontmatter}
\title{A mass-, kinetic energy- and helicity-conserving mimetic dual-field discretization for three-dimensional incompressible Navier-Stokes equations, part I: Periodic domains}
\author[TUD]{Yi Zhang}
\author[TUD]{Artur Palha}
\author[TUD]{Marc Gerritsma}
\address[TUD]{Delft University of Technology, Faculty of Aerospace Engineering, Kluyverweg 1, 2629 HS Delft, the Netherlands}
\author[CU]{Leo G. Rebholz}
\address[CU]{Clemson University, 105 Sikes Hall, Clemson, SC 29634, United States}
\begin{abstract}
We introduce a mimetic dual-field discretization which conserves mass, kinetic energy and helicity for three-dimensional incompressible Navier-Stokes equations. The discretization makes use of a conservative dual-field mixed weak formulation where two evolution equations of velocity are employed and dual representations of the solution are sought for each variable. A temporal discretization, which staggers the evolution equations and handles the nonlinearity such that the resulting discrete algebraic systems are linear and decoupled, is constructed. The spatial discretization is mimetic in the sense that the finite dimensional function spaces form a discrete de Rham complex. Conservation of mass, kinetic energy and helicity in the absence of dissipative terms is proven at the discrete level. Proper dissipation rates of kinetic energy and helicity in the viscous case is also proven. Numerical tests supporting the method are provided.
\end{abstract}
\end{frontmatter}

\section{Introduction}\label{section:introduction}

\subsection{Relevance of structure preserving methods with focus on kinetic energy and helicity conservation}\label{section:introduction:relevance_structure_preserving}
In this work we address the discretization of the incompressible Navier-Stokes equations, defined on a periodic domain $\Omega\subset\mathbb{R}^{3}$ and time interval $(0, t_{F}]$. These well known equations govern the dynamics of an incompressible fluid's velocity, $\boldsymbol{u}:\Omega\times(0, t_{F}] \mapsto \mathbb{R}^{3}$, and \REVO{pressure, $p:\Omega\times(0, t_{F}] \mapsto \mathbb{R}$}, subject to a body force, $\boldsymbol{f}:\Omega\times(0, t_{F}] \mapsto \mathbb{R}^{3}$, and an \MOD{initial} condition, $ \boldsymbol{u}^{0}:\Omega \mapsto \mathbb{R}^{3} $.  A general dimensionless form of these equations is \vspace{-3ex}

\begin{subequations}\label{eq:navier_stokes_continuous}
	\begin{align}
		&\frac{\partial \boldsymbol{u}}{\partial t} + \mathcal{C}(\boldsymbol{u}) - \eR\mathcal{D}(\boldsymbol{u}) + \nabla p = \boldsymbol{f}\,, &&\text{in }\Omega\times(0, t_{F}]\,,\label{eq:navier_stokes_continuous a}\\
		&\nabla\cdot\boldsymbol{u} = 0\,,  &&\text{in } \Omega\times(0, t_{F}]\,, \label{eq:navier_stokes_continuous b}\\
		&\left.\boldsymbol{u}\right|_{t=0} = \boldsymbol{u}^{0},&&\text{in } \Omega\,,
	\end{align} 
\end{subequations}
where $\mathcal{C}(\boldsymbol{u})$ and $\mathcal{D}(\boldsymbol{u})$ represent the nonlinear convective term and the linear dissipative term, respectively, and $\mathrm{Re}$ is the Reynolds number. The operators $\mathcal{C}$ and $\mathcal{D}$ can take different forms, all analytically equivalent at the continuous level, see for example \cite{Zang1991, Ronquist1996, Morinishi1998, Layton2009, Capuano2018, olshanskii2002navier, layton2009accuracy}. 
 
The four most common forms of the nonlinear convective term $\mathcal{C}(\boldsymbol{u})$ present in the literature, e.g., \cite{Morinishi1998, Capuano2018, Palha2017a}, are  \vspace{-3ex}

\begin{subequations}\label{eq: forms for nonlinear convective term}
\begin{align}
	\text{\emph{Advective form:} } & \mathcal{C}(\boldsymbol{u}) := \boldsymbol{u}\cdot\nabla\boldsymbol{u}\,, \\ 
	\text{\emph{Conservative (or divergence) form:} } & \mathcal{C}(\boldsymbol{u}) := \nabla\cdot\left(\boldsymbol{u}\otimes\boldsymbol{u}\right)\,,\\
	\text{\emph{Skew-symmetric form:} } & \mathcal{C}(\boldsymbol{u}) := \frac{1}{2}\nabla\cdot\left(\boldsymbol{u}\otimes\boldsymbol{u}\right) + \frac{1}{2}\boldsymbol{u}\cdot\nabla\boldsymbol{u}\,,\\
	\text{\emph{Rotational (or Lamb) form:} } & \mathcal{C}(\boldsymbol{u}) := \boldsymbol{\omega}\times\boldsymbol{u}\ +\ \frac{1}{2}\nabla\left(\boldsymbol{u}\cdot \boldsymbol{u} \right)  \,, \label{eq: forms for nonlinear convective term rotational}
\end{align}
\end{subequations}
where $\boldsymbol{\omega} := \nabla\times\boldsymbol{u}$ is the vorticity field. Besides these most common forms, it is also possible to construct a wide range of nonlinear convective terms as linear combinations of the above mentioned ones and/or employing vector calculus identities. For example, one such choice with interesting properties is the EMAC scheme \cite{Charnyi2017}. Following similar ideas, it is possible to construct analytically equivalent representations for the dissipative term $\mathcal{D}$, for example,
\begin{equation}\label{Eq: representations for the dissipative term}
	\mathcal{D}(\boldsymbol{u}) := \varDelta \boldsymbol{u}, \qquad \text{and}\qquad \mathcal{D}(\boldsymbol{u}) := -\nabla\times\nabla\times\boldsymbol{u}  = -\nabla\times\boldsymbol{\omega}\,,
\end{equation}
where the latter representation can be derived from the former by using the identity $ \varDelta \boldsymbol{u} = \nabla\left( \nabla\cdot \boldsymbol{u}\right) - \nabla\times\left(\nabla\times\boldsymbol{u}\right) $ and the divergence free condition \eqref{eq:navier_stokes_continuous b}.

As mentioned before, these different forms are equivalent at the continuous level and, therefore, may be used interchangeably. At the discrete level, see for example \cite{Morinishi1998, Capuano2018, Palha2017a}, a particular choice of convective term used as the starting point of the discretization process leads to numerical schemes with substantially different properties.

One interesting aspect of the incompressible Navier-Stokes equations \eqref{eq:navier_stokes_continuous} is the fact that, in the inviscid limit ($\Rn \to \infty$) and when the external body force is conservative (there exists a scalar field $ \varphi $ such that $ \boldsymbol{f}=\nabla\varphi $), its dynamics conserves several invariants. Some of these invariants are the total kinetic energy $\mathcal{K}$ (in 2D and 3D), total enstrophy $\mathcal{E}$ (in 2D), and the total helicity $\mathcal{H}$ (in 3D),
\begin{equation} \label{eq:conservation_laws}
	\mathcal{K} := \frac{1}{2}\int_{\Omega} \boldsymbol{u}\cdot\boldsymbol{u}\,, \qquad \mathcal{E} := \frac{1}{2}\int_{\Omega}\boldsymbol{\omega} \cdot \boldsymbol{\omega}, \qquad\text{and}\qquad \mathcal{H} := \int_{\Omega} \boldsymbol{u}\cdot\boldsymbol{\omega}\,,
\end{equation}
provided there is no net in- or out-flow of kinetic energy, enstrophy or helicity over the domain boundary. Note that, in 2D, vorticity can be regarded as a vector field constrained to the direction orthogonal to the planar 2D domain and velocity can be regarded as a vector field whose component along the direction orthogonal to the planar 2D domain is zero, i.e,
$ \boldsymbol{\omega} = \left[ 0,0,\omega \right]^{\mathsf{T}}$ and $\boldsymbol{u} = \left[u,v,0 \right] ^{\mathsf{T}} $.
Thus helicity is trivially zero in 2D flows.

The proofs for these conservation laws are straightforward. For illustration purposes and as an introduction to some of the ideas discussed later in this work, we present these proofs here for the case of no external force, i.e., $\boldsymbol{f} = \boldsymbol{0}$, and periodic boundary condition. For simplicity, and without loss of generality, we use the rotational (or Lamb) form for the nonlinear convective term, \eqref{eq: forms for nonlinear convective term rotational}. The total (or Bernoulli) pressure is defined as
\[
P := p + \dfrac{1}{2}\boldsymbol{u}\cdot\boldsymbol{u}\,.
\]

Kinetic energy conservation (in 2D and 3D) corresponds to $\dfrac{\mathrm{d}\mathcal{K}}{\mathrm{d}t} = 0$. Differentiating $\mathcal{K}$ as defined in \eqref{eq:conservation_laws} with respect to time and taking \eqref{eq:navier_stokes_continuous} in the inviscid limit, $\Rn \to \infty$, leads to
\[
	\frac{\mathrm{d}\mathcal{K}}{\mathrm{d}t} = \int_{\Omega}\frac{\partial\boldsymbol{u}}{\partial t}\cdot\boldsymbol{u} \stackrel{\eqref{eq:navier_stokes_continuous}}{=} -\int_{\Omega}\mathcal{C}(\boldsymbol{u})\cdot\boldsymbol{u} - \int_{\Omega}\nabla p\cdot\boldsymbol{u} = -\int_{\Omega}\left(\boldsymbol{\omega}\times\boldsymbol{u}\right)\cdot\boldsymbol{u} + \int_{\Omega}P\nabla\cdot\boldsymbol{u} = 0\,,
\]
where we have used (i) the vector calculus relation that the cross product of two vectors is perpendicular to either vector, i.e.,
\begin{equation}\label{Eq: perpendicular inner product}
	\left(\boldsymbol{a}\times\boldsymbol{b}\right)\perp\boldsymbol{a}\quad\text{and}\quad \left(\boldsymbol{a}\times\boldsymbol{b}\right)\perp\boldsymbol{b}\,,
\end{equation}
(ii) integration by parts on the total pressure term and (iii) the divergence free condition \eqref{eq:navier_stokes_continuous b}.

Enstrophy conservation (in 2D) equates to $\dfrac{\mathrm{d}\mathcal{E}}{\mathrm{d}t} = 0$. As done above for kinetic energy, time differentiation of $\mathcal{E}$ as defined in \eqref{eq:conservation_laws} gives
\begin{equation}\label{eq:enstrophy_time_differentiation}
	\frac{\mathrm{d}\mathcal{E}}{\mathrm{d}t} = \int_{\Omega} \frac{\partial\boldsymbol{\omega}}{\partial t}\cdot\boldsymbol{\omega}\,. 
\end{equation}
Computing the curl of the momentum equation in \eqref{eq:navier_stokes_continuous} with $\Rn \to \infty$ and substituting $\boldsymbol{\omega} = \nabla\times\boldsymbol{u}$
into \eqref{eq:enstrophy_time_differentiation} results in
\[
	\frac{\mathrm{d}\mathcal{E}}{\mathrm{d}t} = -\int_{\Omega} \nabla\times\left(\boldsymbol{\omega}\times\boldsymbol{u}\right)\cdot\boldsymbol{\omega} = -\frac{1}{2}\int_{\Omega} \left(\boldsymbol{u}\cdot\nabla\omega\right)\omega - \frac{1}{2}\int_{\Omega} \nabla\cdot\left(\boldsymbol{u}\omega\right)\omega =  \frac{1}{2}\int_{\Omega} \omega\nabla\cdot\left(\boldsymbol{u}\omega\right) - \frac{1}{2}\int_{\Omega} \nabla\cdot\left(\boldsymbol{u}\omega\right)\omega = 0\,,
\]
where we first used the vector calculus identity
\[
	\nabla\times\left(\boldsymbol{\omega}\times\boldsymbol{u}\right) = \frac{1}{2}\left(\boldsymbol{u}\cdot\nabla\omega\right) + \frac{1}{2} \nabla\cdot(\boldsymbol{u}\omega)\,,
\]
followed by integration by parts on the first term of the second equality.

Helicity conservation (in 3D) stands for $\dfrac{\mathrm{d}\mathcal{H}}{\mathrm{d}t} = 0$. Expanding the time derivative of $\mathcal{H}$ as defined in \eqref{eq:conservation_laws} leads to
\begin{equation}\label{eq:enstrophy_conservation_time_derivative_expansion}
	\frac{\mathrm{d}\mathcal{H}}{\mathrm{d}t} = \int_{\Omega}\frac{\partial\boldsymbol{u}}{\partial t} \cdot \boldsymbol{\omega} + \int_{\Omega}\boldsymbol{u} \cdot \frac{\partial\boldsymbol{\omega}}{\partial t}\,.  
\end{equation}
If we now use the momentum equation in \eqref{eq:navier_stokes_continuous} and its curl, \eqref{eq:enstrophy_conservation_time_derivative_expansion} may be rewritten as
\[
\begin{aligned}
	\frac{\mathrm{d}\mathcal{H}}{\mathrm{d}t} 
	&= -\int_{\Omega}\left( \boldsymbol{\omega}\times\boldsymbol{u}\right)  \cdot \boldsymbol{\omega} - \int_{\Omega} \boldsymbol{u} \cdot \nabla\times\left(\boldsymbol{\omega}\times\boldsymbol{u}\right) 
	+ \int_{\Omega}\nabla\times\nabla P\cdot\bu  + \int_{\Omega}\nabla P \cdot\bw\\
	&= -\int_{\Omega}\left( \boldsymbol{\omega}\times\boldsymbol{u}\right)  \cdot \boldsymbol{\omega} - \int_{\Omega} \boldsymbol{\omega} \cdot \left(\boldsymbol{\omega}\times\boldsymbol{u}\right) + \int_{\Omega}\nabla\times\nabla P\cdot\bu  - \int_{\Omega} P \cdot\nabla\cdot\nabla\times\bu = 0\,,
\end{aligned}
\]
where we have used (i) the definition of vorticity $\boldsymbol{\omega} := \nabla\times\boldsymbol{u}$, (ii) integration by parts on the second and fourth terms in the right side of the first identity, (iii) the vector calculus relation \eqref{Eq: perpendicular inner product}, and (iv) the identities $ \nabla\times\nabla\left(\cdot\right)\equiv \MOD{\boldsymbol{0}} $ and $ \nabla\cdot\nabla\times\left(\cdot\right)\equiv 0 $.

\REVO{These conservation laws for kinetic energy (in 2D and 3D), enstrophy (in 2D), and helicity (in 3D), are the expression of a more general structure underlying the incompressible Euler equations: the Hamiltonian structure, \cite{morrison1980noncanonical, morrison1982poisson, Olver1982, Salmon1988, Morrison1998, chandre2012on, morrison2020lagrangian}. A system of partial differential equations (PDEs) is Hamiltonian if it can be cast in the general form, see for example \cite{Olver1982, abraham_diff_geom},
\[
	\frac{\partial \boldsymbol{y}}{\partial t} = \mathcal{S}\frac{\delta H(\boldsymbol{y})}{\delta \boldsymbol{y}}\,,
\]
where $\mathcal{S}$ is a skew-adjoint operator, such that the induced bilinear form must also be a derivation and satisfy the Jacobi-identity, and $H$ is the Hamiltonian functional.}

The system of equations \eqref{eq:navier_stokes_continuous}, in the inviscid limit, is not in Hamiltonian form but may be rewritten in this form if pressure is eliminated and the Hamiltonian functional is set to the kinetic energy $\mathcal{K}$. This can be achieved by either: (i) restricting the momentum equation to divergence free velocity fields (e.g., by making use of the stream function (in 2D) or the stream vector field (in 3D) $\boldsymbol{\psi}$ such that $\boldsymbol{u} = \nabla\times\boldsymbol{\psi}$) \cite{Arnold1969}, or (ii) taking the curl of the momentum equation (transforming it into the vorticity equation) \cite{Olver1982}. 

Noether's theorem establishes a connection between conservation laws of a Hamiltonian system and its underlying symmetries, \cite{Olver1980, Holm2011, Fecko2006, goldsteinClassicalMechanics}, thus highlighting the strong connection between the (geometric) structure of a system of PDEs and its dynamics. For example, spatial translation symmetry gives rise to conservation of linear momentum, and temporal translation symmetry results in energy conservation.

Helicity, on the other hand, is a more subtle quantity. Helicity as introduced in \eqref{eq:conservation_laws} is a particular case of the general concept of helicity of a divergence-free vector field, $\boldsymbol{v}$, tangent to the boundary $\partial\Omega$ of a simply connected domain $\Omega\subset\mathbb{R}^{3}$, see for example \cite{Holm2011, Arnold1992},
\[
	\mathcal{H}(\nabla\times\boldsymbol{v}) := \int_{\Omega} \boldsymbol{v}\cdot\nabla\times\boldsymbol{v}\,,
\]
which measures the average linking of its field lines.

The 19th century works of Helmholtz \cite{Helmholtz1867} and Kelvin \cite{Kelvin2011} contain the seminal ideas for the modern concept of helicity, \cite{Moffatt2008}. A renewed interest in these ideas appeared only later in the mid 20th century, first in the context of \MOD{magnetohydrodynamics (MHD)}, \cite{Woltjer1958}, and then for hydrodynamics, \cite{Moreau1960, Moffatt1969}. Moreau, \cite{Moreau1960}, discovered the law of conservation of helicity, and the term helicity appeared first in the work by Moffatt, \cite{Moffatt1969}, where the topological nature of this quantity was highlighted. For a detailed historical discussion of helicity see the very informative works by Moffatt \cite{Moffatt1992, Moffatt2014}.

It is possible to show, see for example \cite{Arnold1992}, that helicity of any divergence-free vector field is preserved under the action of any volume preserving diffeomorphism. This property shows that helicity is not a dynamical invariant but a topological invariant, since its conservation is independent of the specific diffeomorphism, \cite{Holm2011}. In fact, helicity is associated to the nontrivial kernel of the operator $\mathcal{S}$, and is a Casimir for the Hamiltonian formulation of the inviscid Navier-Stokes equations, \cite{Holm2011}.	\REVT{In the same way, in 2D, enstrophy is also a Casimir of the inviscid Navier-Stokes equations (as are all integral powers of vorticity).}

This very brief digression into Hamiltonian formalism intends to show the connection between the physical properties of a system of PDEs and its underlying geometrical structure. Invariants are not mere incidental features of the dynamics of a system, they are expressions of the underlying structure of the equations. 

Helicity plays \REVO{an} important role in the generation and evolution of turbulence, \cite{capuano2018effects, vallefuoco2019discrete, Yan2020}. The joint cascade of energy and helicity, \cite{Brissaud1973}, is an active field of research, \cite{Biferale2013, Chen2003a, Kessar2015, Sahoo2015, Alexakis2018}. Particularly important is the interaction between the two and how helicity impacts the energy cascade and, therefore, turbulence, \cite{Moffatt1969, capuano2018effects, Yan2020, Kessar2015, Ditlevsen2001a, Chen2003c, Chen2003b}. This complex interaction between the energy cascade and the helicity cascade and especially the suppressive role of helicity motivates the focus on the development of discretization schemes that, besides conserving energy, conserve helicity. In the same way as energy conserving schemes \REVT{have shown to substantially} contribute to a higher fidelity in simulations, see for example \cite{Morinishi1998, Capuano2015c, Duponcheel2008, Mullen2009, Pavlov2011, perot43discrete, Arakawa1966}, due to the connection between the cascades of energy and helicity, helicity conserving schemes should also present a positive impact towards improving the simulation accuracy.

\subsection{Overview of structure\MOD{-}preserving methods for fluid flows}\label{section:introduction:overview_of_structure_preserving}
As highlighted above, the solutions to systems of PDEs (of which the Navier-Stokes are a particular example) satisfy strong constraints \cite{Christiansen2011,Tadmor2012}. These constraints reflect the underlying \emph{mathematical structure} of the equations (e.g., Hamiltonian structure, Poisson structure, de Rham sequence). These fundamental \emph{mathematical structures} have long played an essential role in modern physics and pure mathematics. Owing to the fundamental nature of these structures and their impact on the dynamics of the systems under study, in recent years there has been an increasing interest in the various aspects of structure preservation at the discrete level \cite{Christiansen2011, ArnoldBook2006a, Koren2014}. This interest is rooted in three important points. First, there are well known connections between discrete structure preservation and \emph{standard properties} of numerical methods \cite{Christiansen2011, Hairer2006, arnold2006finite}. Second, \emph{standard properties} only guarantee physical fidelity in the limit of fully (at least highly) resolved discretizations. Reaching this limit requires infeasible computational resources (e.g., \cite{Orszag1970}). In contrast, structure preserving discretizations, by construction, generate solutions that satisfy the underlying physics even in highly under-resolved simulations. This is extremely relevant since most (if not all) simulations are inherently under-resolved. Third, physics preservation is fundamental when coupling systems in multiphysics problems \cite{Dowell2001}. 

The underlying principle behind structure preserving discretizations is to construct discrete approximations that retain as much as possible the structure of the original system of PDEs. A departure from this principle introduces spurious unphysical modes that pollute the physics of the system being modeled \cite{Hairer2006, arnold2010finite, bochev2003discourse}. For example, as seen before, turbulence plays a fundamental role in the dynamics of the flow. A correct representation of the turbulent dynamics of a fluid is paramount in order to achieve accurate simulations. For this reason, if a numerical discretization introduces spurious unphysical energy dissipation into the system, it will fail to accurately capture the energy cascade and consequently the turbulent dynamics, \cite{Rebholz2007b, Verstappen2003, Rebholz2007}. The main focus of structure preserving discretizations for flow problems has been on energy conservation, e.g. \cite{Morinishi1998, Capuano2015c, Duponcheel2008, Verstappen2003}. As noted in the previous section, there is a growing knowledge on the role played by helicity and its impact on the energy cascade. For this reason, more recently,  helicity conservation at the discrete level has been addressed in the literature, see for example \REVT{\cite{Capuano2018, Rebholz2007, hu2021helicity}}.

Most standard structure preserving discretizations can be seen as variations of staggered grid methods which date back to the pioneering works of Harlow and Welch \cite{Harlow1965}, and Arakawa and colleagues \cite{Arakawa1977, Mesinger1976}. These methods employ a discretization that distributes the different physical quantities (pressure, velocity, vorticity, etc) at different locations in the mesh (vertices, faces, cell centres). It can be shown that, by doing so, important conservation properties can be maintained. Since then, much work has been produced and a rich variety of different flavours of structure preserving discretizations have been presented: finite differences/finite volumes \cite{HymanShashkovSteinberg97,BrezziBuffaLipnikov2009,HymanShashkovSteinberg2002,RobidouxAdjointGradients1996, Perot2000}, discrete exterior calculus (DEC) \cite{desbrun2005discrete}, finite element exterior calculus (FEEC) \cite{arnold2006finite, Bossavit1998, hiptmair2001} and the works by the authors \cite{Palha2017a, kreeft::stokes,Lee2018b, Lee2019a, DeDiego2019a}. 

More recently, another approach develops a discretization of the physical field laws based on discrete variational principles. This approach has been used in the past to construct variational integrators for Lagrangian systems, e.g. \cite{Kouranbaeva2000,Marsden2003}. These ideas have been extended to magneto-hydrodynamics \cite{Kraus2015, Kraus2018, Kaltsas2019}, incompressible flows \cite{Gawlik2020}, and geophysical flow \cite{Brecht2019a, Bauer2019}.

\subsection{Objective}\label{SEC: Objective}
In this work, extending the initial ideas introduced for the 2D case, see \cite{Palha2017a}, we combine (i) a particular choice for the formulation of the Navier-Stokes equations with (ii) a structure preserving discretization. Specifically, we will present two velocity evolution equations (dual-field) in \MOD{a} rotational form, discretized by the mimetic spectral element method (MSEM) \cite{kreeft::stokes, palha2014physics, Kreeft2011}. 

\REVT{This formulation attempts to address the dual character of the velocity field in the incompressible Navier-Stokes equations. 
This dual character implies that it is natural to look for a solution for the velocity field in  $ H(\mathrm{div}; \Omega) \cap H(\mathrm{curl}; \Omega)$.
At the continuous level this is easily achievable, but that is not true at the discrete level since the space $H(\mathrm{div}; \Omega) \cap H(\mathrm{curl}; \Omega)$ is hard to discretize}.  The use of two velocity field evolution equations enables the representation of this dual character. It is shown that in this way the resulting discretization conserves mass, kinetic energy, and helicity in 3D. 

The vorticity fields \MOD{in the rotational form of the nonlinear convective term, see \eqref{eq: forms for nonlinear convective term rotational},} serve as a means of exchanging information between the two evolution equations. Additionally, this leads to a leap-frog like scheme that handles the nonlinear rotational term by staggering in time the velocity and vorticity such that the resulting discrete algebraic systems are \MOD{linearized} and decoupled.

Overall, the objective of this novel approach is the construction of a discretization which conserves mass, kinetic energy and helicity for the incompressible Navier-Stokes equations in the absence of dissipative terms and predicts the proper decay rate of kinetic energy and helicity based on the global enstrophy and an integral quantity of vorticity, respectively.

\subsection{Outline of paper}
The outline of the paper is as follows: In Section \ref{Sec: scheme}, we introduce a dual-field mixed weak formulation and prove that it preserves the desired conservation properties. In Section \ref{Sec: temporal discretization}, a conservative staggered temporal discretization scheme is applied to the formulation, which is followed by a mimetic spatial discretization in Section \ref{Sec: spatial discretization}. Numerical results that support the method are presented in Section \ref{Sec: Numerical results}. Finally, a summary is given and potential future work is listed in Section \ref{Sec: Conclusions}.

\section{A mass-, kinetic energy- and helicity-conserving formulation}\label{Sec: scheme}
In this section, we propose a new conservative formulation for the Navier-Stokes equations in periodic domains. As we will only consider periodic domains in this paper, from now on, $ \Omega $ represents a 3D periodic domain. The function spaces are the classic Hilbert spaces which form an exact complex, namely, the well-known de Rham (or Hilbert) complex \cite{Palha2017a, arnold2006finite, arnold2010finite, bochev2003discourse}:
\begin{equation}\label{Eq: de Rham complex}
	\mathbb{R}\hookrightarrow H^1(\Omega)\stackrel{\nabla}{\longrightarrow}H(\mathrm{curl};\Omega)\stackrel{\nabla\times}{\longrightarrow}H(\mathrm{div};\Omega)\stackrel{\nabla\cdot}{\longrightarrow}L^2(\Omega)\rightarrow 0\,.
\end{equation} 
This complex plays a fundamental role in the proofs and analysis of the presented work.

\subsection{The rotational form of the incompressible Navier-Stokes equations} \label{Sub: rotational form}
If in \eqref{eq:navier_stokes_continuous} we use the rotational (or Lamb) form for the nonlinear convective term, \eqref{eq: forms for nonlinear convective term rotational}, and use the representaion $ \mathcal{D}(\boldsymbol{u}) = - \nabla\times\bw $ for the linear dissipative term, \eqref{Eq: representations for the dissipative term}, we obtain the rotational form of the incompressible Navier-Stokes equations,  \vspace{-3ex}

\begin{subequations}\label{Eq: NS rotational formulation}
	\begin{align}
		&\dfrac{\partial\boldsymbol{u}}{\partial t} + \boldsymbol{\omega}\times\boldsymbol{u} + \eR \nabla\times \boldsymbol{\omega} + \nabla P = \boldsymbol{f}\,,\label{Eq: NS rotational formulation a}\\
		&\boldsymbol{\omega} = \nabla\times \boldsymbol{u}\,,\\
		& \nabla\cdot\boldsymbol{u} = 0\,.  \label{Eq: NS rotational formulation c}
	\end{align}
\end{subequations}
We have proven that, in 3D and in the inviscid limit $ \left(\Rn\to\infty \right) $, these equations preserve total kinetic energy and total helicity over time for the case of no external body force, $ \boldsymbol{f}=\boldsymbol{0} $, in Section~\ref{section:introduction:relevance_structure_preserving}. For non-zero conservative external body force, $ \boldsymbol{f}=\nabla\varphi\neq\boldsymbol{0} $, we can include it by replacing the total pressure by an extended total pressure
\begin{equation}\label{Eq: extended total pressure}
	P^{\prime}:= P - \varphi\,.
\end{equation}
All analysis and proofs remain valid.
Without loss of generality, in this paper we will only use zero external body force for the analysis and proofs.

When the flow is viscous, $ \Rn < \infty $, the viscosity dissipates  kinetic energy of the incompressible Navier-Stokes equations at rate
\begin{equation} \label{Eq: continuous K disspation}
	\dfrac{\mathrm{d} \mathcal{K}}{\mathrm{d} t} = - \dfrac{2}{\Rn}  \mathcal{E}\,,
\end{equation}
while it dissipates or generates helicity at rate 
\begin{equation} \label{Eq: continuous H disspation}
	\dfrac{\mathrm{d} \mathcal{H}}{\mathrm{d} t} = - \dfrac{2}{\Rn} \left\langle \bw, \nabla\times\bw \right\rangle_{\Omega}\,,
\end{equation}
where $ \left\langle \cdot, \cdot\right\rangle _{\Omega} $ denotes the inner product, i.e.,  
\[
\left\langle \boldsymbol{a},\boldsymbol{b}\right\rangle _{\Omega}=\int_{\Omega}\boldsymbol{a}\cdot\boldsymbol{b} \quad \text{and}\quad \left\langle c, d\right\rangle _{\Omega} = \int_{\Omega}cd\,,
\]
if $ \boldsymbol{a},\boldsymbol{b} $ are vectors and $ c,d $ are scalars. The viscosity always dissipates kinetic energy because the total enstrophy cannot be negative, $ \mathcal{E}\ge 0 $, see the definition of the total enstrophy in \eqref{eq:conservation_laws}. It either dissipates or generates helicity because the term $ \left\langle \bw, \nabla\times\bw \right\rangle_{\Omega} $ generally can be either positive or negative (or zero). This means the dissipation rate of helicity can be negative.

\subsection{A conservative dual-field mixed weak  formulation}
We propose the following dual-field mixed weak formulation \MOD{for the rotational form of the incompressible Navier-Stokes equations}: Given $ \boldsymbol{f} \in \left[ L^{2}(\Omega) \right] ^{3} $, seek $ \left( \boldsymbol{u}_{1}, \bw_{2},  P_{0}\right)  \in H(\mathrm{curl}; \Omega)\times H(\mathrm{div}; \Omega)\times H^1(\Omega) $ and $ \left( \bu_{2},\bw_{1},P_{3}\right)  \in H(\mathrm{div}; \Omega) \times H(\mathrm{curl}; \Omega)\times  L^2(\Omega) $ such that, \vspace{-3ex}

\begin{subequations}\label{Eq: WF}
	\begin{align}
	&\left\langle \dfrac{\partial \bu_{1}}{\partial t}, \be_{1}\right\rangle _{\Omega}  + \left\langle\bw_{1}\times \bu_{1} , \be_{1}\right\rangle _{\Omega} + \eR\left\langle \bw_{2}, \nabla\times\be_{1}\right\rangle _{\Omega} + \left\langle\nabla P_{0}, \be_{1}\right\rangle _{\Omega} = \left\langle \boldsymbol{f}, \be_{1}\right\rangle _{\Omega} &&\forall\be_{1}\in H(\mathrm{curl};\Omega)\,,\label{Eq: WF d}\\
	&\left\langle \nabla\times \bu_{1}, \be_{2}\right\rangle _{\Omega}  - \left\langle \bw_{2},\be_{2}\right\rangle _{\Omega}  = 0 && \forall\be_{2}\in H(\mathrm{div};\Omega)\,, \label{Eq: WF e}\\
	&\left\langle \bu_{1}, \nabla\epsilon_{0} \right\rangle _{\Omega} = 0 &&\forall\epsilon_{0}\in H^{1}(\Omega) \,,\label{Eq: WF f}\\
	&\left\langle  \dfrac{\partial \bu_{2}}{\partial t}, \be_{2}\right\rangle_{\Omega}   + \left\langle \bw_{2}\times \bu_{2} , \be_{2}\right\rangle _{\Omega} + \eR\left\langle\nabla\times\bw_{1}, \be_{2}\right\rangle _{\Omega}- \left\langle  P_{3}, \nabla\cdot\be_{2}\right\rangle _{\Omega}  = \left\langle \boldsymbol{f}, \be_{2}\right\rangle _{\Omega} &&\forall\be_{2}\in H(\mathrm{div};\Omega)\,, \label{Eq: WF a}\\
	&\left\langle  \bu_{2}, \nabla\times\be_{1}\right\rangle _{\Omega}-\left\langle \bw_{1} , \be_{1}\right\rangle _{\Omega}=0 &&\forall\be_{1}\in H(\mathrm{curl};\Omega)\,, \label{Eq: WF b}\\
	&\left\langle\nabla\cdot \bu_{2}, \epsilon_{3} \right\rangle _{\Omega} = 0 &&\forall\epsilon_{3}\in L^{2}(\Omega)\,. \label{Eq: WF c}
	\end{align}
\end{subequations}
\begin{remark}\label{remark L2 integrable}\REVO{
	In this formulation, the terms $ \left( \boldsymbol{\omega}_{i}\times \boldsymbol{u}_{i} \right) \cdot \boldsymbol{\epsilon}_{i} $ are not known to be $ L^2 $-integrable for the vector fields that belong to the infinite dimensional function spaces $ H(\mathrm{curl}; \Omega) $ ($ i=1 $) and $ H(\mathrm{div}; \Omega) $ ($ i=2 $).  Showing this integrability requires proving additional regularity of the velocity and vorticity variables, which we currently are unable to do.  However, in the finite dimensional case, the  known regularity is sufficient, see Section 4.  Thus, despite the potential mathematical issue, we still write this formulation above for its clear interpretation and to motivate the discrete scheme.
}\end{remark}
The formulation \eqref{Eq: WF} is called \emph{dual-field} because it contains two evolution equations, \eqref{Eq: WF d} and \eqref{Eq: WF a}, and, for each variable, dual representations of its solution are sought: For velocity, we seek $ \left( \boldsymbol{u}_{1}, \boldsymbol{u}_{2}\right) \in H(\mathrm{curl};\Omega) \times H(\mathrm{div};\Omega) $, for vorticity, we seek $ \left( \boldsymbol{\omega}_{2}, \boldsymbol{\omega}_{1}\right) \in H(\mathrm{div};\Omega) \times H(\mathrm{curl};\Omega) $, and, for total pressure, we seek $ \left( P_{0}, P_{3}\right) \in H^1(\Omega) \times L^2(\Omega) $.
If all variables are sufficiently smooth, integration by parts will show that either $ \left( \boldsymbol{u}_{1}, \bw_{2},  P_{0}\right) $ or $ \left( \bu_{2},\bw_{1},P_{3}\right) $ solves the Navier-Stokes equations in rotational form, \eqref{Eq: NS rotational formulation}. 
Note that the de Rham complex \eqref{Eq: de Rham complex} and the constraint \eqref{Eq: WF e} ensure
\begin{equation}\label{Eq: curl exact}
	\bw_{2}=\nabla\times \bu_{1}\,. 
\end{equation}
Therefore, in practice, $ \bw_{2} $ may be dropped from \eqref{Eq: WF} if we replace it by $ \nabla\times \bu_{1} $. We leave in $ \bw_{2} $ above to maintain the clearness of the formulation. 


\subsection{Properties of the formulation} \label{Sub: continuous conservation proof}
We now show that the proposed dual-field formulation \eqref{Eq: WF} conserves (i) the mass in terms of $ \bu_{2} $ and, in the case of conservative external body force and zero viscosity, (ii) the kinetic energy in the formats 
\[\mathcal{K}_{1} =\frac{1}{2} \left\langle \bu_{1},\bu_{1}\right\rangle _{\Omega}\quad\text{and}\quad \mathcal{K}_{2} = \frac{1}{2}\left\langle \bu_{2},\bu_{2}\right\rangle _{\Omega} \,, \] 
and (iii) the helicity in the formats 
\[ \mathcal{H}_{1} = \int_{\Omega}\bu_{1}\cdot\bw_{1} = \left\langle \bu_{1},\bw_{1}\right\rangle _{\Omega} \quad\text{and}\quad \mathcal{H}_{2}= \int_{\Omega}\bu_{2}\cdot\bw_{2} = \left\langle \bu_{2},\bw_{2}\right\rangle _{\Omega} \,.\]
We will also analyze the dissipation rate of kinetic energy and helicity in the viscous case for the proposed formulation.

Note that, in this subsection, everything is still at the continuous level. The purpose is to show that the proposed weak formulation possesses the same properties as the strong formulation does.


\subsubsection{Mass conservation}\label{Subsub: continuous mass conservation}
For the mass conservation, since we have restricted $ \bu_{2}$ to space $ H(\mathrm{div};\Omega) $, the de Rham complex \eqref{Eq: de Rham complex} and the constraint \eqref{Eq: WF c} ensure that the relation
\[
 H(\mathrm{div};\Omega)\ni\bu_{2}\stackrel{\nabla\cdot}{\longrightarrow}0\in  L^2(\Omega)
\] is strongly satisfied; no integration by parts is required. Therefore, the mass conservation is satisfied for velocity $ \bu_{2}$. Such an approach is widely used to construct mass conserving discretizations. While for $ \bu_{1} \in H(\mathrm{curl};\Omega) $, the mass conservation is only weakly satisfied, see \eqref{Eq: WF f}.

\subsubsection{Time rate of change of kinetic energy} \label{Subsub: continuous kinetic energy conservation}
In the inviscid limit $ (\Rn\to\infty) $ and when $ \boldsymbol{f}=\boldsymbol{0} $, the kinetic energy conservation is equivalent to
\[
\dfrac{\mathrm{d} \mathcal{K}_{1}}{\mathrm{d} t} = \left\langle \dfrac{\partial \bu_{1}}{\partial t}, \bu_{1}\right\rangle _{\Omega} = 0\quad \text{and}\quad
\dfrac{\mathrm{d} \mathcal{K}_{2}}{\mathrm{d} t} = \left\langle \dfrac{\partial \bu_{2}}{\partial t}, \bu_{2}\right\rangle _{\Omega} = 0\,.
\]
Because \eqref{Eq: WF d} is valid for all $ \be_{1}\in  H(\mathrm{curl};\Omega)  $, we can select $ \be_{1} $ to be $ \bu_{1}\in H(\mathrm{curl};\Omega) $. As a result, we get
\[
\left\langle \dfrac{\partial \bu_{1}}{\partial t}, \bu_{1}\right\rangle _{\Omega}  + {\left\langle\bw_{1}\times \bu_{1} , \bu_{1}\right\rangle _{\Omega}} + \left\langle\bu_{1}, \nabla P_{0} \right\rangle _{\Omega} =\left\langle \dfrac{\partial \bu_{1}}{\partial t}, \bu_{1}\right\rangle _{\Omega} = 0\,.
\]
The second term vanishes because of \eqref{Eq: perpendicular inner product}.
Meanwhile, from \eqref{Eq: WF f}, we know that $ \left\langle\bu_{1}, \nabla P_{0}\right\rangle _{\Omega} = 0 $ because $ P_{0}\in H^{1}(\Omega) $. Therefore, the third term also vanishes, which accomplishes the proof of kinetic energy conservation for $ \mathcal{K}_{1} $. Similarly, by selecting $ \be_{2} $ of \eqref{Eq: WF a} to be $ \bu_{2} $, we can get
\[
\left\langle  \dfrac{\partial \bu_{2}}{\partial t}, \bu_{2}\right\rangle_{\Omega}   + {\left\langle \bw_{2}\times \bu_{2} , \bu_{2}\right\rangle _{\Omega}} - {\left\langle  P_{3}, \nabla\cdot\bu_{2}\right\rangle _{\Omega}}  = \left\langle  \dfrac{\partial \bu_{2}}{\partial t}, \bu_{2}\right\rangle_{\Omega} =0\,,
\]
where the second and third terms vanish \MOD{because} \eqref{Eq: perpendicular inner product} and \eqref{Eq: WF c}, respectively. Thus we can conclude that $ \mathcal{K}_{2}$  is also preserved over time.

In the viscous case, $ \Rn<\infty $, if we repeat the above analysis, the viscous terms will remain. We will eventually obtain the following kinetic energy dissipation rates,
\begin{equation}\label{Eq: boundedness of K1}
\dfrac{\mathrm{d} \mathcal{K}_{1}}{\mathrm{d} t} = \left\langle \dfrac{\partial \bu_{1}}{\partial t}, \bu_{1}\right\rangle _{\Omega} = - \eR\left\langle\bw_{2}, \nabla\times\bu_{1}\right\rangle _{\Omega} = - \eR\left\langle\bw_{2}, \bw_{2}\right\rangle _{\Omega} = -\dfrac{2}{\Rn}\mathcal{E}_{2} \leq 0\,,
\end{equation}

\begin{equation}\label{Eq: boundedness of K2}
\dfrac{\mathrm{d} \mathcal{K}_{2}}{\mathrm{d} t} = \left\langle \dfrac{\partial \bu_{2}}{\partial t}, \bu_{2}\right\rangle _{\Omega} = - \eR\left\langle\nabla\times\bw_{1}, \bu_{2}\right\rangle _{\Omega}
\stackrel{\eqref{Eq: WF b}}{=} -\eR\left\langle\bw_{1}, \bw_{1}\right\rangle _{\Omega} = -\dfrac{2}{\Rn}\mathcal{E}_{1}\leq 0\,,
\end{equation}
where the total enstrophy $ \mathcal{E}_{1} $ and $ \mathcal{E}_{2} $ are defined as
\[\mathcal{E}_{1}=\dfrac{1}{2}\left\langle \bw_{1},\bw_1\right\rangle_{\Omega} \quad \text{and}\quad \mathcal{E}_{2}=\dfrac{1}{2}\left\langle \bw_{2},\bw_2\right\rangle_{\Omega}\,.\]
This is in agreement with the kinetic energy dissipation rate of the strong formulation, see \eqref{Eq: continuous K disspation}.

\subsubsection{Time rate of change of helicity}\label{Subsub: continuous helicity conservation}
If $ \Rn\to\infty $ and $ \boldsymbol{f}=\boldsymbol{0} $, the helicity conservation is equivalent to
\[
\dfrac{\mathrm{d} \mathcal{H}_{1}}{\mathrm{d} t} = \dfrac{\mathrm{d} }{\mathrm{d} t}\left\langle \bu_{1},\bw_{1}\right\rangle _{\Omega} = \left\langle \dfrac{\partial \bu_{1}}{\partial t}, \bw_{1}\right\rangle _{\Omega} + \left\langle \bu_{1}, \dfrac{\partial \bw_{1}}{\partial t}\right\rangle _{\Omega} = 0\,,
\]

\[
\dfrac{\mathrm{d} \mathcal{H}_{2}}{\mathrm{d} t} = \dfrac{\mathrm{d} }{\mathrm{d} t}\left\langle \bu_{2},\bw_{2}\right\rangle _{\Omega} = \left\langle \dfrac{\partial \bu_{2}}{\partial t}, \bw_{2}\right\rangle _{\Omega} + \left\langle \bu_{2}, \dfrac{\partial \bw_{2}}{\partial t}\right\rangle _{\Omega} = 0\,.
\]
Replacing $ \be_{1} $ in \eqref{Eq: WF d} by $ \bw_{1}\in H(\mathrm{curl};\Omega) $ leads to
\begin{equation}\label{Eq: helicity leg 1}
	\left\langle \dfrac{\partial \bu_{1}}{\partial t}, \bw_{1}\right\rangle _{\Omega}  + {\left\langle\bw_{1}\times \bu_{1} , \bw_{1}\right\rangle _{\Omega}} + {\left\langle\bw_{1},\nabla P_{0} \right\rangle _{\Omega}} = \left\langle \dfrac{\partial \bu_{1}}{\partial t}, \bw_{1}\right\rangle _{\Omega}=0\,.
\end{equation}
The second term vanishes because of \eqref{Eq: perpendicular inner product}. Meanwhile, we have \eqref{Eq: WF b} saying
\begin{equation}\label{Eq: helicity conservation relation 0}
-\left\langle \bu_{2}, \nabla\times\be_{1}\right\rangle _{\Omega}+\left\langle \bw_{1} , \be_{1}\right\rangle _{\Omega}=0\quad \forall\be_{1}\in H(\mathrm{curl};\Omega)\,.
\end{equation}
And because $ P_{0}\in H^{1}(\Omega) $, we have $ \nabla P_{0}\in  H(\mathrm{curl};\Omega) $ (in particular, $ \nabla P_{0} $ is in the null space of $ H(\mathrm{curl};\Omega) $ with respect to $ \nabla\times $ ). Thus we can replace $ \be_{1} $ in \eqref{Eq: helicity conservation relation 0} by $ \nabla P_{0} $ and get  
\begin{equation}\label{Eq: H1 importance}
-\left\langle  \bu_{2}, \nabla\times\nabla P_{0}\right\rangle _{\Omega}+\left\langle \bw_{1} , \nabla P_{0}\right\rangle _{\Omega}=0 \,.
\end{equation}
This implies $ \left\langle \bw_{1} , \nabla P_{0}\right\rangle _{\Omega}=0 $ because $ \nabla\times\nabla(\cdot) \equiv 0$ showing that the third term of \eqref{Eq: helicity leg 1} vanishes.

If we take the time derivative of \eqref{Eq: helicity conservation relation 0}, we have
\begin{equation}\label{Eq: helicity conservation relation 1}
\left\langle\dfrac{ \partial \bu_{2}}{\partial t}, \nabla\times\be_{1}\right\rangle _{\Omega} = \left\langle \dfrac{\partial \bw_{1}}{\partial t} , \be_{1}\right\rangle _{\Omega}\quad \forall\be_{1}\in H(\mathrm{curl};\Omega)\,.
\end{equation}
In addition, we know that, \eqref{Eq: WF a},
\begin{equation}\label{Eq: equation 100}
	\left\langle \dfrac{\partial \bu_{2}}{\partial t}, \be_{2}\right\rangle _{\Omega}  + \left\langle\bw_{2}\times \bu_{2} , \be_{2}\right\rangle _{\Omega} - \left\langle  P_{3}, \nabla\cdot\be_{2}\right\rangle _{\Omega}  = 0\quad \forall \be_{2}\in H(\mathrm{div};\Omega)\,.
\end{equation}
Therefore, given any $\be_{1} \in H(\mathrm{curl};\Omega) $, \eqref{Eq: equation 100} must hold for $ \nabla\times \be_{1}\in H(\mathrm{div};\Omega) $, i.e.,
\begin{equation}\label{Eq: helicity conservation relation 2}
\left\langle \dfrac{\partial \bu_{2}}{\partial t}, \nabla\times\be_{1}\right\rangle _{\Omega}  + \left\langle\bw_{2}\times \bu_{2} , \nabla\times\be_{1}\right\rangle _{\Omega} - \left\langle  P_{3}, \nabla\cdot\nabla\times\be_{1}\right\rangle _{\Omega}  = 0\quad \forall\be_{1}\in H(\mathrm{curl};\Omega)\,.
\end{equation}
If we insert \eqref{Eq: helicity conservation relation 1} into \eqref{Eq: helicity conservation relation 2}, we obtain
\[
\left\langle \dfrac{\partial \bw_{1}}{\partial t} , \be_{1}\right\rangle _{\Omega}  + \left\langle\bw_{2}\times \bu_{2} , \nabla\times\be_{1}\right\rangle _{\Omega} - \left\langle  P_{3}, \nabla\cdot\nabla\times\be_{1}\right\rangle _{\Omega}  = 0\quad \forall\be_{1}\in H(\mathrm{curl};\Omega)\,.
\]
Because $ \bu_{1}\in H(\mathrm{curl};\Omega) $, we now replace $ \be_{1} $ in above equation with $ \bu_{1} $ and obtain
\begin{equation}\label{Eq: helicity leg 2}
\left\langle \dfrac{\partial \bw_{1}}{\partial t} , \bu_{1}\right\rangle _{\Omega}  + {\left\langle\bw_{2}\times \bu_{2} , \nabla\times\bu_{1}\right\rangle _{\Omega}} - {\left\langle  P_{3}, \nabla\cdot\nabla\times\bu_{1}\right\rangle _{\Omega}}  = \left\langle \dfrac{\partial \bw_{1}}{\partial t} , \bu_{1}\right\rangle _{\Omega}  = 0\,.
\end{equation}
Since $ \bw_{2}=\nabla\times\bu_{1} $, see \eqref{Eq: curl exact}, is exactly satisfied, the second term of \eqref{Eq: helicity leg 2} vanishes due to \eqref{Eq: perpendicular inner product}, and the third term is zero because $ \nabla\cdot\nabla\times(\cdot)\equiv 0 $. Overall, \eqref{Eq: helicity leg 1} and \eqref{Eq: helicity leg 2} together prove that helicity $ \mathcal{H}_{1} $ is preserved over time. 

We now reuse \eqref{Eq: WF b} and select $ \be_{1} $ to be $ \bu_{1}\in H(\mathrm{curl};\Omega) $. As a result, we get
\[ -\left\langle  \bu_{2}, \nabla\times\bu_{1}\right\rangle _{\Omega}+\left\langle \bw_{1} , \bu_{1}\right\rangle _{\Omega}=0\,,\]
which implies 
\[
\mathcal{H}_{2} = \left\langle \bu_{2},\bw_{2}\right\rangle _{\Omega} = \left\langle  \bu_{2}, \nabla\times\bu_{1}\right\rangle _{\Omega} = \left\langle \bw_{1} , \bu_{1}\right\rangle _{\Omega} = \mathcal{H}_{1}\,.
\]
Thus both $ \mathcal{H}_{1} $ and $ \mathcal{H}_{2} $ are preserved over time.

In the viscous case, $ \Rn<\infty  $, if we repeat above analysis, the viscous contribution will not cancel and we will obtain the following helicity dissipation rate,
\[\dfrac{\partial \mathcal{H}_{1}}{\partial t} = \dfrac{\partial \mathcal{H}_{2}}{\partial t} = -\dfrac{2}{\Rn}\left\langle \bw_{2},\nabla\times\bw_{1} \right\rangle _{\Omega}\,,\]
which is consistent with that of the strong formulation, see \eqref{Eq: continuous H disspation}.

\section{Temporal discretization}\label{Sec: temporal discretization}
Inspired by a mass, energy, enstrophy and vorticity conserving (MEEVC) \cite{Palha2017a} scheme for the 2D incompressible Navier-Stokes equations, we construct a staggered temporal discretization for the two evolution equations in the dual-field formulation \eqref{Eq: WF}.
The MEEVC scheme, as well as the presented method, starts with a formulation of two evolution equations. The two evolution equations are discretized temporally at two sequences of time steps respectively using a Gauss integrator. The two sequences of time steps are staggered such that the endpoints of time steps in one sequence are exactly the midpoints of time steps in the other sequence. Thus at each time step either discrete evolution equation can use the solution from the other one as known variable at the midpoint, see Fig.~\ref{fig: temporal scheme}. 

We use a lowest order Gauss integrator as the time integrator \cite{sanderse2013energy, hairer2006geometric,steinberg2016explicit}. For example, if we apply the integrator to an ordinary differential equation (ODE) of the form
\[\dfrac{\mathrm{d}f(t)}{\mathrm{d}t} = h\left( f(t), t\right) \]
at a time step from time instant $ t^{k-1} $ to time instant $ t^{k} $, we obtain 
\begin{equation}\label{Eq: gauss integrator}
\dfrac{f^{k}-f^{k-1}}{\varDelta t} = h\left(  f^{k-\frac{1}{2}}, t^{k-1}+\dfrac{\varDelta t}{2}\right)\,,
\end{equation}
where $ \varDelta t = t^{k}-t^{k-1} $, $ f^{k} = f(t^{k}) $. Additionally, we will use the midpoint rule, namely,
\begin{equation}\label{Eq: mpr}
f^{k-\frac{1}{2}} = f(t^{k-1}+\dfrac{\varDelta t}{2}) := \dfrac{f^{k}+f^{k-1}}{2}\,.
\end{equation}
We further introduce two time sequences, the \emph{integer time steps} and the \emph{half-integer time steps}. The integer time steps use time instants indicated with integer superscripts. For example, $ k $th $ (k=1,2,\cdots) $ integer time step (denoted by $ S_{k} $) is from $ t^{k-1}$ to $ t^{k} $. The half-integer time steps use time instants indicated with half-integer superscripts. For example, $ k $th $ (k=1,2,\cdots) $ half-integer time step (denoted by $ \hat{S}_{k} $) is from $ t^{k-\frac{1}{2}}$ to $ t^{k+\frac{1}{2}} $. These time steps satisfy
\[\varDelta t = t^{i} - t^{i-1} = t^{j+\frac{1}{2}} - t^{j-\frac{1}{2}}\ \mathrm{and}\ t^{k-\frac{1}{2}}=\dfrac{t^{k}+t^{k-1}}{2}\quad  \forall i,j,k = 1,2,\cdots\,.\]
In other words, we restrict ourselves to constant time intervals equal for both time sequences.

\subsection{Temporal discretizations at staggered time steps} 
We now apply the time integrator \eqref{Eq: gauss integrator} to evolution equations \eqref{Eq: WF a} and \eqref{Eq: WF d} at integer and half-integer time steps, respectively.

\subsubsection{Temporal discretization at integer time steps}
If we apply the time integrator \eqref{Eq: gauss integrator} to the evolution equation for $ \bu_{2} $ \eqref{Eq: WF a} at integer time steps, with the midpoint rule, see \eqref{Eq: mpr}, and constraints \eqref{Eq: WF b} and \eqref{Eq: WF c}, we can obtain a semi-discrete weak formulation at, for example, $ k $th integer time step $ S_{k} $: Given $ \left( \bw^{k-1}_{1}, \bu^{k-1}_{2},\boldsymbol{f}^{k-\frac{1}{2}}, \bw^{k-\frac{1}{2}}_{2}\right) \in H(\mathrm{curl}; \Omega)\times H(\mathrm{div}; \Omega)\times\left[L^2(\Omega)  \right] ^{3}  \times  H(\mathrm{div}; \Omega) $,  find $ \left( \bw^{k}_{1}, \bu^{k}_{2}, P^{k-\frac{1}{2}}_{3}\right) \in H(\mathrm{curl}; \Omega)\times H(\mathrm{div}; \Omega) \times L^2(\Omega)  $ such that \vspace{-2ex}

\begin{subequations}\label{Eq: TD1}
	\begin{align}
		&\left\langle \dfrac{\bu^{k}_{2}-\bu^{k-1}_{2}}{\varDelta t}, \be_{2}\right\rangle _{\Omega}  + \left\langle\bw^{k-\frac{1}{2}}_{2}\times \dfrac{\bu^{k}_{2}+\bu^{k-1}_{2}}{2} , \be_{2}\right\rangle _{\Omega}  \label{Eq: TD1 a}\\
		&\hspace{1.5cm}+ \eR\left\langle\nabla\times\dfrac{\bw^{k}_{1}+\bw^{k-1}_{1}}{2}, \be_{2}\right\rangle _{\Omega}- \left\langle  P^{k-\frac{1}{2}}_{3}, \nabla\cdot\be_{2}\right\rangle _{\Omega}  = \left\langle \boldsymbol{f}^{k-\frac{1}{2}}, \be_{2}\right\rangle _{\Omega} &&\forall\be_{2}\in H(\mathrm{div};\Omega)\,, \nonumber \\
		&\left\langle  \bu^{k}_{2}, \nabla\times\be_{1}\right\rangle _{\Omega}-\left\langle \bw^{k}_{1} , \be_{1}\right\rangle _{\Omega}=0 &&\forall\be_{1}\in H(\mathrm{curl};\Omega)\,,\label{Eq: TD1 b}\\
		&\left\langle\nabla\cdot \bu^{k}_{2}, \epsilon_{3} \right\rangle _{\Omega} = 0 &&\forall\epsilon_{3}\in L^{2}(\Omega)\,, \label{Eq: TD1 c}
	\end{align}
\end{subequations}
where $ \bw^{k-\frac{1}{2}}_{2} $ is borrowed from the other time sequence, in particular, is the solution of $ \bw_{2} $ at $ (k-1) $th half-integer time step and, therefore, is known. 

\subsubsection{Temporal discretization at half-integer time steps} 
Similarly, we apply the time integrator \eqref{Eq: gauss integrator} to the evolution equation for $ \bu_{1} $ \eqref{Eq: WF d} at half-integer time steps. With the midpoint rule, see \eqref{Eq: mpr}, and constraints \eqref{Eq: WF e} and \eqref{Eq: WF f}, we can get a second semi-discrete weak formulation at, for example, $ k $th half-integer time step $ \hat{S}_{k} $: Given $ \left( \bu^{k-\frac{1}{2}}_{1}, \bw^{k-\frac{1}{2}}_{2},\boldsymbol{f}^{k}, \bw^{k}_{1}\right) \in H(\mathrm{curl}; \Omega)\times H(\mathrm{div}; \Omega)\times\left[ L^2(\Omega)\right] ^{3} \times H(\mathrm{curl}; \Omega) $, seek $ \left( P^{k}_{0}, \bu^{k+\frac{1}{2}}_{1}, \bw^{k+\frac{1}{2}}_{2}\right) \in H^1(\Omega)\times H(\mathrm{curl}; \Omega)\times H(\mathrm{div}; \Omega)  $ such that \vspace{-2ex}

\begin{subequations}\label{Eq: TD2}
	\begin{align}
		&\left\langle \dfrac{\bu^{k+\frac{1}{2}}_{1}-\bu^{k-\frac{1}{2}}_{1}}{\varDelta t}, \be_{1}\right\rangle _{\Omega}  + \left\langle\bw^{k}_{1}\times \dfrac{\bu^{k+\frac{1}{2}}_{1}+\bu^{k-\frac{1}{2}}_{1}}{2} , \be_{1}\right\rangle _{\Omega} \label{Eq: TD2 a}\\
		&\hspace{1.5cm}+ \eR\left\langle \dfrac{\bw^{k+\frac{1}{2}}_{2}+\bw^{k-\frac{1}{2}}_{2}}{2}, \nabla\times\be_{1}\right\rangle _{\Omega} + \left\langle\nabla P^{k}_{0}, \be_{1}\right\rangle _{\Omega} = \left\langle \boldsymbol{f}^{k}, \be_{1}\right\rangle _{\Omega} &&\forall\be_{1}\in H(\mathrm{curl};\Omega)\,,\nonumber\\
		&\left\langle \nabla\times \bu^{k+\frac{1}{2}}_{1}, \be_{2}\right\rangle _{\Omega} - \left\langle \bw^{k+\frac{1}{2}}_{2},\be_{2}\right\rangle _{\Omega}  = 0 && \forall\be_{2}\in H(\mathrm{div};\Omega)\,,\label{Eq: TD2 b}\\
		&\left\langle \bu^{k+\frac{1}{2}}_{1}, \nabla\epsilon_{0} \right\rangle _{\Omega} = 0 &&\forall\epsilon_{0}\in H^{1}(\Omega)\,, \label{Eq: TD2 c}
	\end{align}
\end{subequations}
where $ \bw^{k}_{1} $ is borrowed from the other time sequence and, more specifically, is the solution of $ \bw_{1} $ at $ k $th integer time step, see \eqref{Eq: TD1}. Thus it is known. The solution $ \bw^{k+\frac{1}{2}}_{2} $ can be sequentially used for the next, the $ (k+1) $st, integer time step. Thus iterations can proceed.

\subsubsection{Overall temporal discretization}
\begin{figure}[h!]
	\centering
	\includegraphics[width=0.9\linewidth]{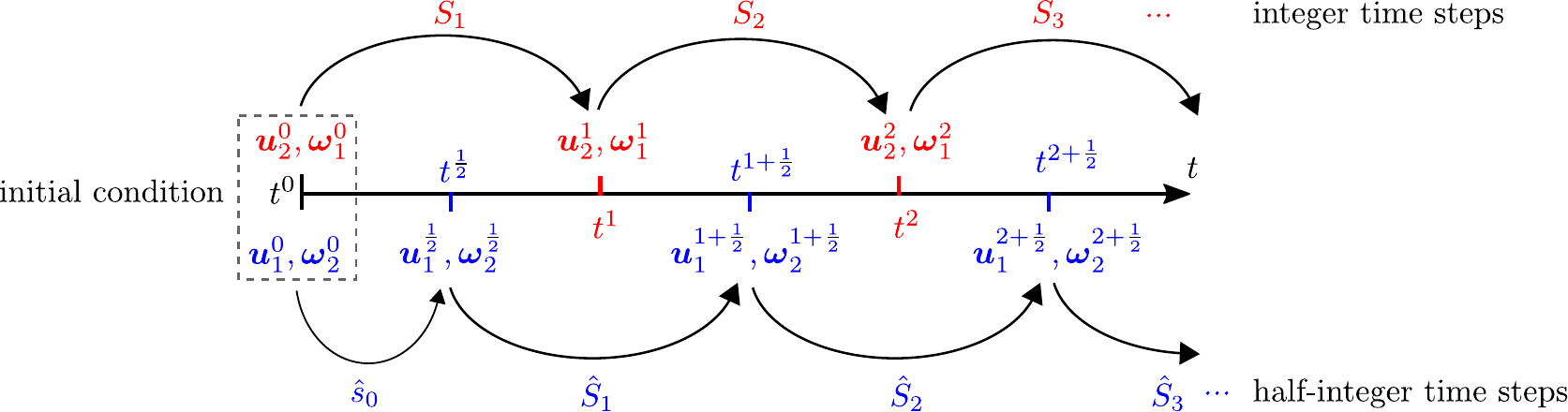}
	\caption{An illustration of the proposed staggered temporal discretization scheme. Integer time steps are denoted by $ S_{k} $, half-integer time steps are denoted by $ \hat{S}_{k} $, and the $ 0 $th time step is denoted by $ \hat{s}_{0} $. The iterations proceed in a sequence: $ \hat{s}_{0}\to S_{1} \to \hat{S}_{1}\to S_{2} \to \hat{S}_{2} \to \cdots $. The $ k $th integer time step also computes $ P_{3}^{k-\frac{1}{2}} $ and the $ k $th half-integer time step also computes $ P_{0}^{k} $.}
	\label{fig: temporal scheme}
\end{figure}

One may notice that to start the iterations we need to know $ \bu_{1}^{\frac{1}{2}}, \bw_{2}^{\frac{1}{2}} $. Therefore we need a $ 0 $th time step, $ \hat{s}_{0} $, computing from $ t^{0} $ to $ t^{\frac{1}{2}} $ for $ \bu_{1}^{\frac{1}{2}}, \bw_{2}^{\frac{1}{2}} $. The simplest approach for the $ 0 $th time step is applying the explicit Euler method to evolution equation \eqref{Eq: WF d} which, together with constraints \eqref{Eq: WF e} and \eqref{Eq: WF f} at $ t^{\frac{1}{2}} $, leads to a semi-discrete system similar to \eqref{Eq: TD2}. More accurate approaches, like directly applying the Gauss integrator \eqref{Eq: gauss integrator} or other (higher order) integrators to formulation \eqref{Eq: WF}, could also be used. These methods will eventually lead to nonlinear discrete algebraic systems for which more expensive iterative methods like the Newton–Raphson method are needed. After the $ 0 $th time step, $ \hat{s}_{0} $, standard iterations, $ S $ and $ \hat{S} $, can proceed.\footnote{It is also fine to switch time sequences for the \MOD{evolution} equations.} The overall temporal scheme is illustrated in Fig.~\ref{fig: temporal scheme}.

It is easy to see that, instead of applying a standard temporal discretization directly to the dual-field mixed weak formulation \eqref{Eq: WF}, using the presented staggered temporal discretization can greatly reduce the computational cost. \REVT{Although} the dual-field formulation doubles the variables, we will only solve for half of them at each time step as the staggered temporal discretization decouples the dual-field formulation. Meanwhile, since each semi-discrete formulation borrows the solution from the other for the nonlinear terms, see the second terms of \eqref{Eq: TD1 a} and \eqref{Eq: TD2 a}, the semi-discrete formulations will lead to \MOD{linearized} discrete algebraic systems.

\subsection{Properties after temporal discretization}\label{Subsub: semi discrete proofs}
In this part, we check whether the conservation (in the inviscid case) and dissipation (in the viscous case) properties proven at the continuous level, see Section~\ref{Sub: continuous conservation proof}, are preserved after the proposed staggered temporal discretization. Note that we have not yet applied a spatial discretization; the function spaces are still the infinite dimensional Hilbert spaces in the de Rham complex, \eqref{Eq: de Rham complex}.
 
\subsubsection{Mass conservation after temporal discretization} \label{Sec: mass conservation semi discrete}
The mass conservation is not influenced by the temporal discretization. \REVO{Due to the same proof as given in Section~\ref{Subsub: continuous mass conservation}, for $ \bu^{k}_{2}\in H(\mathrm{div};\Omega) $, $ \nabla\cdot\bu^{k}_{2}=0 $ is exactly satisfied at all integer time instants, and for $ \bu^{k+\frac{1}{2}}_{1}\in H(\mathrm{curl};\Omega) $ the mass conservation is only weakly imposed through integration by parts, see \eqref{Eq: TD2 c}.}

\subsubsection{Time rate of change of kinetic energy after temporal discretization} Let $ \Rn\to\infty $ and $ \boldsymbol{f} = \boldsymbol{0} $. We replace $ \be_{2} $ by $ \dfrac{\bu^{k-1}_{2} + \bu^{k}_{2}}{2}$ in \eqref{Eq: TD1 a} and replace $ \be_{1} $ by $ \dfrac{\bu^{k+\frac{1}{2}}_{1} + \bu^{k-\frac{1}{2}}_{1}}{2}$ in \eqref{Eq: TD2 a}. Following the same process used for the proof at the continuous level, see Section~\ref{Subsub: continuous kinetic energy conservation}, one can get
\[
\left\langle \dfrac{\bu^{k+\frac{1}{2}}_{1}-\bu^{k-\frac{1}{2}}_{1}}{\varDelta t}, \dfrac{\bu^{k+\frac{1}{2}}_{1} + \bu^{k-\frac{1}{2}}_{1}}{2}\right\rangle _{\Omega} = 0
\quad \text{and}\quad 
\left\langle \dfrac{\bu^{k}_{2}-\bu^{k-1}_{2}}{\varDelta t}, \dfrac{\bu^{k}_{2} + \bu^{k-1}_{2}}{2}\right\rangle _{\Omega} = 0\,,
\]
which then leads to
\begin{equation}\label{Eq: K1 conserving}
	\mathcal{K}_{1}^{k+\frac{1}{2}} = \dfrac{1}{2}\left\langle\bu^{k+\frac{1}{2}}_{1}, \bu^{k+\frac{1}{2}}_{1} \right\rangle _{\Omega} =  \dfrac{1}{2}\left\langle\bu^{k-\frac{1}{2}}_{1}, \bu^{k-\frac{1}{2}}_{1} \right\rangle _{\Omega} = \mathcal{K}_{1}^{k-\frac{1}{2}}\,,
\end{equation}

\begin{equation}\label{Eq: K2 conserving}
\mathcal{K}_{2}^{k} = \dfrac{1}{2}\left\langle\bu^{k}_{2}, \bu^{k}_{2} \right\rangle _{\Omega} =  \dfrac{1}{2}\left\langle\bu^{k-1}_{2}, \bu^{k-1}_{2} \right\rangle _{\Omega} = \mathcal{K}_{2}^{k-1}\,.
\end{equation}
Thus the kinetic energy is preserved at both integer and half-integer time steps.

If $ \Rn<\infty $, with the same analysis, we will obtain 
\[
	\left\langle \dfrac{\bu^{k+\frac{1}{2}}_{1}-\bu^{k-\frac{1}{2}}_{1}}{\varDelta t}, \dfrac{\bu^{k+\frac{1}{2}}_{1} + \bu^{k-\frac{1}{2}}_{1}}{2}\right\rangle _{\Omega} = -\eR\left\langle \dfrac{\bw^{k-\frac{1}{2}}_{2}+\bw^{k+\frac{1}{2}}_{2}}{2}, \nabla\times\dfrac{\bu^{k+\frac{1}{2}}_{1} + \bu^{k-\frac{1}{2}}_{1}}{2}\right\rangle _{\Omega}\,,
\]

\[
\left\langle \dfrac{\bu^{k}_{2}-\bu^{k-1}_{2}}{\varDelta t}, \dfrac{\bu^{k-1}_{2} + \bu^{k}_{2}}{2}\right\rangle _{\Omega} =-\eR\left\langle\nabla\times\dfrac{\bw^{k}_{1}+\bw^{k-1}_{1}}{2}, \dfrac{\bu^{k}_{2} + \bu^{k-1}_{2}}{2}\right\rangle _{\Omega}\,.
\]
With \eqref{Eq: TD1 b}, \eqref{Eq: TD2 b}, we can conclude that
\begin{equation}\label{Eq: energy boundedness K1}
	\dfrac{\mathcal{K}_{1}^{k+\frac{1}{2}}-\mathcal{K}_{1}^{k-\frac{1}{2}}}{\varDelta t} = -\eR\left\langle\dfrac{\bw^{k+\frac{1}{2}}_{2}+\bw^{k-\frac{1}{2}}_{2}}{2},\dfrac{\bw^{k+\frac{1}{2}}_{2}+\bw^{k-\frac{1}{2}}_{2}}{2} \right\rangle _{\Omega} \stackrel{\eqref{Eq: mpr}}{=} - \dfrac{2}{\Rn}\mathcal{E}_{2}^{k}\leq 0\,,
\end{equation}

\begin{equation}\label{Eq: energy boundedness K2}
\dfrac{\mathcal{K}_{2}^{k}-\mathcal{K}_{2}^{k-1}}{\varDelta t} = -\eR\left\langle\dfrac{\bw^{k}_{1}+\bw^{k-1}_{1}}{2}, \dfrac{\bw^{k}_{1}+\bw^{k-1}_{1}}{2} \right\rangle _{\Omega} \stackrel{\eqref{Eq: mpr}}{=} - \dfrac{2}{\Rn}\mathcal{E}_{1}^{k+\frac{1}{2}}\leq 0\,.
\end{equation}
This shows that the boundedness of the kinetic energy \eqref{Eq: boundedness of K1} and \eqref{Eq: boundedness of K2} is preserved by this staggered temporal discretization, which can contribute to the stability of the scheme.

\subsubsection{Time rate of change of helicity after temporal discretization} Let $ \Rn\to\infty $ and $ \boldsymbol{f} = \boldsymbol{0} $. We select $ \be_{1} $ in \eqref{Eq: TD2 a} to be $ \bw_{1}^{k} $ and perform the same process for proving \eqref{Eq: helicity leg 1}. We will get 
\begin{equation}\label{Eq: semi discrete helicity conservation leg 1}
\left\langle \dfrac{\bu_{1}^{k+\frac{1}{2}} - \bu_{1}^{k-\frac{1}{2}}}{\varDelta t}, \bw_{1}^{k}\right\rangle _{\Omega} = 0\, .
\end{equation}
Analogously, by repeating the proof for \eqref{Eq: helicity leg 2}, we can obtain
\begin{equation}\label{Eq: semi discrete helicity conservation leg 2}
\left\langle \dfrac{\bw_{1}^{k} - \bw_{1}^{k-1}}{\varDelta t}, \bu_{1}^{k-\frac{1}{2}}\right\rangle _{\Omega} = 0 \,.
\end{equation}
Equations \eqref{Eq: semi discrete helicity conservation leg 1} and \eqref{Eq: semi discrete helicity conservation leg 2} together imply
\begin{equation}\label{Eq: semi-discrete helicity conservation 1}
	\left\langle \bu_{1}^{k+\frac{1}{2}}, \bw_{1}^{k}\right\rangle _{\Omega}=\left\langle \bu_{1}^{k-\frac{1}{2}}, \bw_{1}^{k}\right\rangle _{\Omega} = \left\langle \bu_{1}^{k-\frac{1}{2}}, \bw_{1}^{k-1}\right\rangle _{\Omega}\,.
\end{equation}
At the half-integer time step $ \hat{S}_{k-1} $ (assume $ k\geq 2 $), \eqref{Eq: semi discrete helicity conservation leg 1} reads
\begin{equation}\label{Eq: special 1}
	\left\langle \dfrac{\bu_{1}^{k-\frac{1}{2}} - \bu_{1}^{k-\frac{3}{2}}}{\varDelta t}, \bw_{1}^{k-1}\right\rangle _{\Omega} = 0 \,.
\end{equation}
With this relation, we can extend \eqref{Eq: semi-discrete helicity conservation 1} to
\[
\left\langle \bu_{1}^{k+\frac{1}{2}}, \bw_{1}^{k}\right\rangle _{\Omega}=\left\langle \bu_{1}^{k-\frac{1}{2}}, \bw_{1}^{k}\right\rangle _{\Omega} = \left\langle \bu_{1}^{k-\frac{1}{2}}, \bw_{1}^{k-1}\right\rangle _{\Omega} = \left\langle \bu_{1}^{k-\frac{3}{2}}, \bw_{1}^{k-1}\right\rangle _{\Omega}\,.
\]
If we further apply the midpoint rule, \eqref{Eq: mpr}, to the first two terms and the last two terms of above \REVT{equation}, we obtain
\begin{equation}\label{Eq: semi-discrete helicity conservation 2}
	\begin{aligned}
		\left\langle \dfrac{\bu_{1}^{k+\frac{1}{2}}+\bu_{1}^{k-\frac{1}{2}}}{2}, \bw_{1}^{k}\right\rangle _{\Omega} \stackrel{\eqref{Eq: mpr}}{=} \left\langle \bu_{1}^{k}, \bw_{1}^{k}\right\rangle _{\Omega} &= \mathcal{H}^{k}_{1} \\
		&= \mathcal{H}^{k-1}_{1} = \left\langle \bu_{1}^{k-1}, \bw_{1}^{k-1}\right\rangle _{\Omega}\stackrel{\eqref{Eq: mpr}}{=} \left\langle \dfrac{\bu_{1}^{k-\frac{1}{2}}+\bu_{1}^{k-\frac{3}{2}}}{2}, \bw_{1}^{k-1}\right\rangle _{\Omega}\,.
	\end{aligned}
\end{equation}
In addition, since \eqref{Eq: TD1 b} holds for all $ \be_{1}\in H(\mathrm{curl};\Omega) $, we can fill $ \bu^{k}_{1} \stackrel{\eqref{Eq: mpr}}{=} \dfrac{\bu^{k+\frac{1}{2}}_{1} + \bu^{k-\frac{1}{2}}_{1}}{2} \in H(\mathrm{curl};\Omega)$ in it and obtain
\begin{equation}\label{Eq: semi-discrete helicity conservation 3}
	\left\langle \bw^{k}_{1} , \bu^{k}_{1}\right\rangle _{\Omega} = \left\langle  \bu^{k}_{2}, \nabla\times\bu^{k}_{1}\right\rangle _{\Omega} = \left\langle  \bu^{k}_{2}, \nabla\times\dfrac{\bu^{k+\frac{1}{2}}_{1} + \bu^{k-\frac{1}{2}}_{1}}{2}\right\rangle _{\Omega} = \left\langle  \bu^{k}_{2}, \dfrac{\bw^{k+\frac{1}{2}}_{2} + \bw^{k-\frac{1}{2}}_{2}}{2}\right\rangle _{\Omega}\,.
\end{equation}
Again, as $ \bw_{2} $ is only solved at half-integer time instants, see \eqref{Eq: TD2}, we use the midpoint rule, \eqref{Eq: mpr}, to bring it to the integer time instants, namely,
$ \bw_{2}^{k} = \dfrac{\bw^{k+\frac{1}{2}}_{2} + \bw^{k-\frac{1}{2}}_{2}}{2} $. As a result, \eqref{Eq: semi-discrete helicity conservation 3} implies
\begin{equation}\label{Eq: H1 equal to H2}
 \mathcal{H}^{k}_{1} = \left\langle \bu_{1}^{k}, \bw_{1}^{k}\right\rangle _{\Omega}= \left\langle \bu_{2}^{k}, \bw_{2}^{k}\right\rangle _{\Omega} = \mathcal{H}^{k}_{2}\, .
\end{equation}
And with \eqref{Eq: semi-discrete helicity conservation 2}, we can finally conclude that
\[
\mathcal{H}^{k}_{1} = \mathcal{H}^{k}_{2} = \mathcal{H}^{k-1}_{1} = \mathcal{H}^{k-1}_{2} = \mathcal{C}\,.
\]


In the viscous case, $ \Rn<\infty $, repeating above analysis at the half-integer time step $ \hat{S}_{k} $ and at the integer time step $ S_{k} $ (see \eqref{Eq: semi discrete helicity conservation leg 1} and \eqref{Eq: semi discrete helicity conservation leg 2}) leads to
\begin{equation}\label{Eq: viscous D HC 1}
	\left\langle \dfrac{\bu_{1}^{k+\frac{1}{2}} - \bu_{1}^{k-\frac{1}{2}}}{\varDelta t}, \bw_{1}^{k}\right\rangle _{\Omega} = -\eR\left\langle \dfrac{\bw^{k+\frac{1}{2}}_{2}+\bw^{k-\frac{1}{2}}_{2}}{2}, \nabla\times\bw_{1}^{k}\right\rangle _{\Omega}\,,
\end{equation}

\begin{equation}\label{Eq: viscous D HC 2}
	\left\langle \dfrac{\bw_{1}^{k} - \bw_{1}^{k-1}}{\varDelta t}, \bu_{1}^{k-\frac{1}{2}}\right\rangle _{\Omega} = - \eR\left\langle\nabla\times\dfrac{\bw^{k}_{1}+\bw^{k-1}_{1}}{2}, \bw_{2}^{k-\frac{1}{2}}\right\rangle _{\Omega}\,.
\end{equation}
\clearpage
\noindent If we combine the above two equations, i.e., \eqref{Eq: viscous D HC 1} + \eqref{Eq: viscous D HC 2}, and use the midpoint rule, \eqref{Eq: mpr}, we obtain
\begin{equation}\label{Eq: viscous D HC --- 0}
	\dfrac{1}{\varDelta t}\left\langle \bu_{1}^{k+\frac{1}{2}}, \bw_{1}^{k}\right\rangle_{\Omega} - \dfrac{1}{\varDelta t}\left\langle\bu_{1}^{k-\frac{1}{2}}, \bw_{1}^{k-1} \right\rangle_{\Omega} =  -\eR\left\langle\bw^{k}_{2}, \nabla\times\bw_{1}^{k}\right\rangle _{\Omega}-\eR\left\langle\nabla\times\bw^{k-\frac{1}{2}}_{1}, \bw_{2}^{k-\frac{1}{2}}\right\rangle _{\Omega}\,.
\end{equation}
Again, \eqref{Eq: viscous D HC 1} is still valid at $ (k-1) $st half-integer time step, $ \hat{S}_{k-1} $, (assume $ k\geq2 $) where it reads
\begin{equation}\label{Eq: viscous D HC 3}
	\left\langle \dfrac{\bu_{1}^{k-\frac{1}{2}} - \bu_{1}^{k-\frac{3}{2}}}{\varDelta t}, \bw_{1}^{k-1}\right\rangle _{\Omega} = -\eR\left\langle \dfrac{\bw^{k-\frac{1}{2}}_{2}+\bw^{k-\frac{3}{2}}_{2}}{2}, \nabla\times\bw_{1}^{k-1}\right\rangle _{\Omega}\,.
\end{equation}
If we now combine \eqref{Eq: viscous D HC 2} and \eqref{Eq: viscous D HC 3}, and use the midpoint rule, \eqref{Eq: mpr}, we get
\begin{equation}\label{Eq: viscous D HC --- 1}
	\dfrac{1}{\varDelta t}\left\langle \bu_{1}^{k-\frac{1}{2}}, \bw_{1}^{k}\right\rangle_{\Omega} - \dfrac{1}{\varDelta t}\left\langle \bu_{1}^{k-\frac{3}{2}}, \bw_{1}^{k-1}\right\rangle_{\Omega} =  -\eR\left\langle\bw^{k-1}_{2}, \nabla\times\bw_{1}^{k-1}\right\rangle _{\Omega}-\eR\left\langle\nabla\times\bw^{k-\frac{1}{2}}_{1}, \bw_{2}^{k-\frac{1}{2}}\right\rangle _{\Omega}\,.
\end{equation}
We now can combine \eqref{Eq: viscous D HC --- 0} and \eqref{Eq: viscous D HC --- 1} and obtain
\[
	\begin{aligned}
		\dfrac{\left\langle \bu_{1}^{k+\frac{1}{2}}+ \bu_{1}^{k-\frac{1}{2}}, \bw_{1}^{k}\right\rangle_{\Omega}}{\varDelta t} - 
		&\dfrac{\left\langle  \bu_{1}^{k-\frac{1}{2}}+ \bu_{1}^{k-\frac{3}{2}},\bw_{1}^{k-1}\right\rangle_{\Omega}}{\varDelta t} \\
		&\stackrel{\eqref{Eq: mpr}}{=} 
		2\dfrac{\left\langle \bu_{1}^{k}, \bw_{1}^{k}\right\rangle_{\Omega} - \left\langle  \bu_{1}^{k-1},\bw_{1}^{k-1}\right\rangle_{\Omega}}{\varDelta t}\\
		&= 2\dfrac{\mathcal{H}^{k}_{1}-\mathcal{H}^{k-1}_{1}}{\varDelta t} \\
		& = -\dfrac{2}{\Rn}\left\langle\nabla\times\bw^{k-\frac{1}{2}}_{1}, \bw_{2}^{k-\frac{1}{2}}\right\rangle _{\Omega} -\dfrac{\left\langle\bw^{k}_{2}, \nabla\times\bw_{1}^{k}\right\rangle _{\Omega}+\left\langle\bw^{k-1}_{2}, \nabla\times\bw_{1}^{k-1}\right\rangle _{\Omega}}{\Rn} \,.
	\end{aligned}
\]
Finally, because \eqref{Eq: H1 equal to H2} still holds, the viscosity dissipates $ \mathcal{H}^{k}_{1} $ and $ \mathcal{H}^{k}_{2} $ at the same rate, denoted by
\begin{equation}\label{Eq: discrete helicity disspation H1 H2 1}
	\begin{aligned}
	\mathscr{D}(\bw_{1},\bw_{2}):=\dfrac{\mathcal{H}^{k}_{1}-\mathcal{H}^{k-1}_{1}}{\varDelta t} = \dfrac{\mathcal{H}^{k}_{2}-\mathcal{H}^{k-1}_{2}}{\varDelta t}= 
	-&\dfrac{\left\langle\nabla\times\bw^{k-\frac{1}{2}}_{1}, \bw_{2}^{k-\frac{1}{2}}\right\rangle _{\Omega}}{\Rn} \\
	&-\dfrac{\left\langle\bw^{k}_{2}, \nabla\times\bw_{1}^{k}\right\rangle _{\Omega}+\left\langle\bw^{k-1}_{2}, \nabla\times\bw_{1}^{k-1}\right\rangle _{\Omega}}{2\Rn}\,.
	\end{aligned}
\end{equation}

\section{Mimetic spatial discretization}\label{Sec: spatial discretization}
It has been shown that the de Rham complex plays an essential role in the proofs and analysis of the conservation properties and the dissipation rates for the proposed dual-field formulation at both continuous and semi-discrete levels. For example, \eqref{Eq: H1 importance} is valid because we have chosen $ P_{0}\in H^1(\Omega) $ such that $ \nabla P_{0} \in H(\mathrm{curl};\Omega) $ is guaranteed. Choosing $ \bu_{1}\in H(\mathrm{curl};\Omega) $ and $ \bw_{2}\in H(\mathrm{div};\Omega) $ ensures that the relation $ \bw_{2} = \nabla\times\bu_{1} $ is satisfied exactly. In addition, as shown in Section~\ref{Subsub: continuous mass conservation}, the de Rham complex is essential for the mass conservation $ \nabla\cdot\bu_{2}=0 $ where $ \bu_{2}\in H(\mathrm{div};\Omega) $. 

In this work we consider a set of discrete function spaces,
\[\left\lbrace G(\Omega), C(\Omega), D(\Omega), S(\Omega) \right\rbrace \,, \] 
where
\[
 G(\Omega) \subset H^1(\Omega)\,,\quad
 C(\Omega) \subset H(\mathrm{curl};\Omega)\,,\quad
 D(\Omega) \subset H(\mathrm{div};\Omega)\,,\quad
 S(\Omega) \subset L^2(\Omega)\,,
\]
such that
\begin{equation}\label{eq: discrete de Rham complex}
		\mathbb{R}\hookrightarrow
	G(\Omega) \stackrel{\nabla}{\longrightarrow} 
	C(\Omega) \stackrel{\nabla\times}{\longrightarrow} 
	D(\Omega) \stackrel{\nabla\cdot}{\longrightarrow} 
	S(\Omega)
	\rightarrow 0\,,
\end{equation}
i.e., they constitute a discrete de Rham complex. In order to enable the validity of the proofs and analysis at the fully discrete level, we need to employ such a set of discrete spaces for the spatial discretization.

Any sequence of discrete function spaces that satisfies \eqref{eq: discrete de Rham complex} is equally valid. 
One possible choice is to employ $G(\Omega) = \mathrm{CG}_{N}$, $C(\Omega) = \mathrm{NED}^{1}_{N}$, $D(\Omega) = \mathrm{RT}_{N}$, and $S(\Omega) = \mathrm{DG}_{N-1}$, where $\mathrm{CG}_{N}$ are the Lagrange polynomials of degree $N$, $\mathrm{NED}^{1}_{N}$ are the N\'{e}d\'{e}lec $H(\mathrm{curl})$-conforming spaces of the first kind of degree $N$, see \cite{nedelec}, $\mathrm{RT}_{N}$ are the Raviart-Thomas spaces of degree $N$, see \cite{nedelec, raviart606mixed}, and $\mathrm{DG}_{N-1}$ are the discontinuous Lagrange spaces of degree $(N-1)$. Another possible exact sequence of discrete function spaces employing b-splines is employed in the works by Hiemstra et al. \cite{Hiemstra2014}, Buffa et al. \cite{BuffaSangalliRivasVazquez2011}, and Ratnani and Sonnendr\"{u}cker \cite{Ratnani2012}. We call these spaces \emph{structure-preserving} or \emph{mimetic} spaces, see another example, the \emph{mimetic polynomial spaces} \cite{ Palha2017a, kreeft::stokes, palha2014physics, Kreeft2011, gerritsma2011edge}. \REVO{Note that variables in the finite dimensional spaces $ C(\Omega) $ and $ D(\Omega) $ possess the regularity that ensures the $ L^2 $-integrability of the convective terms in the weak formulation \eqref{Eq: WF}, see Remark~\ref{remark L2 integrable}.}

\subsection{Fully discrete systems} \label{Sub: 3d discretization}
Applying a particular set of mimetic spaces to the semi-discrete problems \eqref{Eq: TD1} and \eqref{Eq: TD2} leads to two local fully discrete linear algebraic systems, one for the $ k $th integer time step $ S_{k} $, i.e., \vspace{-3ex}

\begin{subequations} \label{Eq: full dis 1}
	\begin{align}
		&\mathsf{N} \dfrac{\vec{\bu}^{k}_{2}-\vec{\bu}^{k-1}_{2}}{\varDelta t} + \mathsf{R}^{k-\frac{1}{2}} \dfrac{\vec{\bu}^{k}_{2}+\vec{\bu}^{k-1}_{2}}{2} + \eR\mathsf{C}\dfrac{\vec{\bw}^{k}_{1}+\vec{\bw}^{k-1}_{1}}{2}- \mathsf{D} ^{\mathsf{T}} \vec{P}^{k-\frac{1}{2}}_{3} = {\boldsymbol{\mathsf{f}}}^{k-\frac{1}{2}}\,, \\
		&\mathsf{C}^{\mathsf{T}} \vec{\bu}^{k}_{2}-\mathsf{M} \vec{\bw}^{k}_{1} =\boldsymbol{0} \,,\\
		&\mathsf{D}\vec{\bu}^{k}_{2} = \boldsymbol{0} \,, 
	\end{align}
\end{subequations}
and one for the $ k $th half-integer time step $ \hat{S}_{k} $, namely,\vspace{-3ex}

\begin{subequations}\label{Eq: full dis 2}
\begin{align}
	&\mathsf{M}\dfrac{\vec{\bu}^{k+\frac{1}{2}}_{1}-\vec{\bu}^{k-\frac{1}{2}}_{1}}{\varDelta t}  + \mathsf{R}^{k} \dfrac{\vec{\bu}^{k+\frac{1}{2}}_{1}+\vec{\bu}^{k-\frac{1}{2}}_{1}}{2} + \eR\mathsf{C}^{\mathsf{T}} \dfrac{\vec{\bw}^{k+\frac{1}{2}}_{2}+\vec{\bw}^{k-\frac{1}{2}}_{2}}{2} + \mathsf{G}\vec{P}^{k}_{0} = {\boldsymbol{\mathsf{f}}}^{k} \,,\\
	&\mathsf{C} \vec{\bu}^{k+\frac{1}{2}}_{1} - \mathsf{N} \vec{\bw}^{k+\frac{1}{2}}_{2}  = \boldsymbol{0}\,, \label{Eq: full dis 2 b}\\
	&\mathsf{G}^{\mathsf{T}} \vec{\bu}^{k+\frac{1}{2}}_{1} = \boldsymbol{0}\,, 
\end{align}
\end{subequations}
where we have used the vector sign to indicate the vector of the expansion coefficients of a discrete variable. And, if $ \gamma, \boldsymbol{\tau}, \boldsymbol{\sigma},  \chi $ are basis functions of mimetic spaces 
\[
G(\Omega)\,,\quad 
C(\Omega)\,,\quad
D(\Omega)\,,\quad
S(\Omega)\,,
\]
respectively, $ \mathsf{M}$ and  $\mathsf{N}$ are the symmetric mass (or stiffness) matrices of spaces $ C(\Omega)$ and $ D(\Omega)$,
\[
\mathsf{M}_{ij} = \left\langle \boldsymbol{\tau}_{j},\  \boldsymbol{\tau}_{i}\right\rangle _{\Omega}\quad \text{and}\quad 
\mathsf{N}_{ij} = \left\langle \boldsymbol{\sigma}_{j},\  \boldsymbol{\sigma}_{i}\right\rangle _{\Omega}\,,
\]
and the entries of matrices $ \mathsf{C}, \mathsf{D}, \mathsf{G},\mathsf{R}^{k-\frac{1}{2}}, \mathsf{R}^{k} $ and vectors $  {\boldsymbol{\mathsf{f}}}^{k-\frac{1}{2}} , {\boldsymbol{\mathsf{f}}}^{k}$ are
\[
\begin{aligned}
	&\mathsf{C}_{ij} = \left\langle\nabla\times\boldsymbol{\tau}_{j},\  \boldsymbol{\sigma}_{i}\right\rangle _{\Omega}\,,
	\\
	&\mathsf{D}_{ij} = \left\langle\nabla\cdot \boldsymbol{\sigma}_{j},\  \chi_{i} \right\rangle _{\Omega}\,,
	\\
	&\mathsf{G}_{ij} = \left\langle\nabla \gamma_{j},\  \boldsymbol{\tau}_{i}\right\rangle _{\Omega}\,,
	\\
	&\mathsf{R}^{k-\frac{1}{2}}_{ij}  = \left\langle  \bw_{2}^{k-\frac{1}{2}}\times \boldsymbol{\sigma}_{j},\ \boldsymbol{\sigma}_{i} \right\rangle _{\Omega}\,,
	\\
	&\mathsf{R}^{k}_{ij}  =\left\langle  \bw_{1}^{k}\times \boldsymbol{\tau}_{j},\ \boldsymbol{\tau}_{i} \right\rangle _{\Omega}\,,
	\\
	&\boldsymbol{\mathsf{f}}^{k-\frac{1}{2}}_{i} = \left\langle \boldsymbol{f}^{k-\frac{1}{2}},\boldsymbol{\sigma}_{i} \right\rangle  _{\Omega}\,,
	\\
	&\boldsymbol{\mathsf{f}}^{k}_{i} = \left\langle \boldsymbol{f}^{k},\boldsymbol{\tau}_{i} \right\rangle  _{\Omega}.
\end{aligned}
\]
If we rearrange the systems \eqref{Eq: full dis 1} and \eqref{Eq: full dis 2} and write them in linear algebra format, we can obtain following linear systems, 
\[
	\begin{bmatrix}
		\dfrac{1}{\varDelta t}\mathsf{N} + \dfrac{1}{2}\mathsf{R}^{k-\frac{1}{2}} & \dfrac{1}{2\Rn}\mathsf{C} & -\mathsf{D}^{\mathsf{T}}\\
		\mathsf{C}^{\mathsf{T}} & -\mathsf{M} & \boldsymbol{0}\\
		\mathsf{D} & \boldsymbol{0} & \boldsymbol{0}
	\end{bmatrix}
\begin{bmatrix}
	\vec{\bu}_{2}^{k} \\
	\vec{\bw}_{1}^{k} \\
	\vec{P}_{3}^{k-\frac{1}{2}} \\
\end{bmatrix}
=
\begin{bmatrix}
	\left( \dfrac{1}{\varDelta t}\mathsf{N}-\dfrac{1}{2}\mathsf{R}^{k-\frac{1}{2}}\right) \vec{\bu}_{2}^{k-1}-\dfrac{1}{2\Rn}\mathsf{C}\ \vec{\bw}_{1}^{k-1}+ {\boldsymbol{\mathsf{f}}}^{k-\frac{1}{2}} \\
	\boldsymbol{0} \\
	\boldsymbol{0} \\
\end{bmatrix}\,,
\]

\[
	\begin{bmatrix}
		\dfrac{1}{\varDelta t}\mathsf{M} + \dfrac{1}{2}\mathsf{R}^{k} & \dfrac{1}{2\Rn}\mathsf{C}^{\mathsf{T}} & \mathsf{G}\\
		\mathsf{C} & -\mathsf{N} & \boldsymbol{0} \\
		\mathsf{G}^{\mathsf{T}} & \boldsymbol{0} & \boldsymbol{0}
	\end{bmatrix}
	\begin{bmatrix}
		\vec{\bu}_{1}^{k+\frac{1}{2}} \\
		\vec{\bw}_{2}^{k+\frac{1}{2}} \\
		\vec{P}_{0}^{k} \\
	\end{bmatrix}
	=
	\begin{bmatrix}
		\left( \dfrac{1}{\varDelta t}\mathsf{M}-\dfrac{1}{2}\mathsf{R}^{k}\right) \vec{\bu}_{1}^{k-\frac{1}{2}} - \dfrac{1}{2\Rn}\mathsf{C}^{\mathsf{T}}\vec{\bw}_{2}^{k-\frac{1}{2}} + {\boldsymbol{\mathsf{f}}}^{k}\\
		\boldsymbol{0} \\
		\boldsymbol{0} \\
	\end{bmatrix}\,.
\]
\MOD{A} similar spatial discretization can be applied to the semi-discrete system for the $ 0 $th time step $ \hat{s}_{0} $, see Fig.~\ref{fig: temporal scheme}. 

Suppose a mesh has been generated in the computational domain $ \Omega $. We can perform such discretizations in all elements. After applying the initial condition and assembling the local systems, we will eventually obtain global linear systems ready to be solved in the sequence shown in Fig.~\ref{fig: temporal scheme}.

\subsection{Properties of the fully discrete systems} \label{Sub: proofs of fully discrete system}
Since we have used a sequence of function spaces which form a discrete de Rham complex, the proofs for the conservation properties and the analysis for the dissipation rates of kinetic energy and helicity at the semi-discrete level, see Section~\ref{Subsub: semi discrete proofs}, remain valid at the fully discrete level.

\section{Numerical experiments} \label{Sec: Numerical results}
We now test the proposed mimetic dual-field method with two manufactured solutions and a more general flow, the well-known Taylor-Green vortex. 

For all tests, we use the mimetic polynomial spaces as our mimetic spaces and do the spatial discretization under the framework of the 
MSEM. Meshes are uniform orthogonal structured hexahedral meshes. The mesh size, namely, the edge length of the cubic element cell, is denoted by $ h $. The degree of the mimetic polynomials is denoted by $ N $. And we use the explicit Euler method for the temporal discretization of the $ 0 $th time step, i.e., $ \hat{s}_{0} $ in Fig.~\ref{fig: temporal scheme}. The implementation is conducted in \emph{Python}.

\subsection{Manufactured solution tests} \label{Sub: TESTS anufactured solutions}
Two manufactured solutions are taken from \cite{Rebholz2007}; one for testing the conservation properties and one for investigating the convergence rate of the method. The domain is selected to be the periodic unit cube $ \Omega:=[0,1]^3 $. 
\subsubsection{Conservation properties and dissipation rates}
For these first tests, we select the initial condition 
\[\left.\boldsymbol{u}\right|_{t=0} = \left[\cos(2\pi z),\ \sin(2\pi z),\ \sin(2\pi x)\right] ^{\mathsf{T}}\,.\]
Such an initial condition possesses kinetic energy $ \left.\mathcal{K}\right|_{t=0}=0.75 $ and helicity $ \left.\mathcal{H}\right|_{t=0}=-6.283 $. The problem is solved until $ t=10 $ on an extremely coarse mesh of $ h=1/3 $ and $ N=2 $.

\begin{figure}[h!]
	\centering{
		\subfloat{
			\begin{minipage}[b]{0.5\textwidth}
				\centering
				\includegraphics[width=0.9\linewidth]{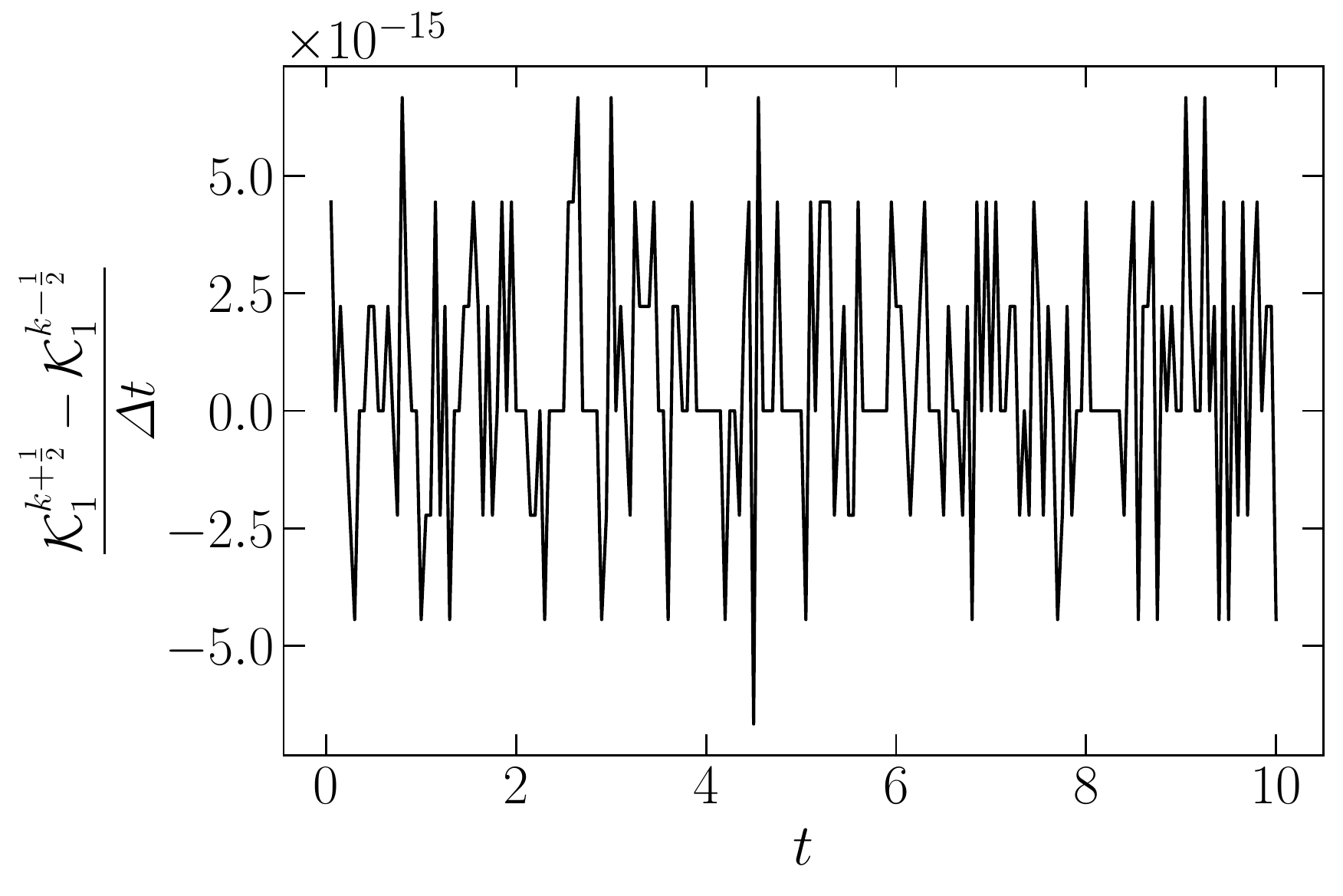}
			\end{minipage}
		}
		\subfloat{
			\begin{minipage}[b]{0.5\textwidth}
				\centering
				\includegraphics[width=0.83\linewidth]{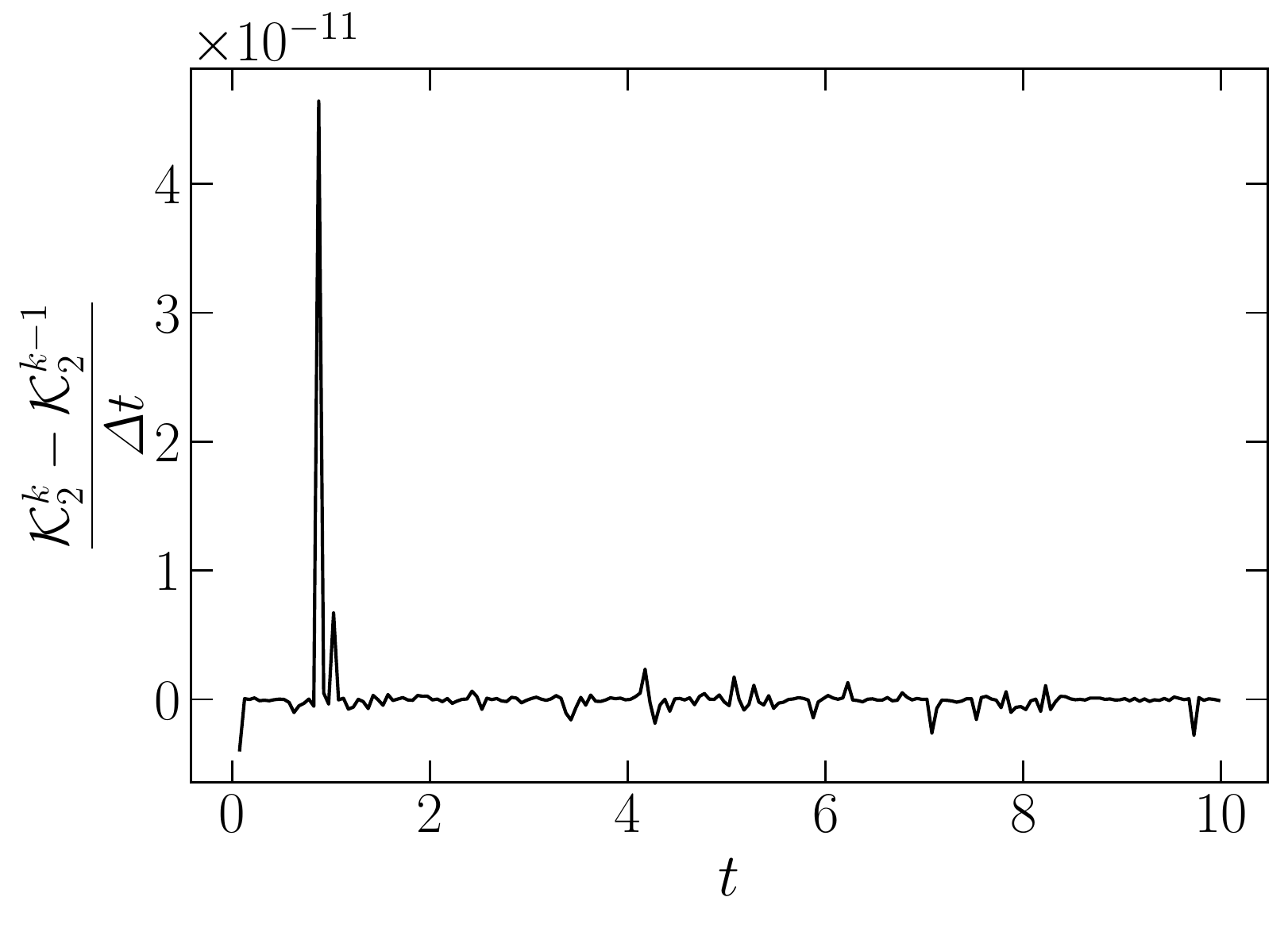}
			\end{minipage}
		}\\
		\subfloat{
			\begin{minipage}[b]{0.5\textwidth}
				\centering
				\includegraphics[width=0.87\linewidth]{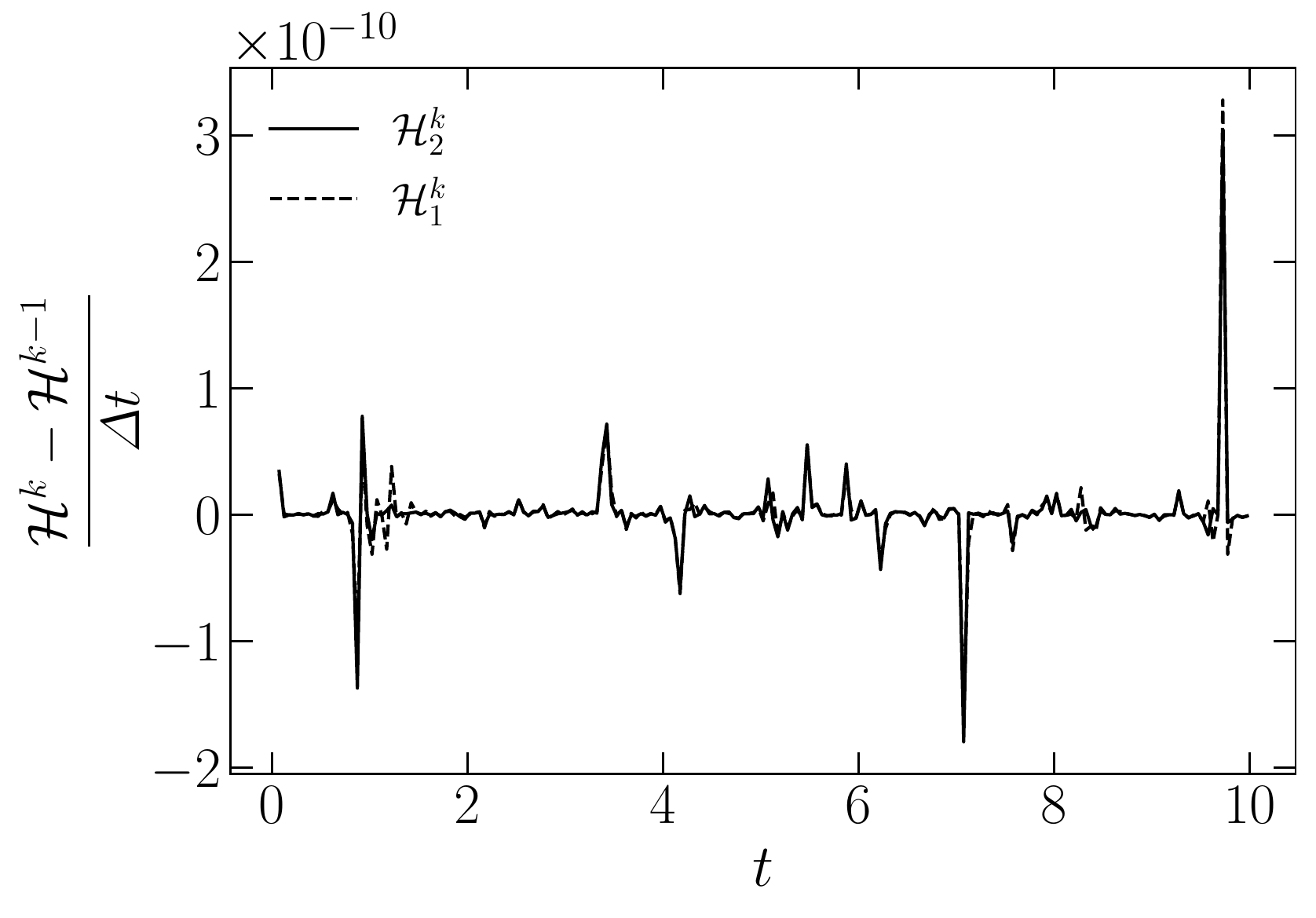}
			\end{minipage}
		}
		\subfloat{
			\begin{minipage}[b]{0.5\textwidth}
				\centering
				\includegraphics[width=0.81\linewidth]{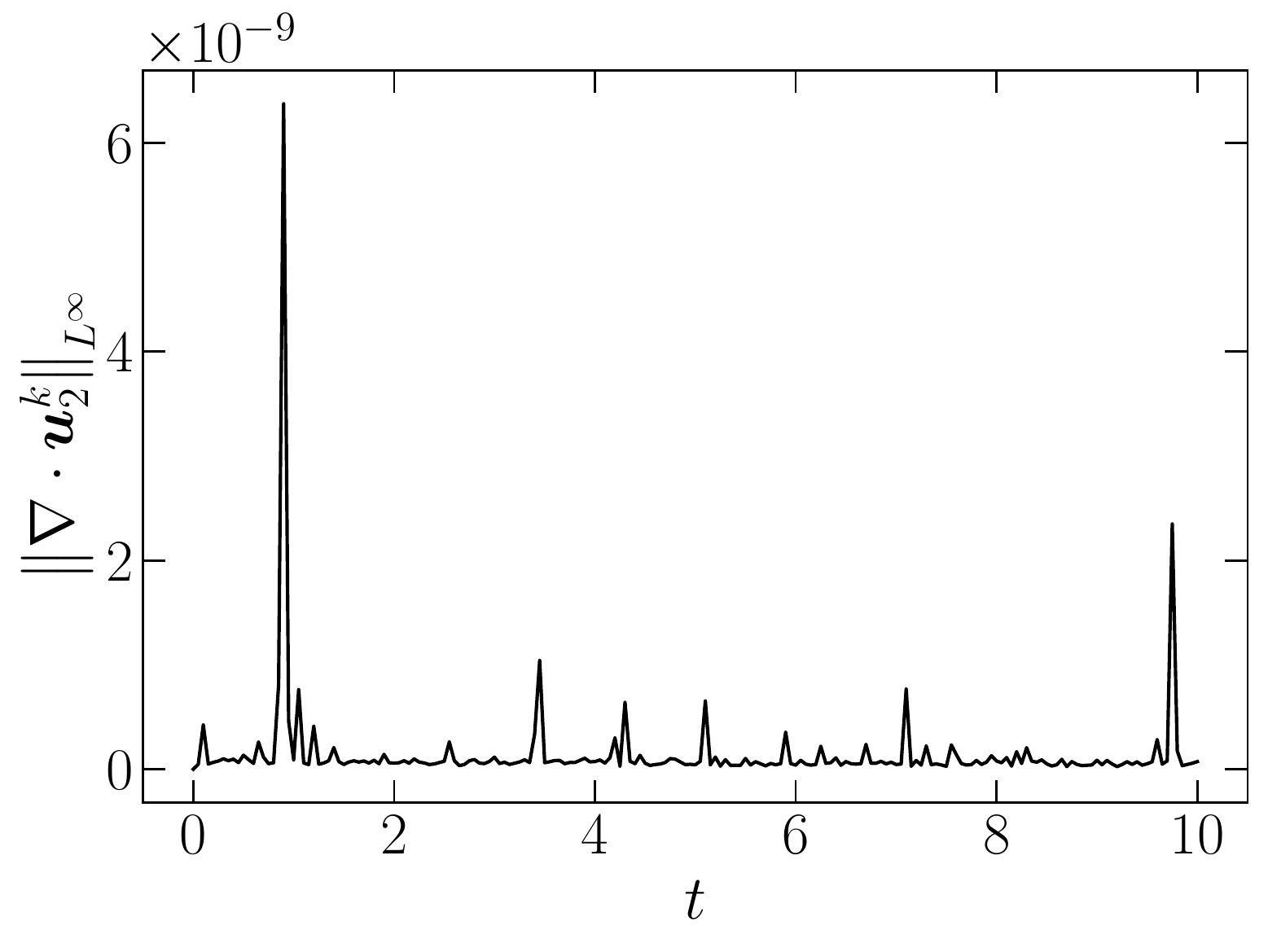}
			\end{minipage}
		}\\
		\caption{Some results of the conservation test for $ h=1/3 $, $ N=2 $ and $ \varDelta t=1/20 $.}
		\label{fig: conservation test 1}}
\end{figure}

We first try to verify that the proposed method does preserve mass, kinetic energy and helicity if in the inviscid limit $ \Rn\to\infty $ and $ \boldsymbol{f}=\boldsymbol{0} $. In Fig.~\ref{fig: conservation test 1} some results are presented. The results of $ \left\|\nabla\cdot\bu_{2}^{k}\right\|_{L^\infty} $ in the bottom-right diagram imply that the pointwise mass conservation is always satisfied. In the bottom-left diagram, the results show that both $ \mathcal{H}^{k}_{1} $ and $ \mathcal{H}^{k}_{2} $ are preserved. The fact that the two lines coincide with each other up to $ \mathcal{O}(10^{-10}) $ verifies \eqref{Eq: H1 equal to H2}. As for kinetic energy, the results are present in the top diagrams where the discrete conservation for both $ \mathcal{K}^{k-\frac{1}{2}}_{1} $ and $ \mathcal{K}^{k}_{2} $ at their corresponding time steps are shown.

\begin{figure}[h!]
	\centering{
		\subfloat{
			\begin{minipage}[b]{0.5\textwidth}
				\centering
				\includegraphics[width=0.88\linewidth]{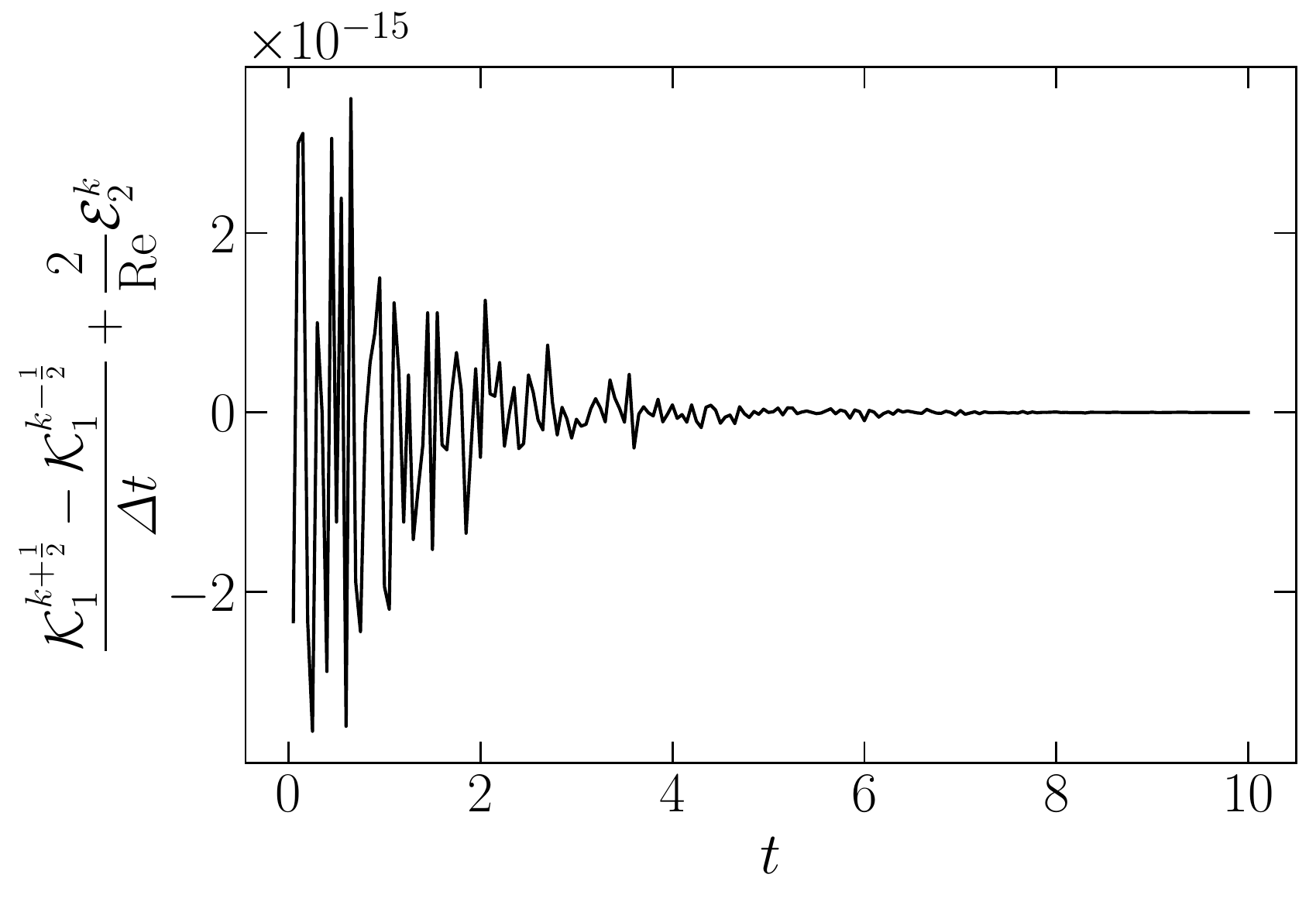}
			\end{minipage}
		}
		\subfloat{
			\begin{minipage}[b]{0.5\textwidth}
				\centering
				\includegraphics[width=0.87\linewidth]{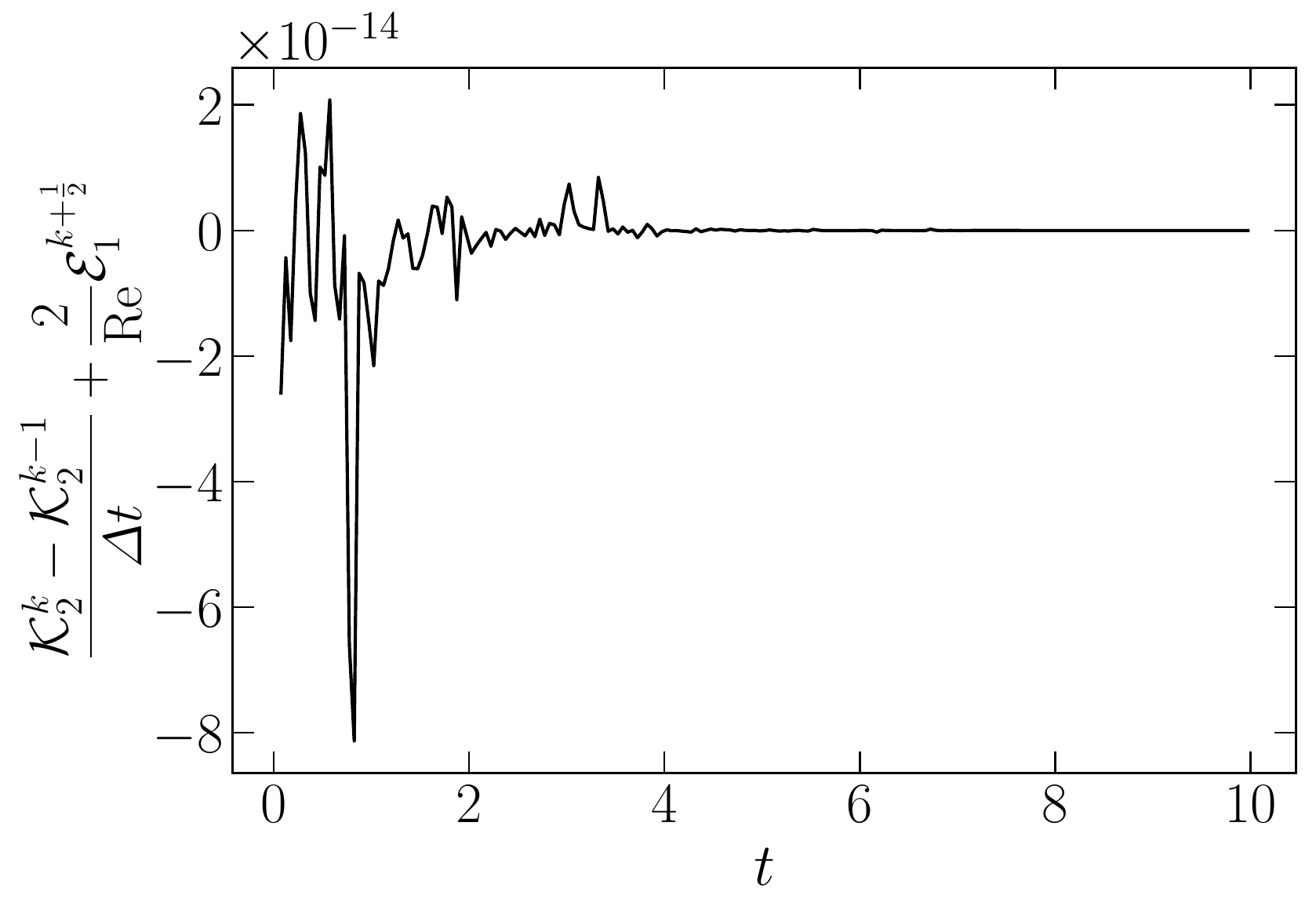}
			\end{minipage}
		}\\
		\subfloat{
			\begin{minipage}[b]{0.5\textwidth}
				\centering
				\includegraphics[width=0.87\linewidth]{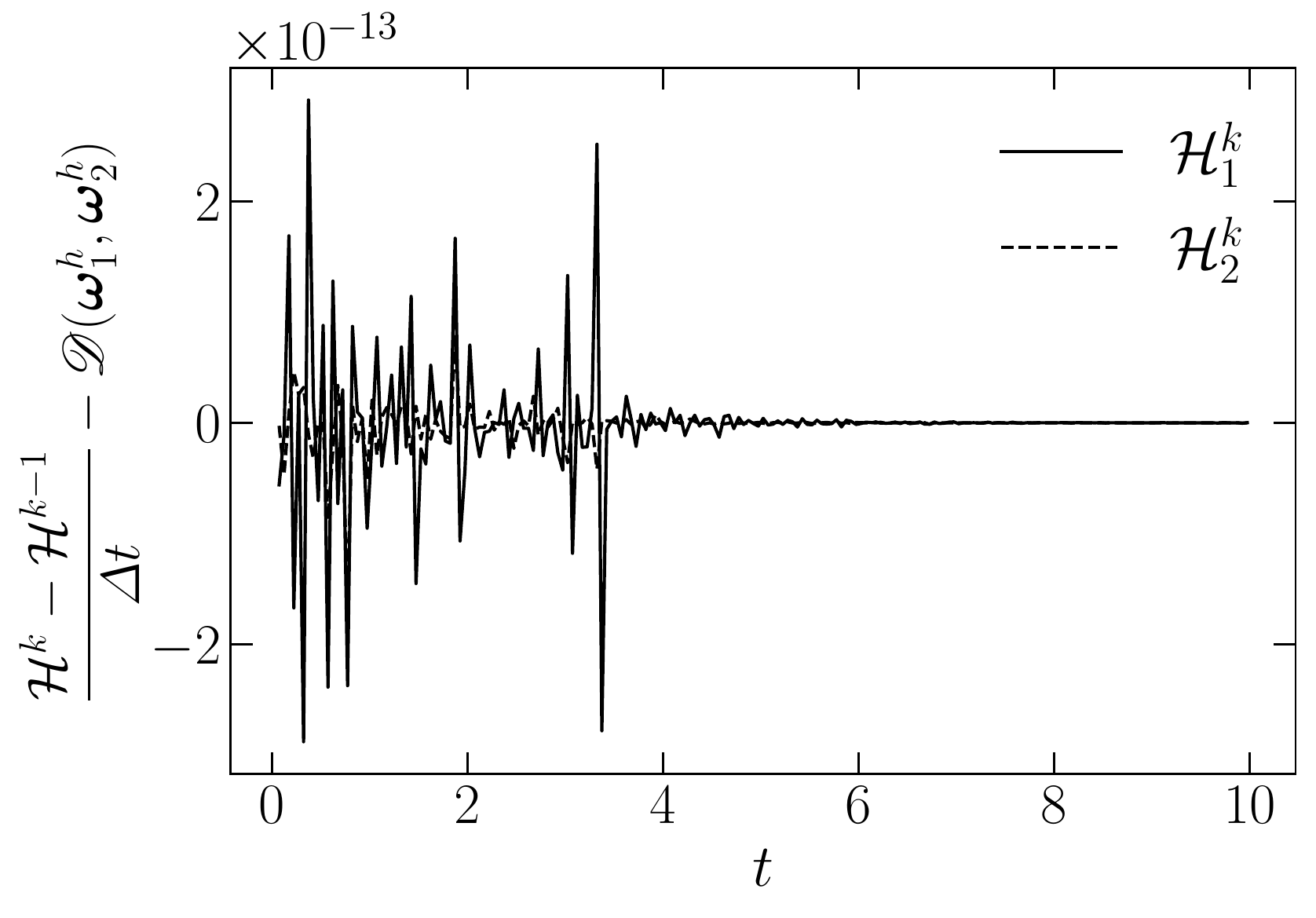}
			\end{minipage}
		}
		\subfloat{
			\begin{minipage}[b]{0.5\textwidth}
				\centering
				\includegraphics[width=0.85\linewidth]{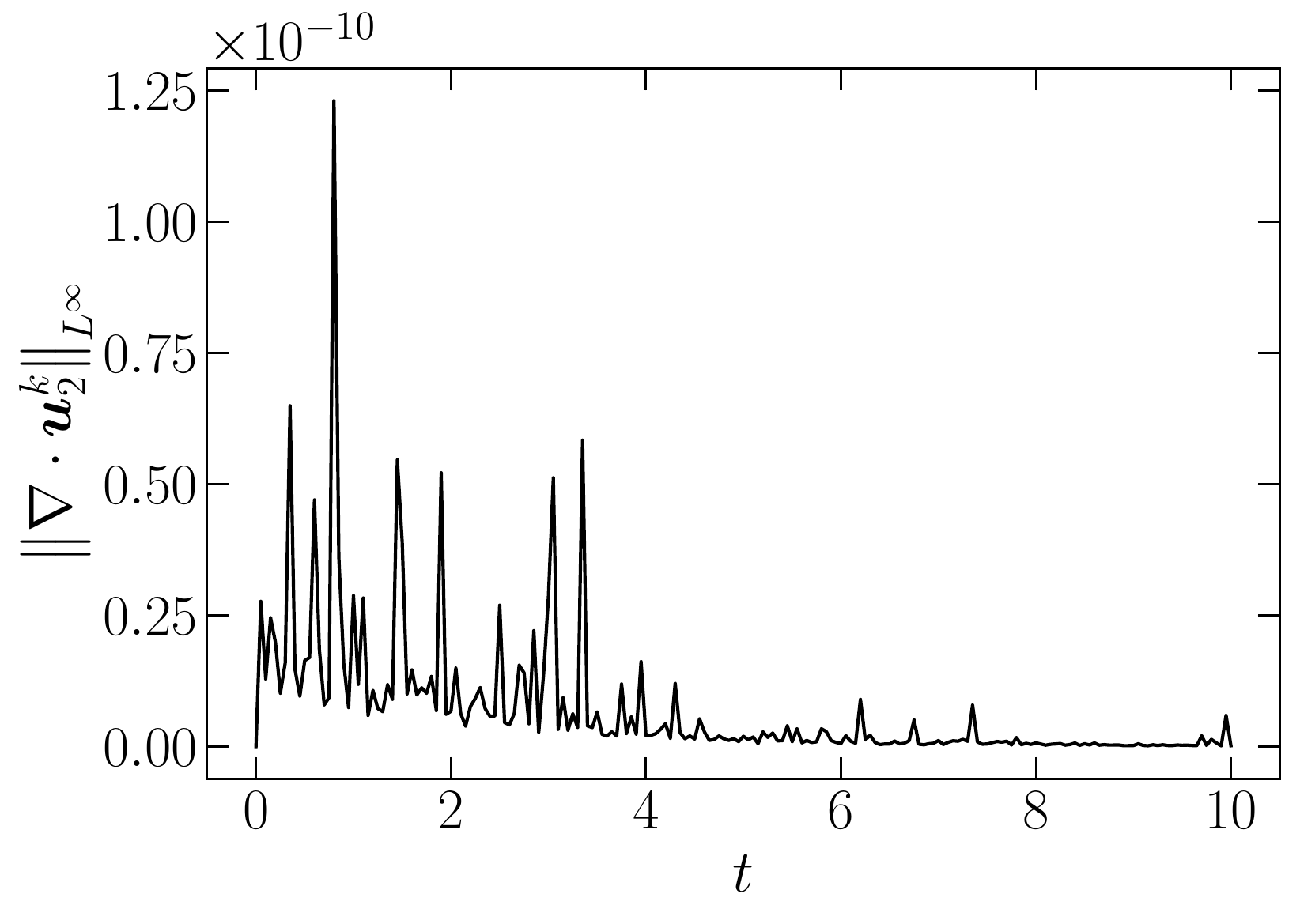}
			\end{minipage}
		}\\
		\caption{\REVT{Some results of the dissipation test for $ \Rn=100 $, $ h=1/3 $, $ N=2 $ and $ \varDelta t=1/20 $. $ \mathscr{D}(\bw_{1}^{h},\bw_{2}^{h}) $ is the dissipation rate of helicity, see \eqref{Eq: discrete helicity disspation H1 H2 1}.}}
		\label{fig: conservation test 2 viscous}}
\end{figure}

We then keep $ \boldsymbol{f}=\boldsymbol{0} $ and use \MOD{a} $ \Rn<\infty $; we let the \MOD{viscosity} dissipate kinetic energy and helicity. Some results for $ \Rn=100 $ are presented in Fig.~\ref{fig: conservation test 2 viscous} where the results shown in the top diagrams verify the dissipation rate of kinetic energy derived in \eqref{Eq: energy boundedness K1} and \eqref{Eq: energy boundedness K2} and the results in the bottom-left diagram are in agreement with the dissipation rate of helicity, see \eqref{Eq: discrete helicity disspation H1 H2 1}. The pointwise conservation of mass is still satisfied at all time steps as shown in the bottom-right diagram of Fig.~\ref{fig: conservation test 2 viscous}.

\begin{figure}[h!]
	\centering{
		\subfloat[$ \Rn\to\infty $]{
			\begin{minipage}[b]{0.5\textwidth}
				\centering
				\includegraphics[width=0.9\linewidth]{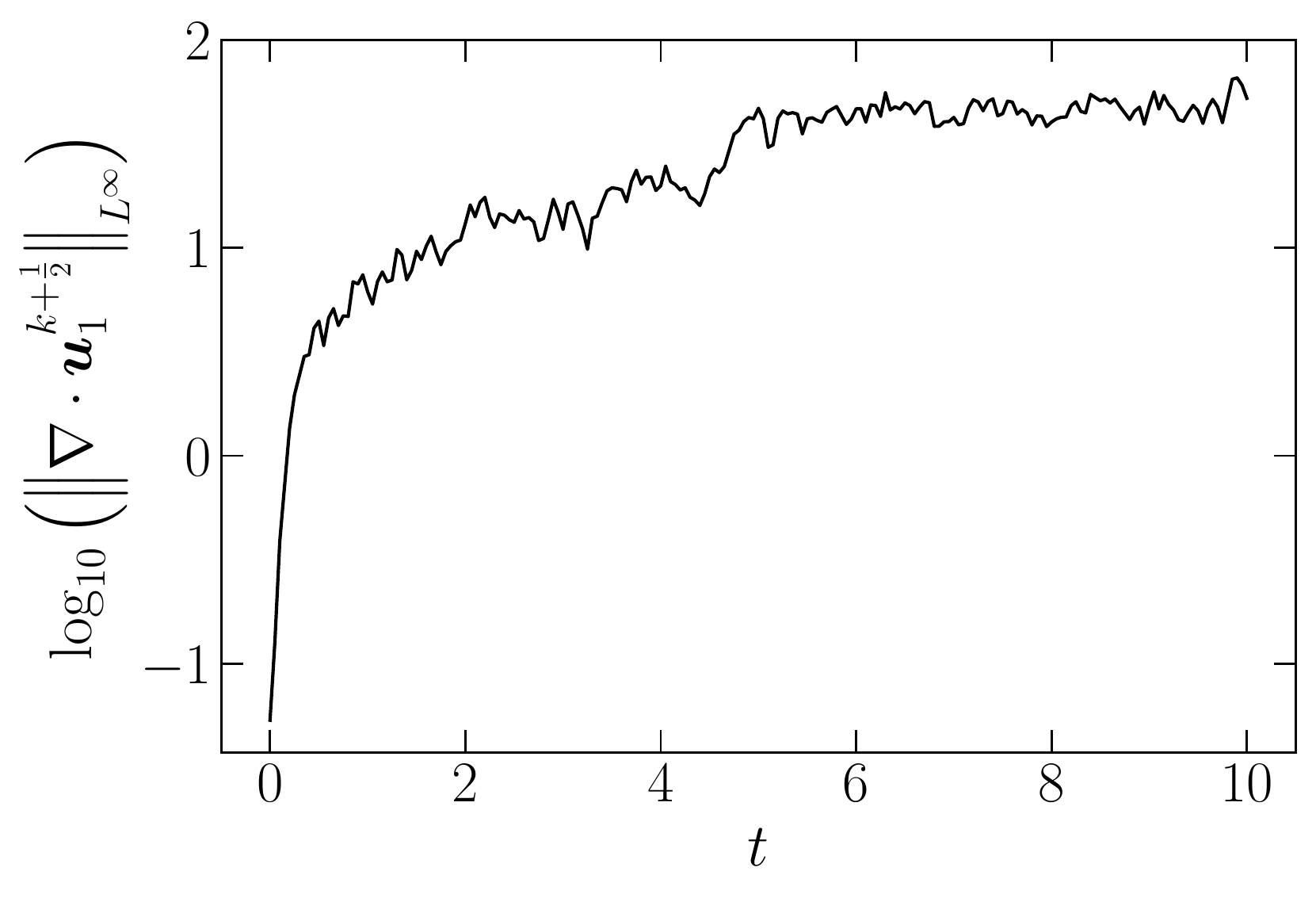}
			\end{minipage}
		}
		\subfloat[$ \Rn=100 $]{
			\begin{minipage}[b]{0.5\textwidth}
				\centering
				\includegraphics[width=0.9\linewidth]{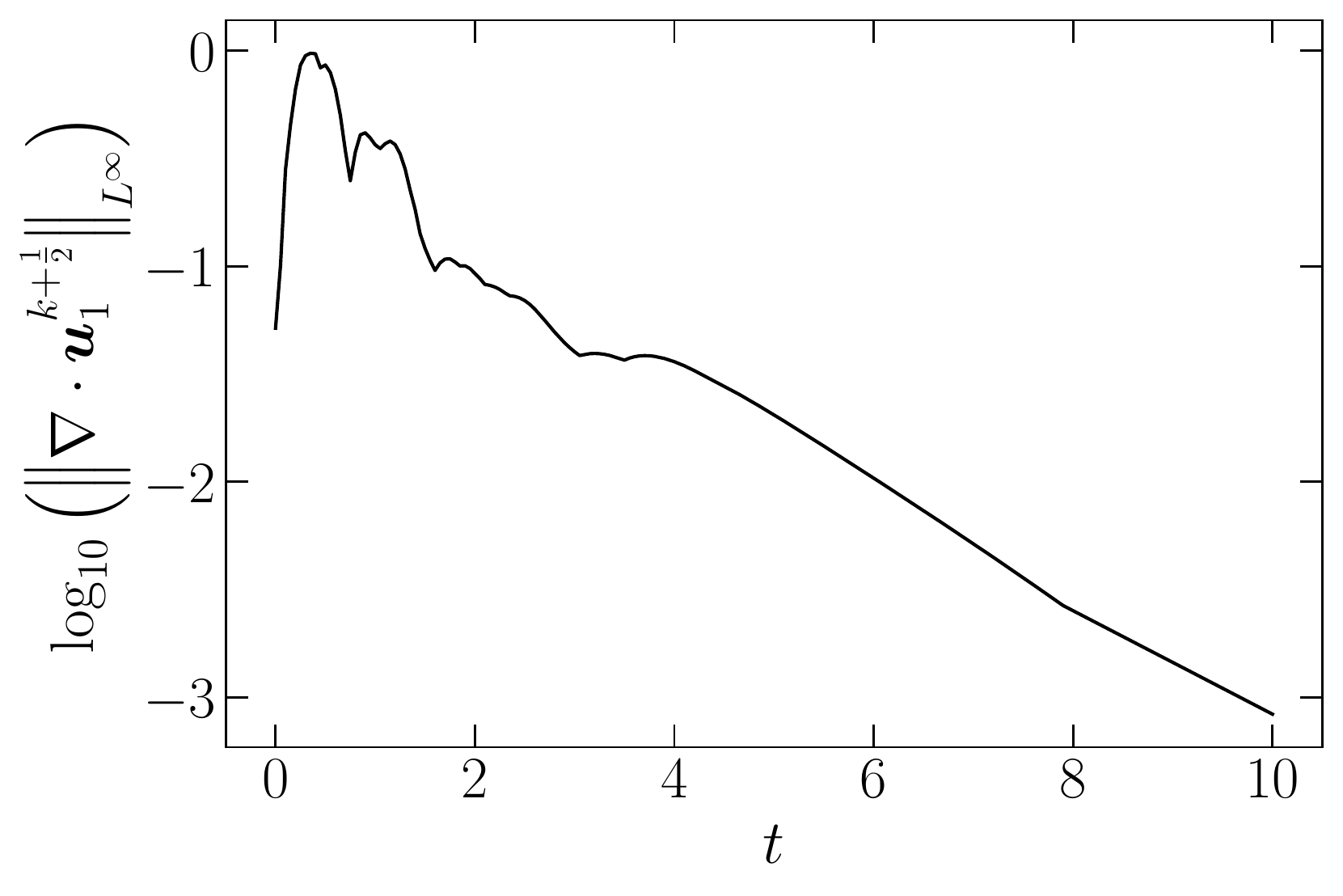}
			\end{minipage}
		}
		\caption{\REVO{Some results of the magnitude of $ \left\|\nabla\cdot\bu^{k+\frac{1}{2}}_{1}\right\|_{L^\infty} $ for both the conservation and dissipation tests at $ h=1/3 $, $ N=2 $ and $ \varDelta t=1/20 $. 
				The divergence is computed per element.
			}}
		\label{fig: conservation test weak conservation}}
\end{figure}

\REVO{In Fig.~\ref{fig: conservation test weak conservation}, some results of the magnitude of $ \left\|\nabla\cdot\bu^{k+\frac{1}{2}}_{1}\right\|_{L^\infty} $ are presented. It is seen that for both the convergence and dissipation tests the conservation of mass is not satisfied for $ \bu^{k+\frac{1}{2}}_{1} $. It is not surprising that the error is large especially for the inviscid case as we have used an extremely coarse mesh. This is consistent with the analysis that the constraint of mass conservation is only weakly imposed for $\bu^{k+\frac{1}{2}}_{1} \in C(\Omega)\subset H(\mathrm{curl};\Omega) $, see Section~\ref{Sec: mass conservation semi discrete}. Also see the analysis at the continuous level in Section~\ref{Subsub: continuous mass conservation}.}

Note that these tests are also valid for the non-zero conservative external body force. If $ \varphi $ is known and $ \boldsymbol{f}=\nabla\varphi\neq \boldsymbol{0} $, we can still first conduct the test with $ \boldsymbol{f}=\boldsymbol{0} $ and get the same results. The only difference is that we now obtain the solution for the extended total pressure $ P^{\prime} $, see \eqref{Eq: extended total pressure}. We can post-process $ P^{\prime} $ with the known $ \varphi $ to retrieve the solution for total pressure $ P $.

\subsubsection{Convergence tests} 
We now investigate whether the proposed method produces  converging solutions and, if yes, what is the convergence rate of the proposed method with a manufactured solution. Assume 
\[\boldsymbol{u}=\left[(2-t)\cos(2\pi z),\ (1+t)\sin(2\pi z),\ (1-t)\sin(2\pi x)\right] ^{\mathsf{T}}\]
and
\[p=\sin(2\pi(x+y+t))\]
solve the Navier-Stokes equations with Reynolds number $ \Rn=1 $ and the body force $ \boldsymbol{f} $ which can be calculated from $ \boldsymbol{u} $, $ p $ and $ \Rn $ using the Navier-Stokes equations. The exact solutions of vorticity $ \boldsymbol{\omega} $ and total pressure $ P $ can also be calculated. We use $ \left. \boldsymbol{u}\right| _{t=0} $ as initial condition and let the flow evolve for different mesh element sizes and polynomial space degrees. Errors are then measured at $ t = 2 $.

\begin{figure}[h!]
	\centering{
		\subfloat{
			\begin{minipage}[b]{0.5\textwidth}
				\centering
				\includegraphics[width=0.85\linewidth]{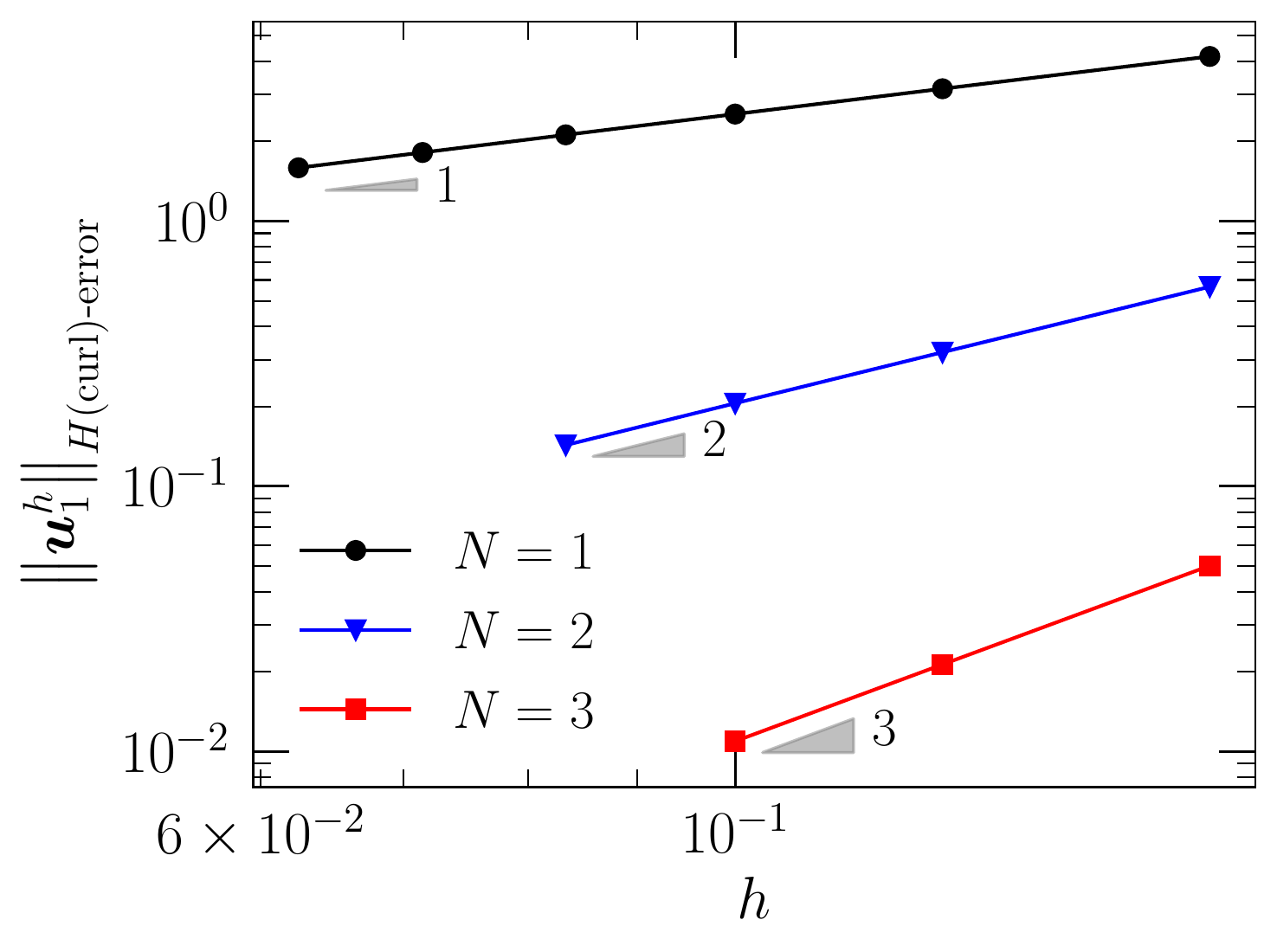}
			\end{minipage}
		}
		\subfloat{
			\begin{minipage}[b]{0.5\textwidth}
				\centering
				\includegraphics[width=0.85\linewidth]{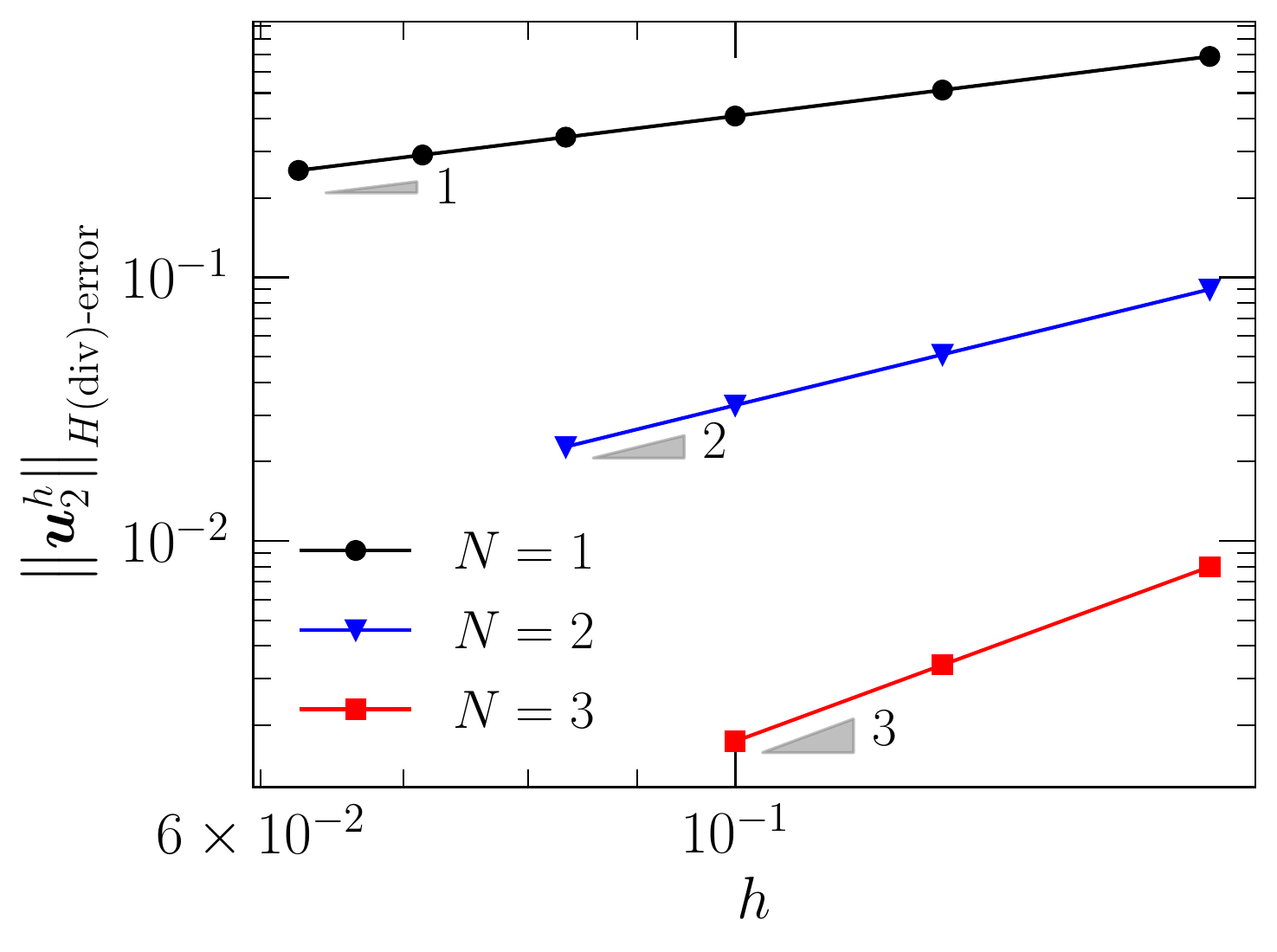}
			\end{minipage}
		}\\
		\subfloat{
		\begin{minipage}[b]{0.5\textwidth}
			\centering
			\includegraphics[width=0.85\linewidth]{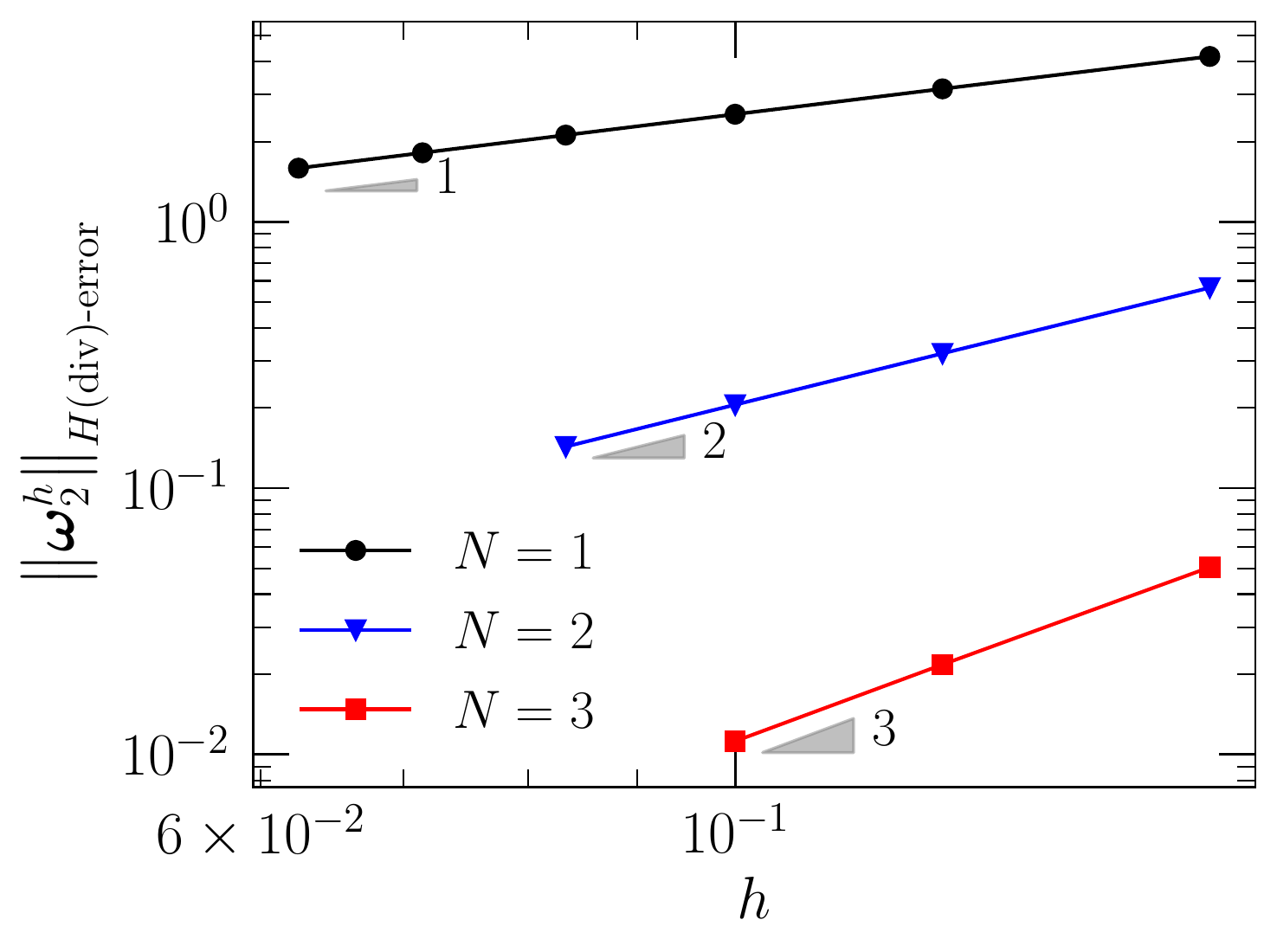}
		\end{minipage}
		}
		\subfloat{
		\begin{minipage}[b]{0.5\textwidth}
			\centering
			\includegraphics[width=0.85\linewidth]{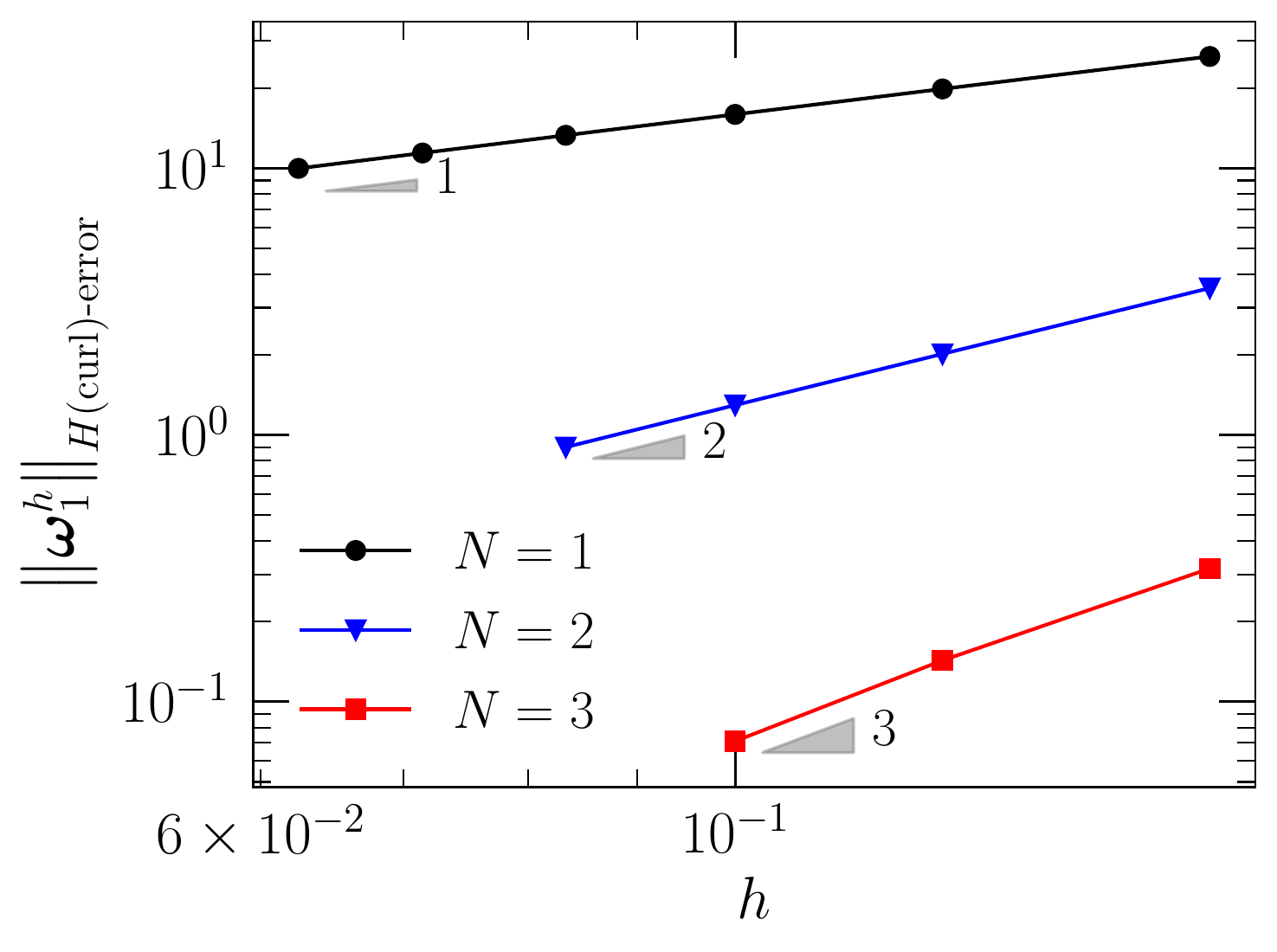}
		\end{minipage}
		}\\
		\subfloat{
			\begin{minipage}[b]{0.5\textwidth}
				\centering
				\includegraphics[width=0.85\linewidth]{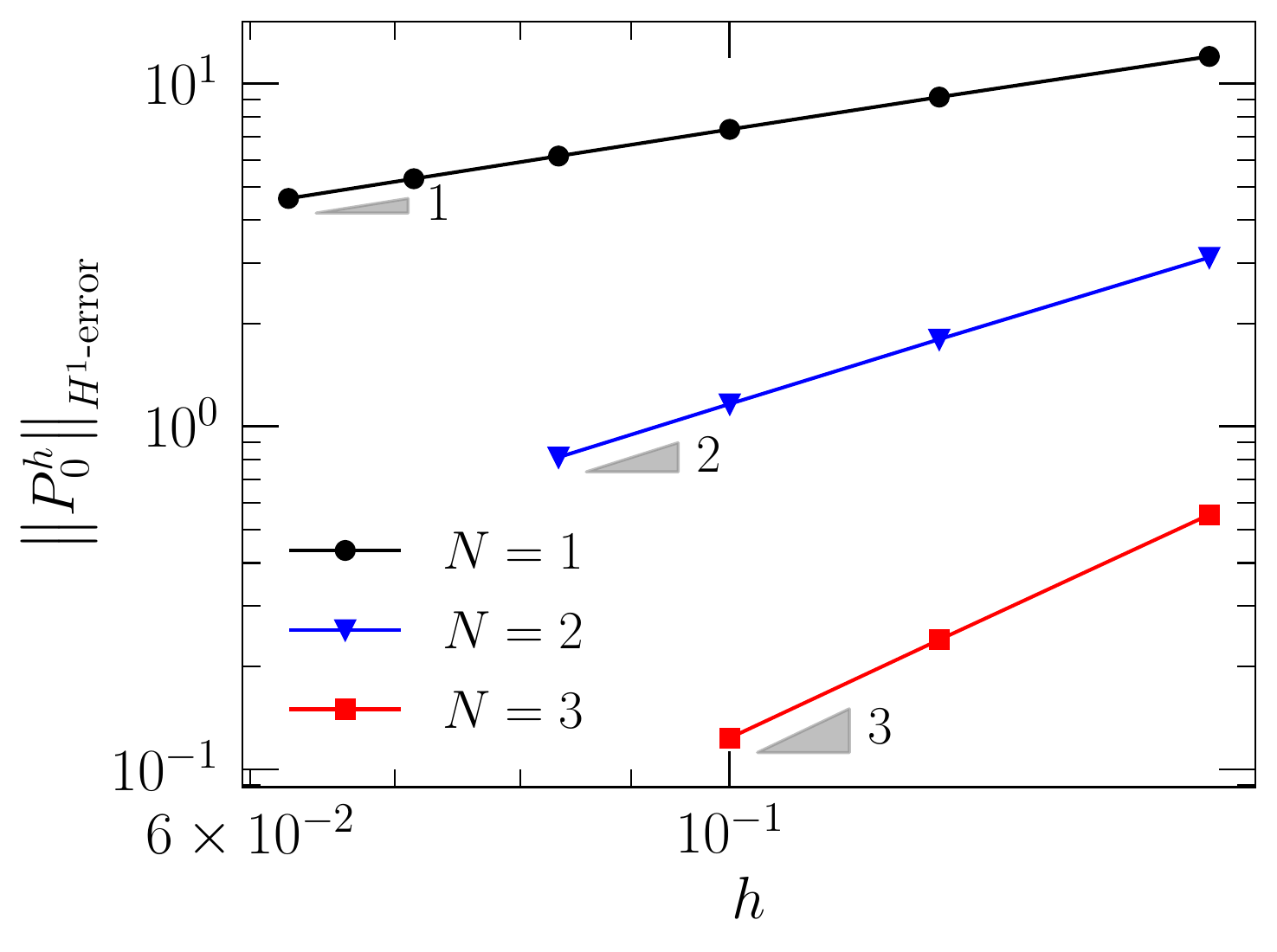}
			\end{minipage}
		}
		\subfloat{
			\begin{minipage}[b]{0.5\textwidth}
				\centering
				\includegraphics[width=0.85\linewidth]{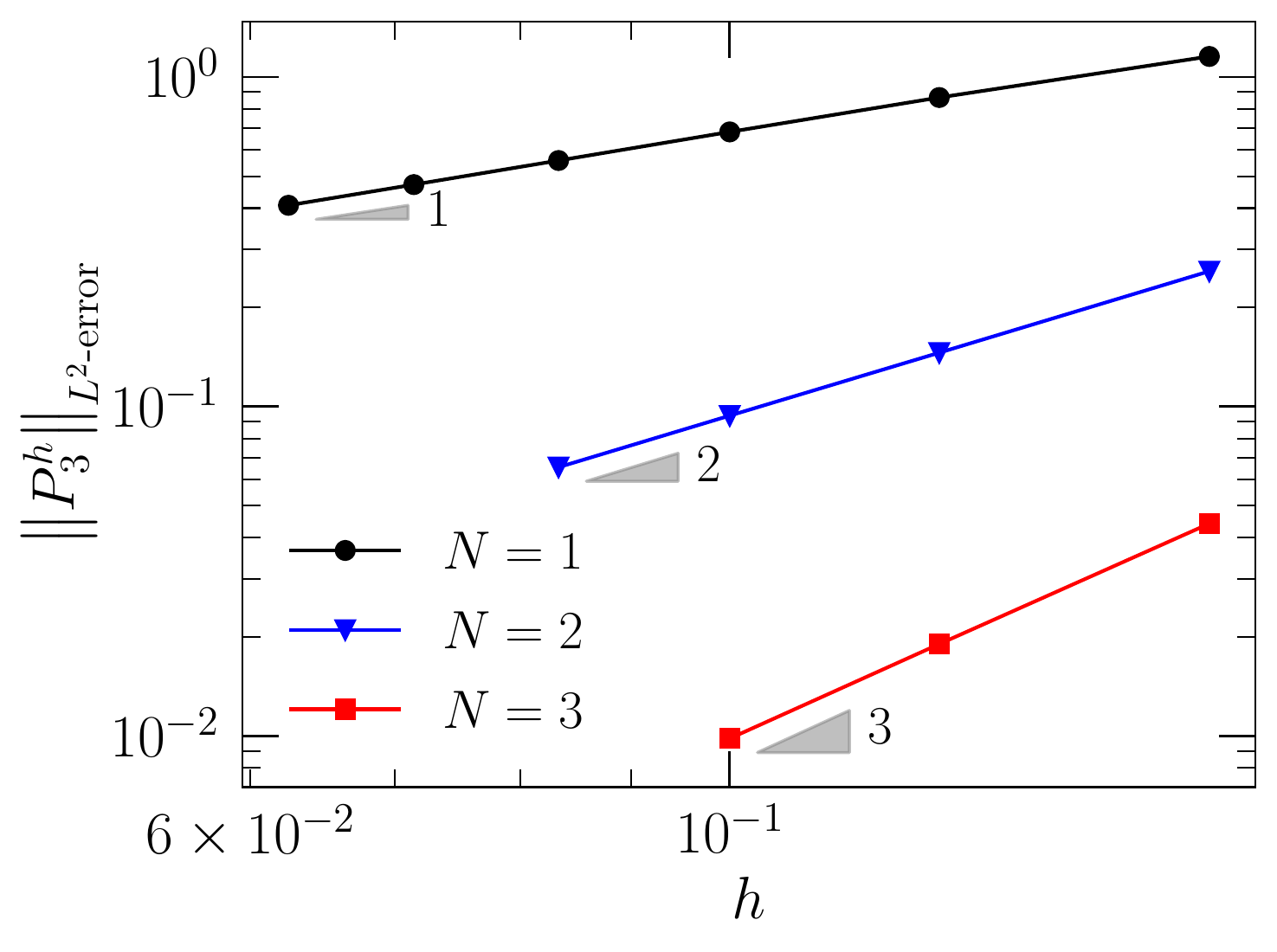}
			\end{minipage}
		}\\
		\subfloat{
		\begin{minipage}[b]{0.5\textwidth}
			\centering
			\includegraphics[width=0.85\linewidth]{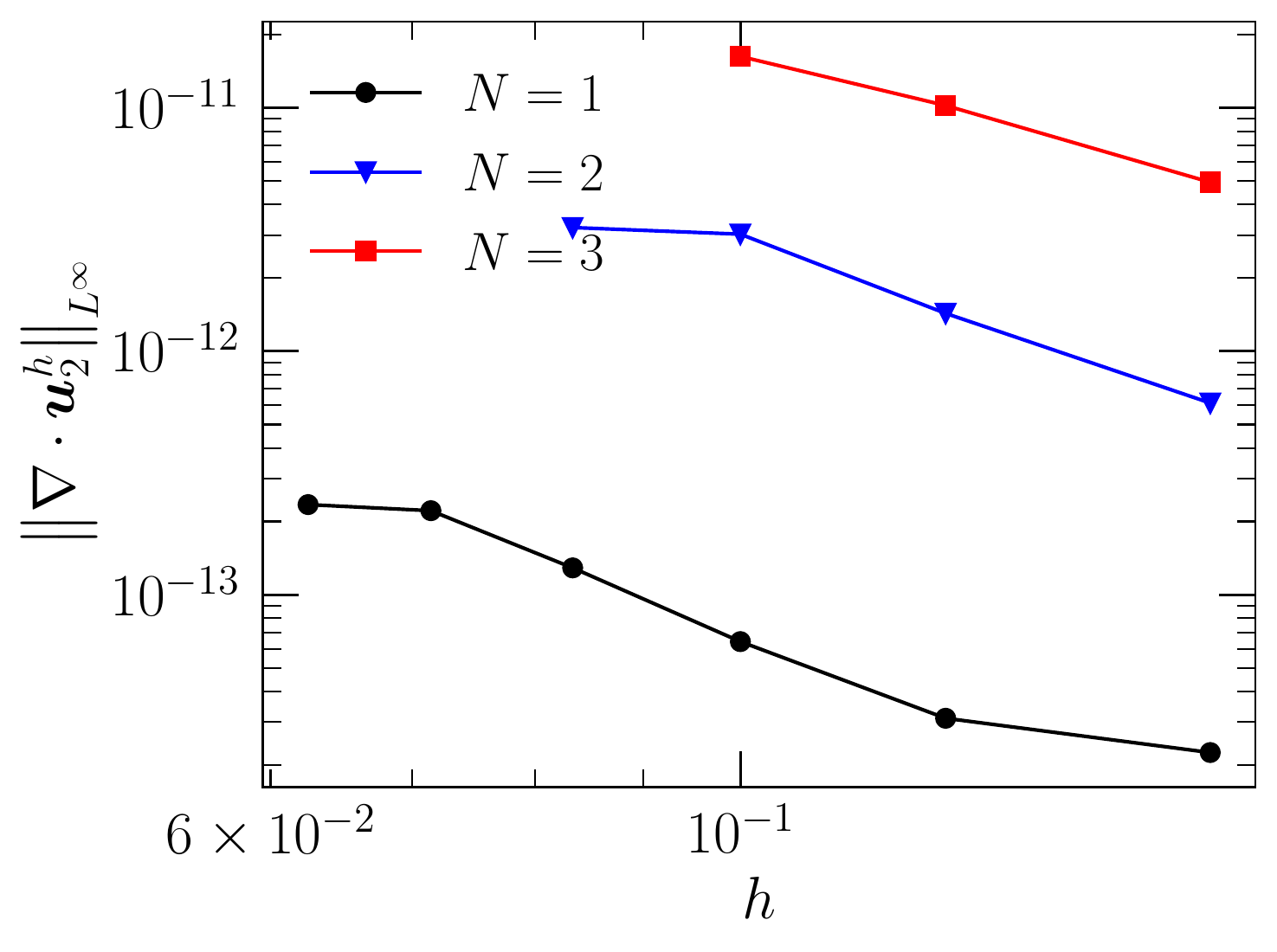}
		\end{minipage}
		}
		\caption{$ ph $-convergence and mass conservation results of the convergence tests. $ \Delta t = 1/50 $.}
		\label{fig: convergence test}
	}
\end{figure}

Results are presented in Fig.~\ref{fig: convergence test} where the optimal convergence rates are observed for all variables of the dual-field formulation when the mesh is \MOD{$ h $-refined} under different polynomial degrees. The plot of $ \left\|\nabla\cdot\bu_{2}^{h}\right\|_{L^\infty} $ shows that the pointwise conservation of mass is satisfied up to the machine precision in all cases.

\begin{figure}[h!]
	\centering{
		\subfloat{
			\begin{minipage}[b]{0.5\textwidth}
				\centering
				\includegraphics[width=0.85\linewidth]{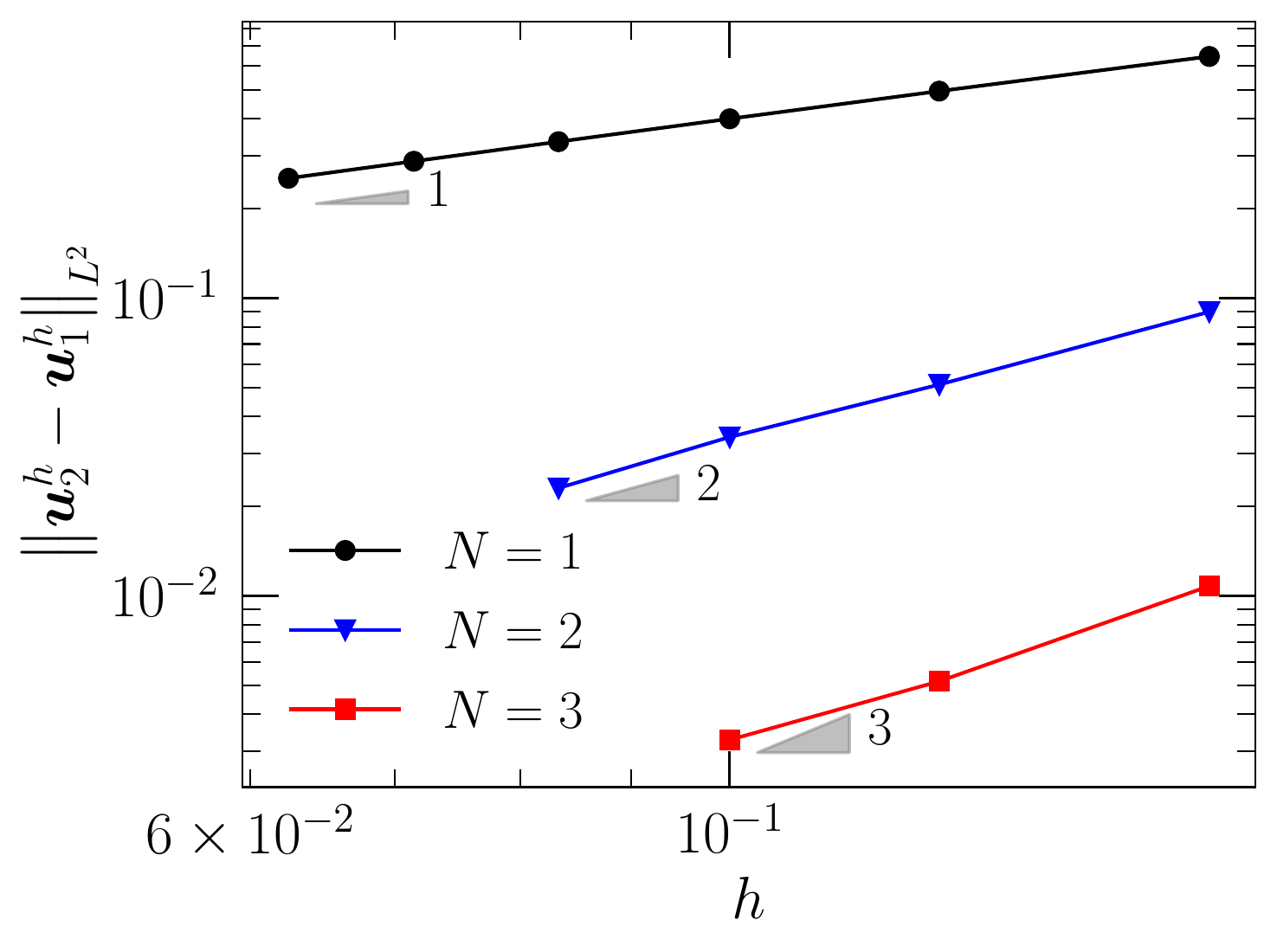}
			\end{minipage}
		}
		\subfloat{
			\begin{minipage}[b]{0.5\textwidth}
				\centering
				\includegraphics[width=0.85\linewidth]{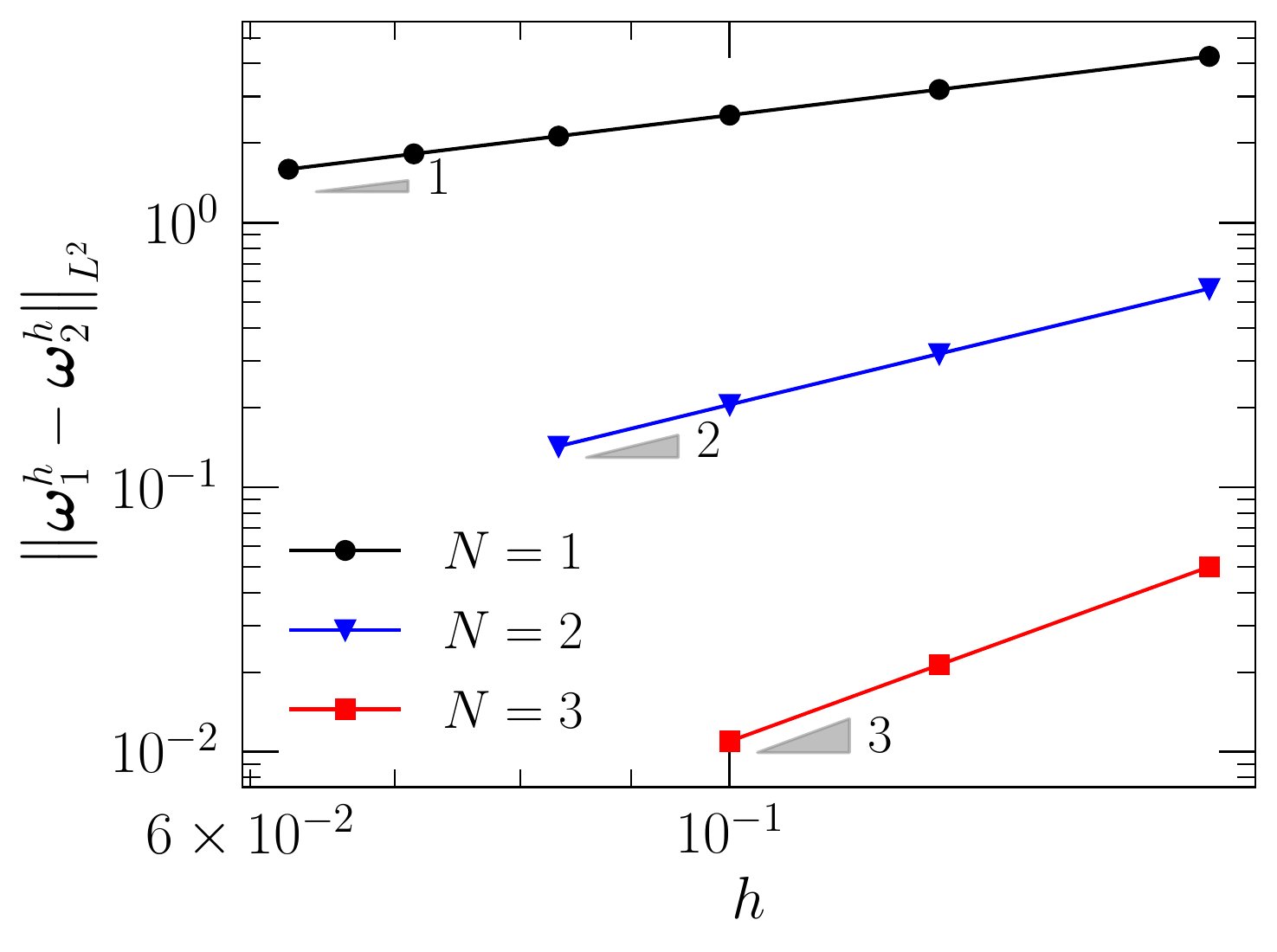}
			\end{minipage}
		}\\
		\caption{$ ph $-convergence results of the $ L^2 $-differences between the dual representations of the solution of velocity (Left) and vorticity (Right) for the convergence tests. $ \Delta t = 1/50 $.}
		\label{fig: convergence test u1u2 w1w2 diff}
	}
\end{figure}

We can also measure the difference between the two dual solutions of one physical variable. The results of $ \left\|\boldsymbol{u}_{2}^{h}-\boldsymbol{u}_{1}^{h}\right\|_{L^2} $ and $ \left\|\boldsymbol{\omega}_{2}^{h}-\boldsymbol{\omega}_{1}^{h}\right\|_{L^2} $ \REVO{at $ t=2 $} are shown in Fig.~\ref{fig: convergence test u1u2 w1w2 diff}. \REVO{Note that, since $ \bu^h_{1} $ and $ \bu^h_{2} $ ($ \bw^h_{1} $ and $ \bw^h_{2} $) are staggered in time, we have used the midpoint rule, \eqref{Eq: mpr}, to $ \bu^{h}_{1} $ ($ \bw^{h}_{1} $) such that it can compared to $ \bu^{h}_{2} $ at $ t=2 $, an integer time instant.} It is not surprising that they converge under $ p $- or $ h $-refinement. This suggests that we can use them for \MOD{accuracy} indicators, for example,
\[
	\dfrac{\left\|\boldsymbol{u}_{2}^{h}-\boldsymbol{u}_{1}^{h}\right\|_{L^2}}{\left\|\boldsymbol{u}_{1}^{h}\right\|_{L^2}}\,,\quad\dfrac{\left\|\boldsymbol{u}_{2}^{h}-\boldsymbol{u}_{1}^{h}\right\|_{L^2}}{\left\|\boldsymbol{u}_{2}^{h}\right\|_{L^2}}\quad \text{or}\quad \dfrac{2\left\|\boldsymbol{u}_{2}^{h}-\boldsymbol{u}_{1}^{h}\right\|_{L^2}}{\left\|\boldsymbol{u}_{1}^{h}+\boldsymbol{u}_{2}^{h}\right\|_{L^2}}\,,
\]
which can be very helpful for general (non-manufactured) simulations. More interestingly, one can measure the local difference of the dual solutions and use it as an indicator for mesh adaptivity, which is outside of the scope of the current paper. From this aspect, the existence of dual representations of the solution for one variable can be regarded as an advantage for the proposed method. \REVO{Despite the existence of the difference between the dual representations, both of them should be considered as equally important solutions of the variable. Recall the dual character of the velocity field which is hard to capture in one discrete space, see Section~\ref{SEC: Objective}. The dual representations together can be regarded as a discretization of its dual character.}

\subsection{Taylor-Green vortex} \label{Sub: TESTS TGV}
We now test the method with a more general flow, the Taylor-Green vortex (TGV) flow. The domain is given as $ \Omega:=[-\pi, \pi ]^3 $ and is periodic. $ V=8\pi^3 $ denotes the volume of the domain. The body force is set to $ \boldsymbol{f}=\boldsymbol{0} $ and the initial condition is selected to be
\[\left.\bu\right|_{t=0}=\left[
\sin(x)\cos(y)\cos(z),
-\cos(x)\sin(y)\cos(z),
0
\right]^{\mathsf{T}}.\]
Such an initial condition possesses \MOD{kinetic} energy $ \left.\mathcal{K}\right|_{t=0} =0.125$ and zero helicity. We solve the flow using the proposed mimetic dual-field method at $ \Rn=500 $. 

\begin{figure}[h!]
	\centering{
		\subfloat[$ t=0 $]{
			\begin{minipage}[b]{0.5\textwidth}
				\centering
				\includegraphics[width=0.89\linewidth]{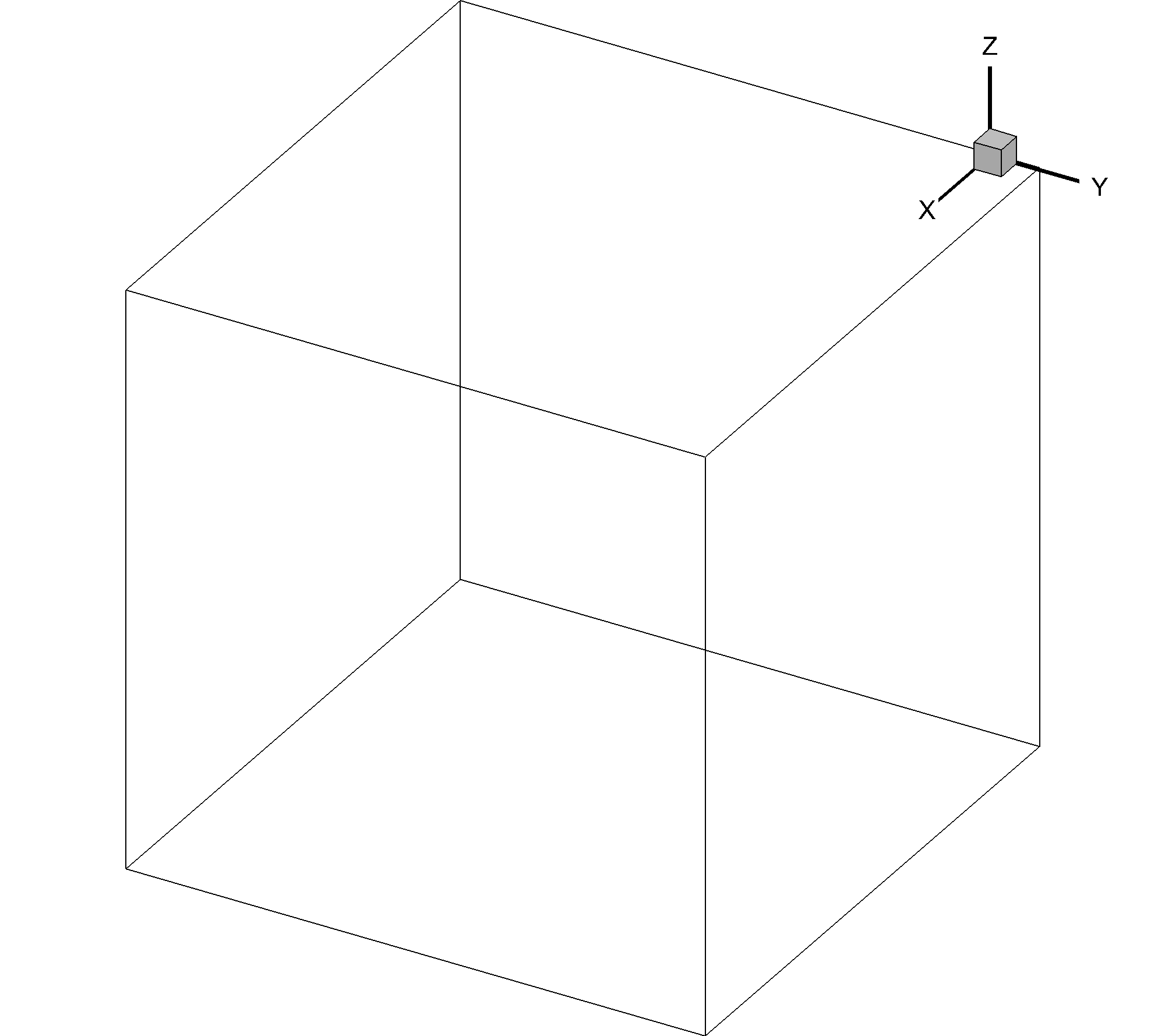}
			\end{minipage}
		}
		\subfloat[$ t=3 $]{
			\begin{minipage}[b]{0.5\textwidth}
				\centering
				\includegraphics[width=0.89\linewidth]{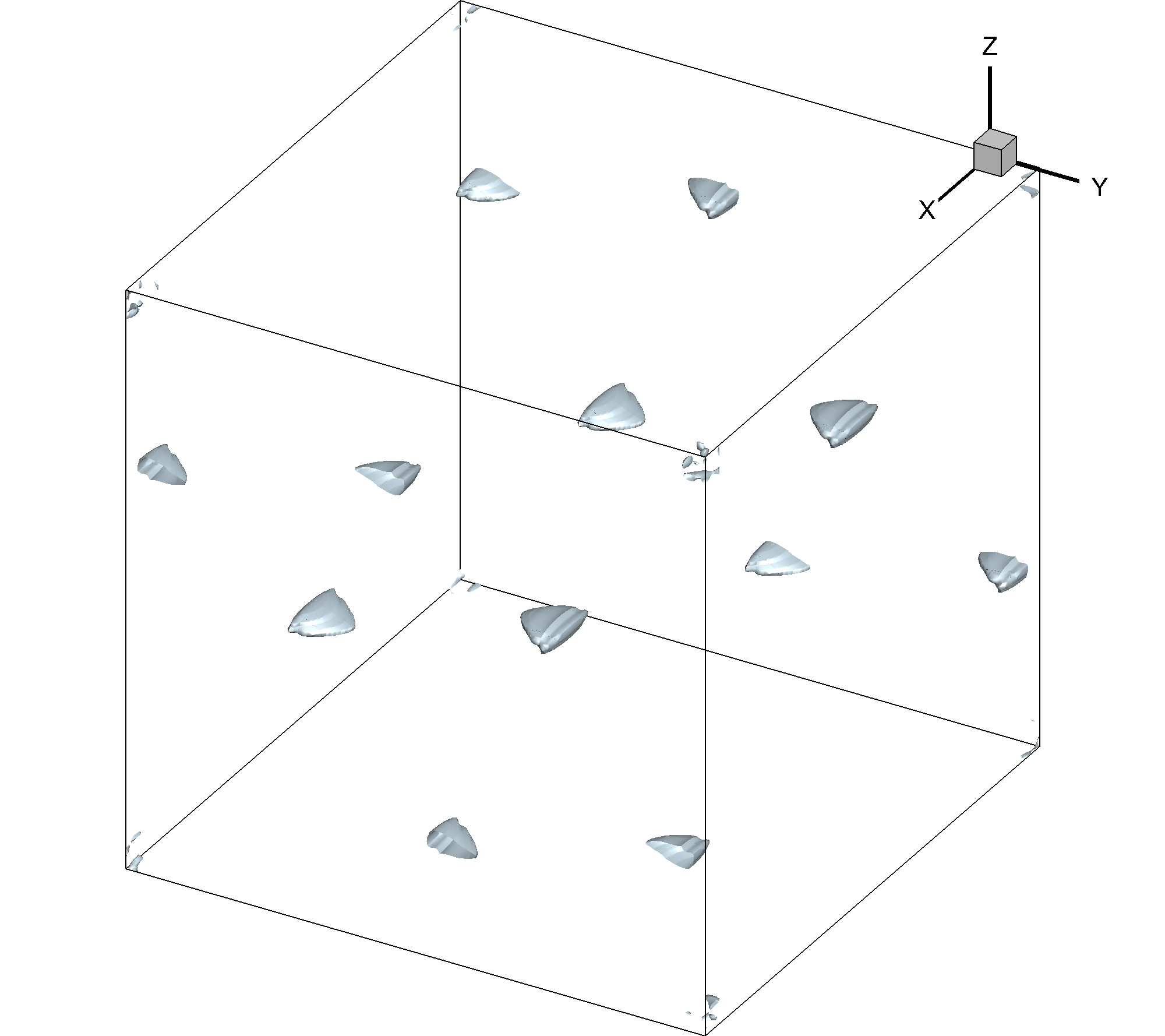}
			\end{minipage}
		}\\
		\subfloat[$ t=6 $]{
			\begin{minipage}[b]{0.5\textwidth}
				\centering
				\includegraphics[width=0.89\linewidth]{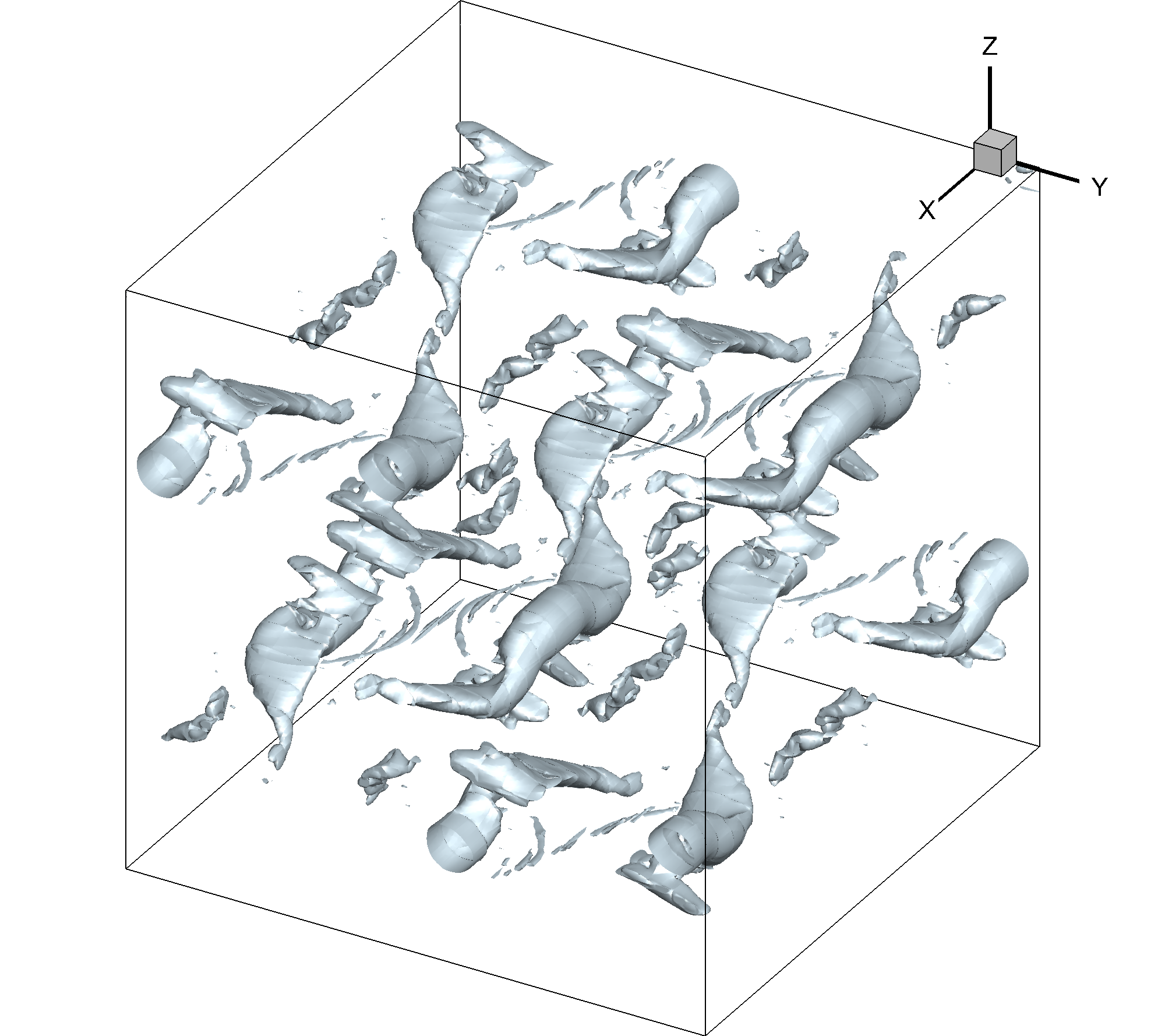}
			\end{minipage}
		}
		\subfloat[$ t=9 $]{
			\begin{minipage}[b]{0.5\textwidth}
				\centering
				\includegraphics[width=0.89\linewidth]{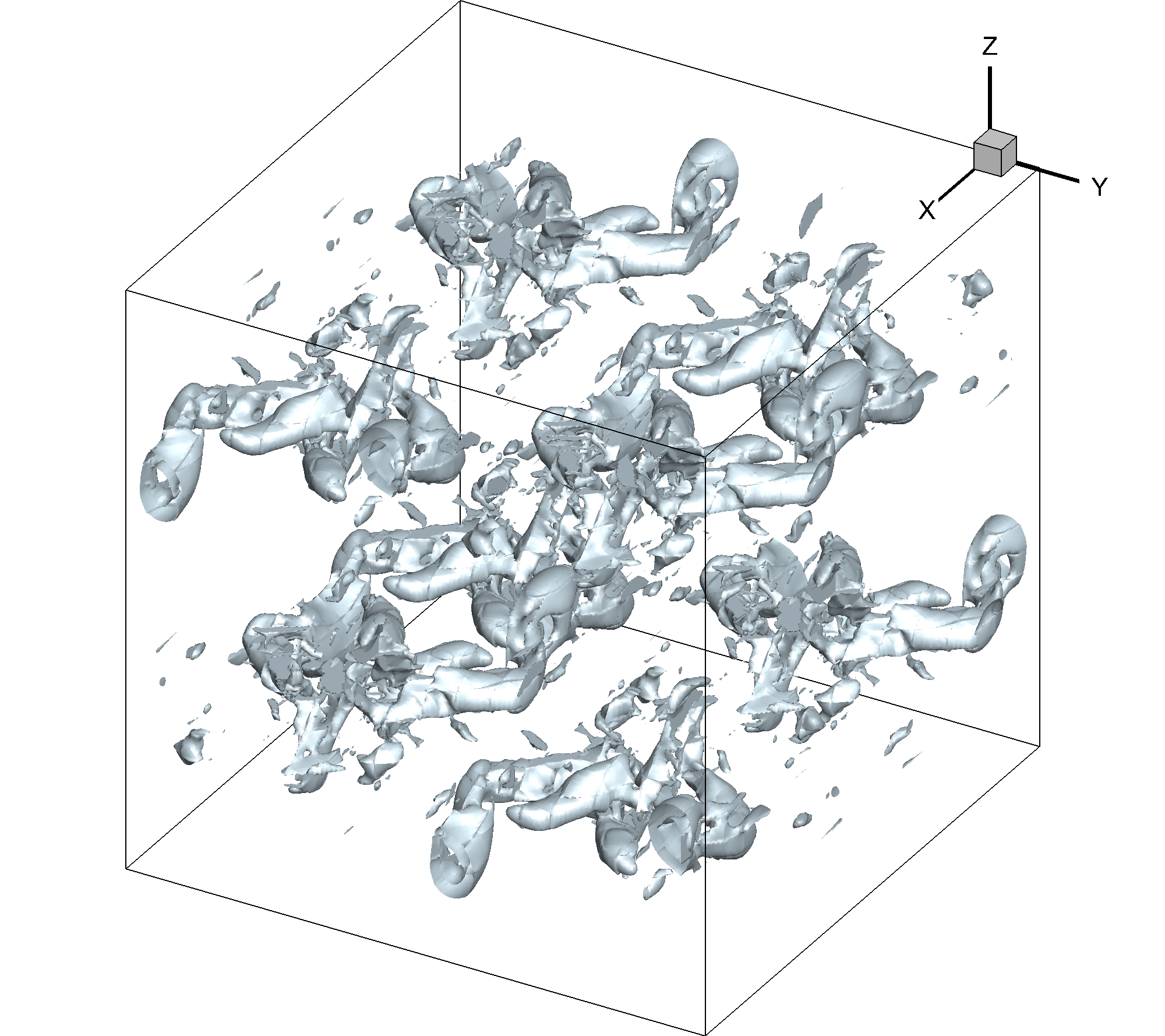}
			\end{minipage}
		}\\
		\subfloat[$ t=12 $]{
			\begin{minipage}[b]{0.5\textwidth}
				\centering
				\includegraphics[width=0.89\linewidth]{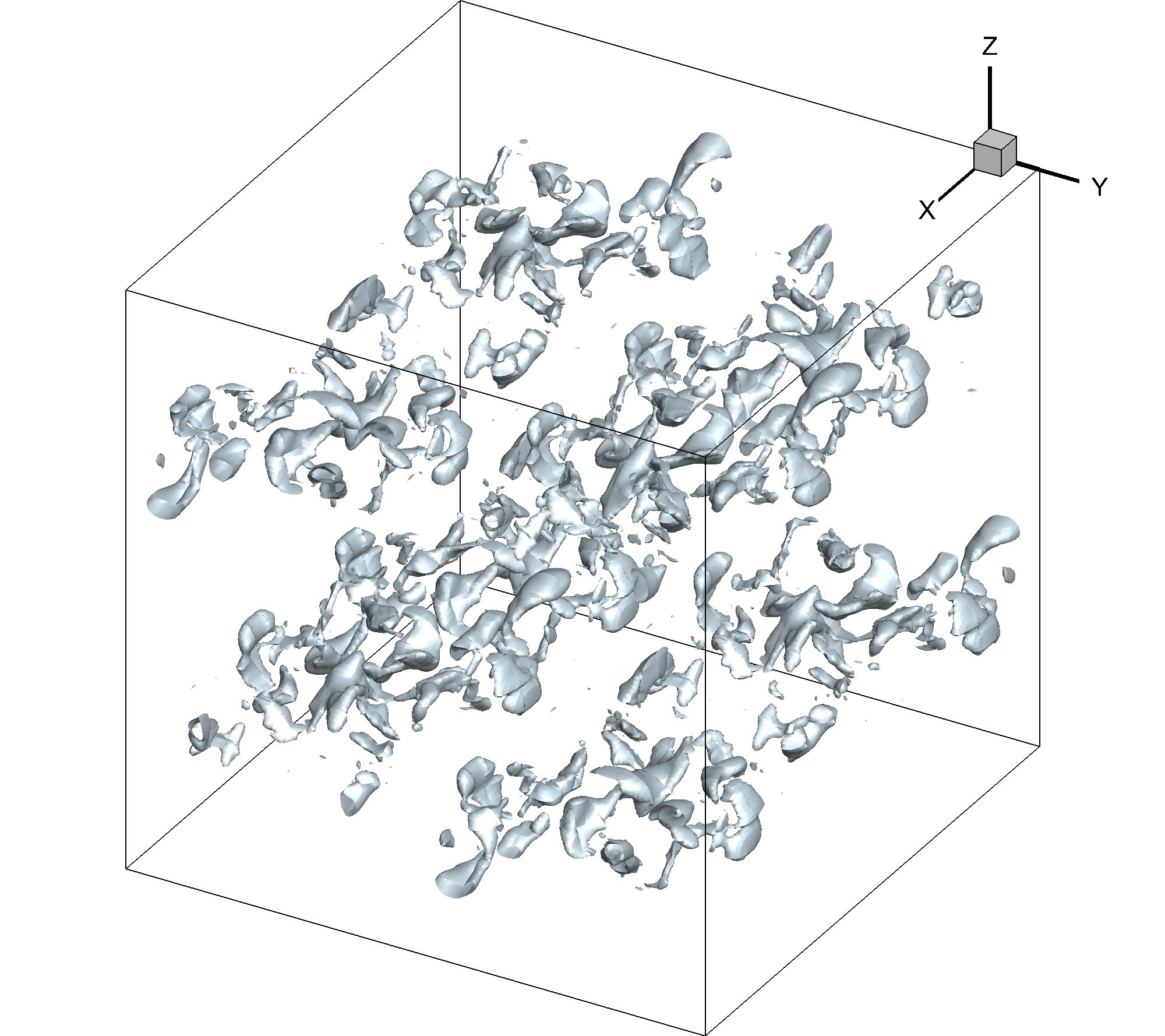}
			\end{minipage}
		}
		\subfloat[$ t=15 $]{
			\begin{minipage}[b]{0.5\textwidth}
				\centering
				\includegraphics[width=0.89\linewidth]{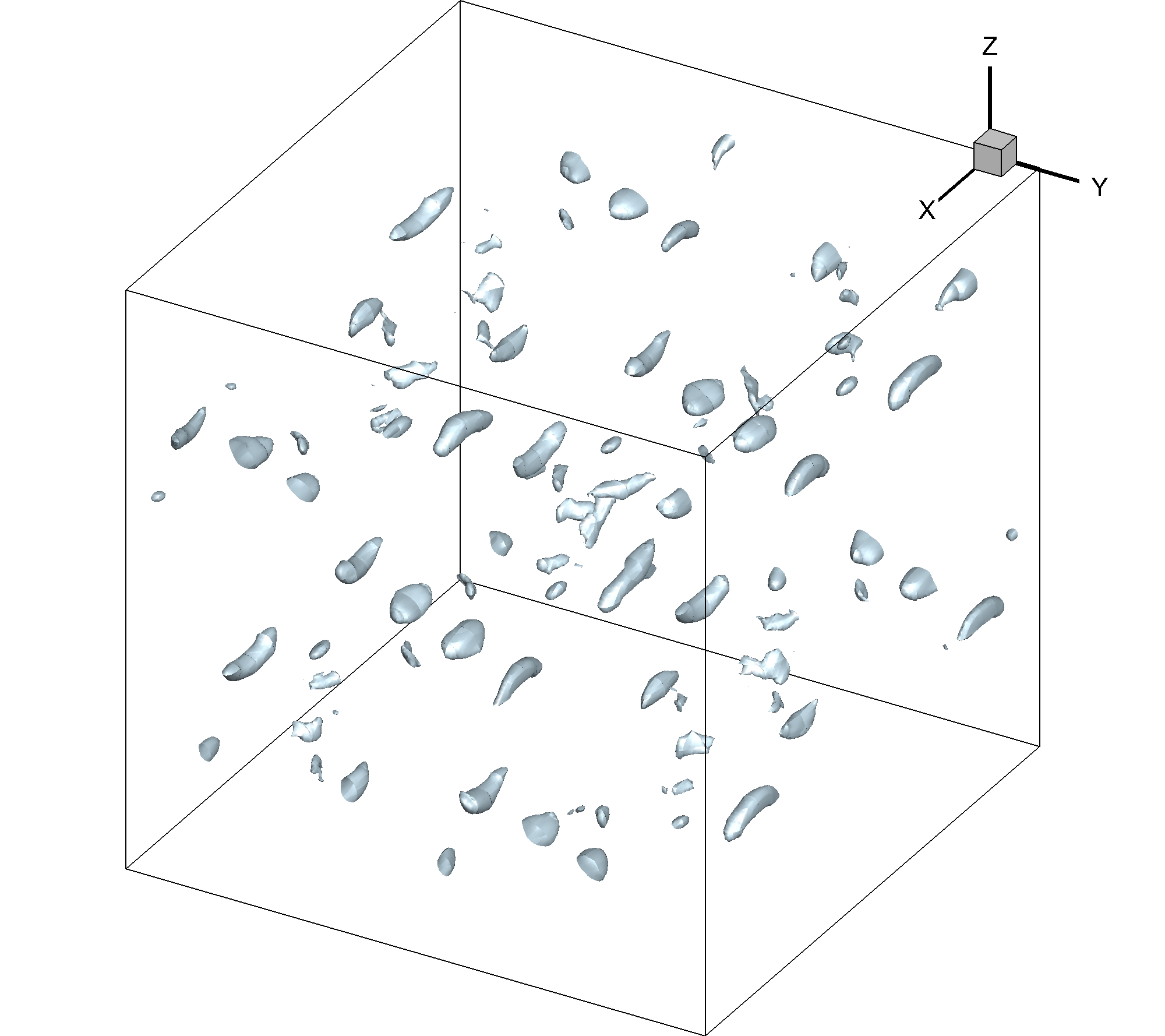}
			\end{minipage}
		}
		\caption{Iso-surface of $ \omega_{1}^{x}=-3 $ $\left( \bw_{1}^{h}=\left( \omega_{1}^{x},\omega_{1}^{y},\omega_{1}^{z}\right) \right)  $ for the TGV test using MDF-$ 24p3 $.  MDF-$ 24p3 $ stands for the mimetic dual-field method at $ h=1/24 $ using polynomial spaces of degree 3. The time step interval is selected to be $ \varDelta t =1/50 $.}
		\label{fig: TGV iso surface wx m3}}
\end{figure}

Iso-surfaces of $ \omega_{1}^{x}=-3 $ $\left( \bw_{1}^{h} =\left(  \omega_{1}^{x}, \omega_{1}^{y}, \omega_{1}^{z}\right)  \right) $ at some time instances are shown in Fig.~\ref{fig: TGV iso surface wx m3}. It is seen that the flow initially induces vortices of clear structures which then break down and finally are dissipated by the \MOD{viscosity}. 

\begin{figure}[h!]
	\centering{
		\subfloat{
			\begin{minipage}[b]{0.5\textwidth}
				\centering
				\includegraphics[width=0.78\linewidth]{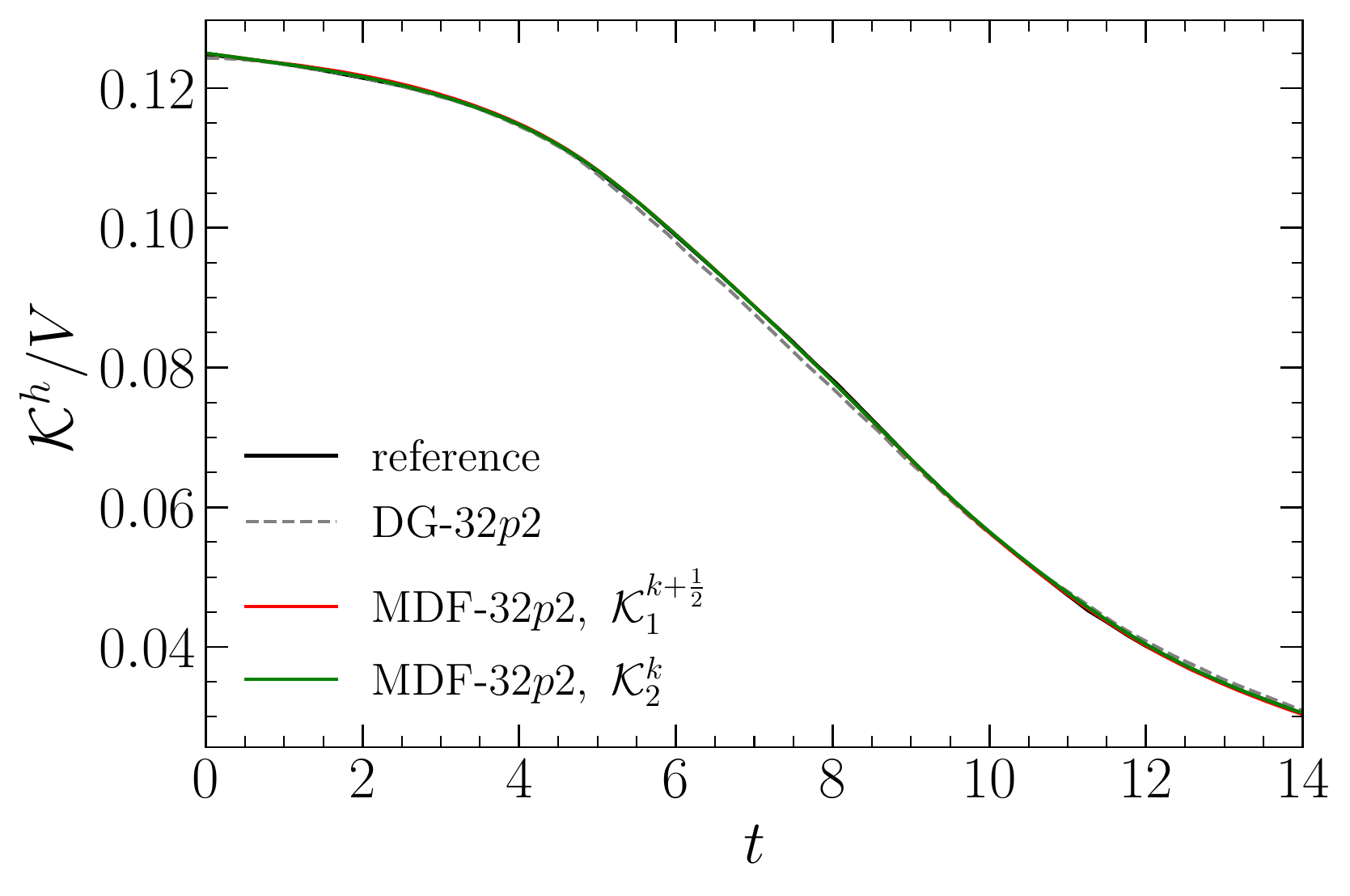}
			\end{minipage}
		}
		\subfloat{
			\begin{minipage}[b]{0.5\textwidth}
				\centering
				\includegraphics[width=0.77\linewidth]{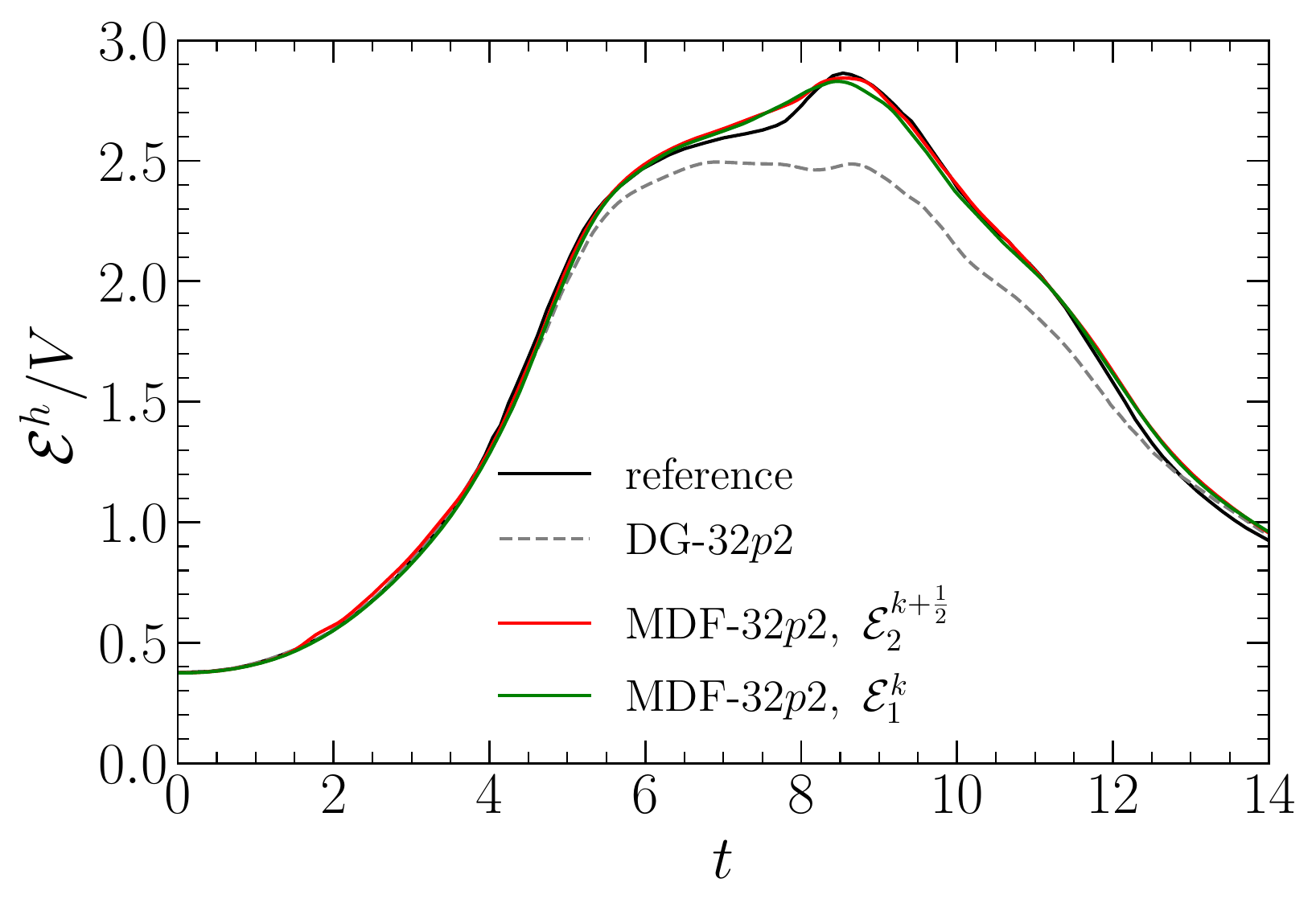}
			\end{minipage}
		}\\
		\subfloat{
			\begin{minipage}[b]{0.5\textwidth}
				\centering
				\includegraphics[width=0.76\linewidth]{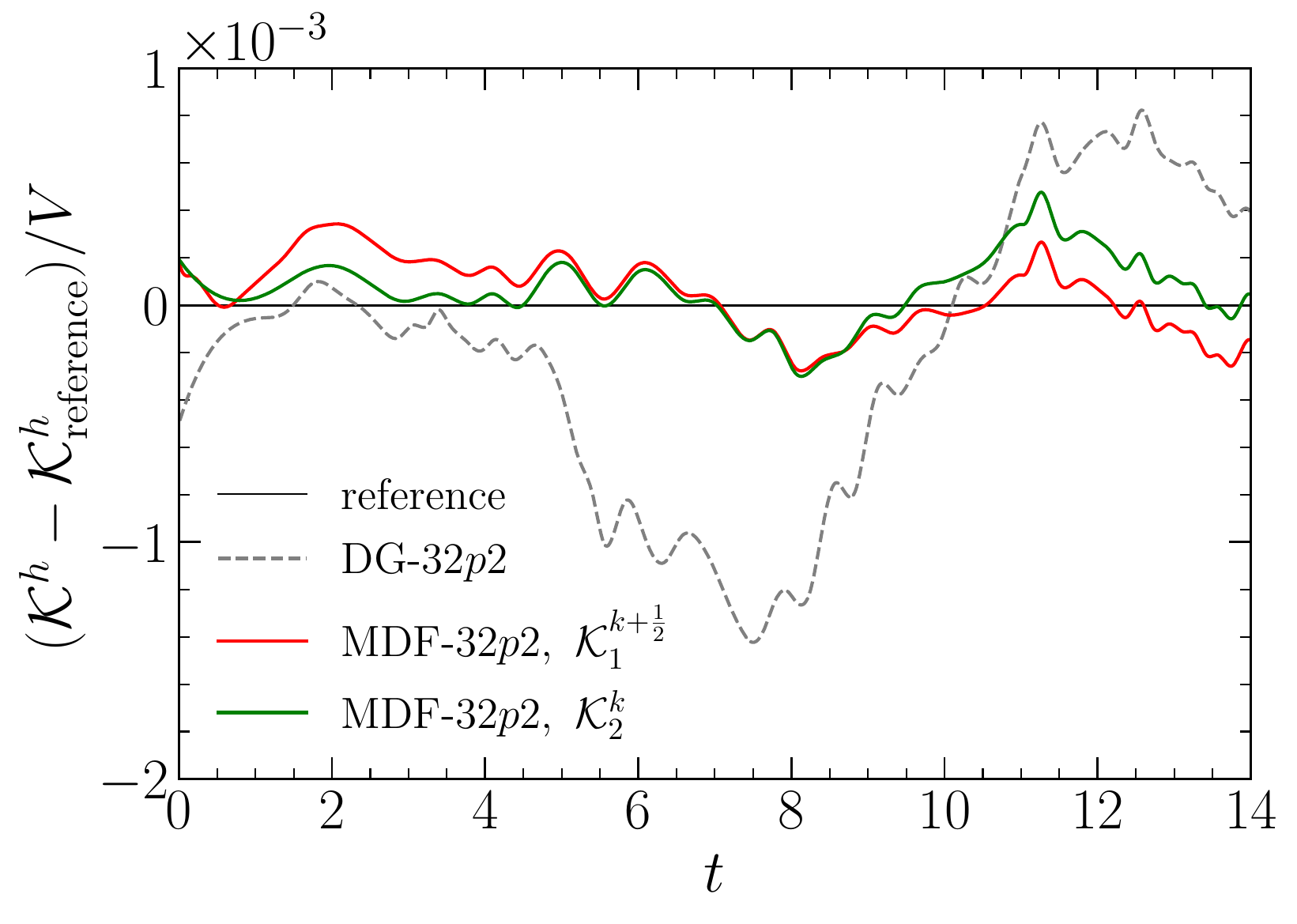}
			\end{minipage}
		}
		\subfloat{
			\begin{minipage}[b]{0.5\textwidth}
				\centering
				\includegraphics[width=0.77\linewidth]{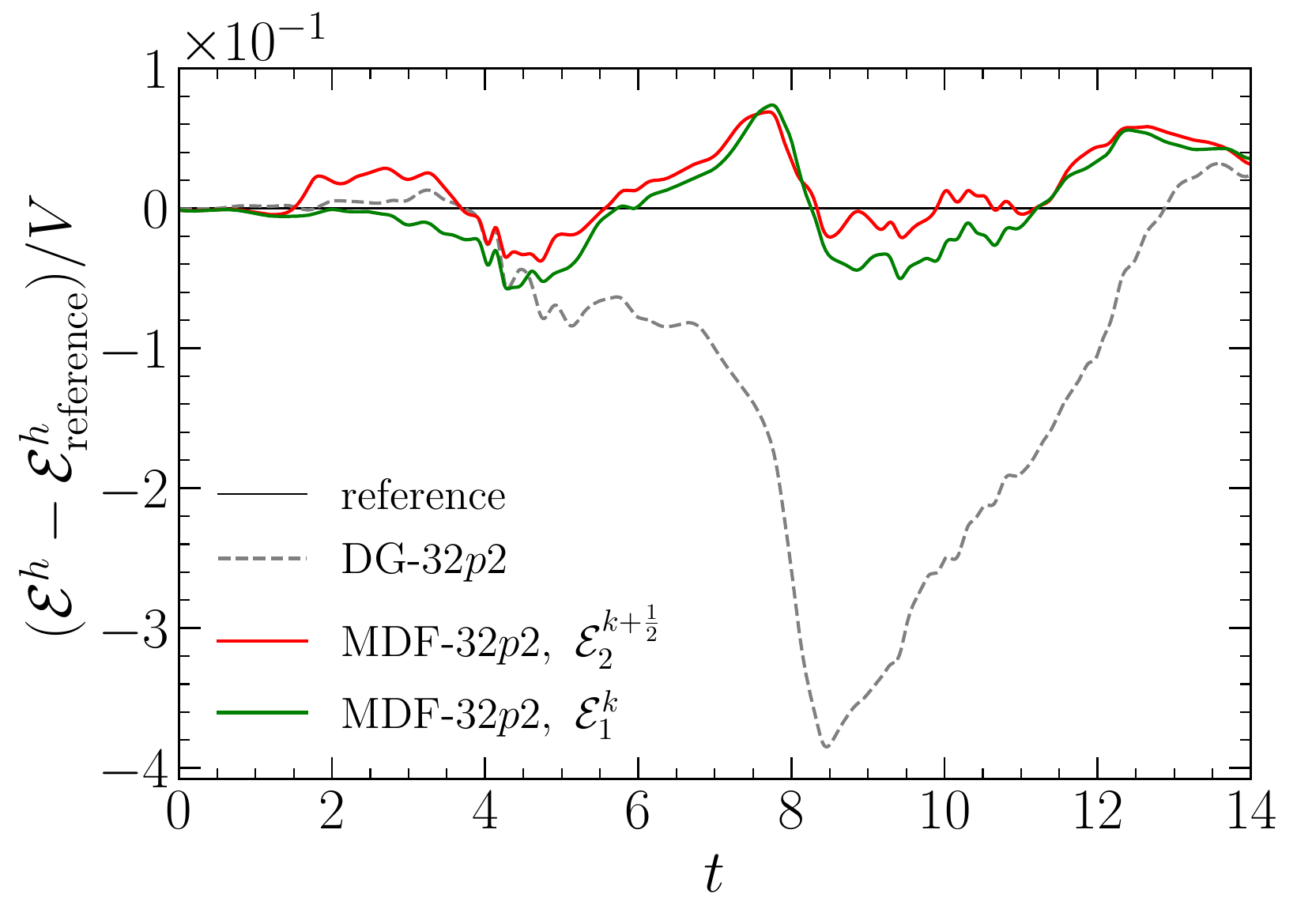}
			\end{minipage}
		}
		\caption{\REVT{A comparison of kinetic energy and enstrophy results of the TGV test. The reference and the DG-$ 32p2 $ results are taken from \cite{chapelier2012inviscid}. The reference method is a $ 128 $th order spectral method. DG stands for a discontinuous Galerkin method. MDF stands for the mimetic dual-field method. $ 32p2 $ represents $ h=1/32 $ and the degree of the polynomial spaces is 2. The time step interval is selected to be $ \varDelta t =1/50 $ for MDF-$32p2 $.}}
		\label{fig: TGV 1}}
\end{figure}

\begin{figure}[h!]
	\centering{
		\subfloat{
			\begin{minipage}[b]{0.5\textwidth}
				\centering
				\includegraphics[width=0.78\linewidth]{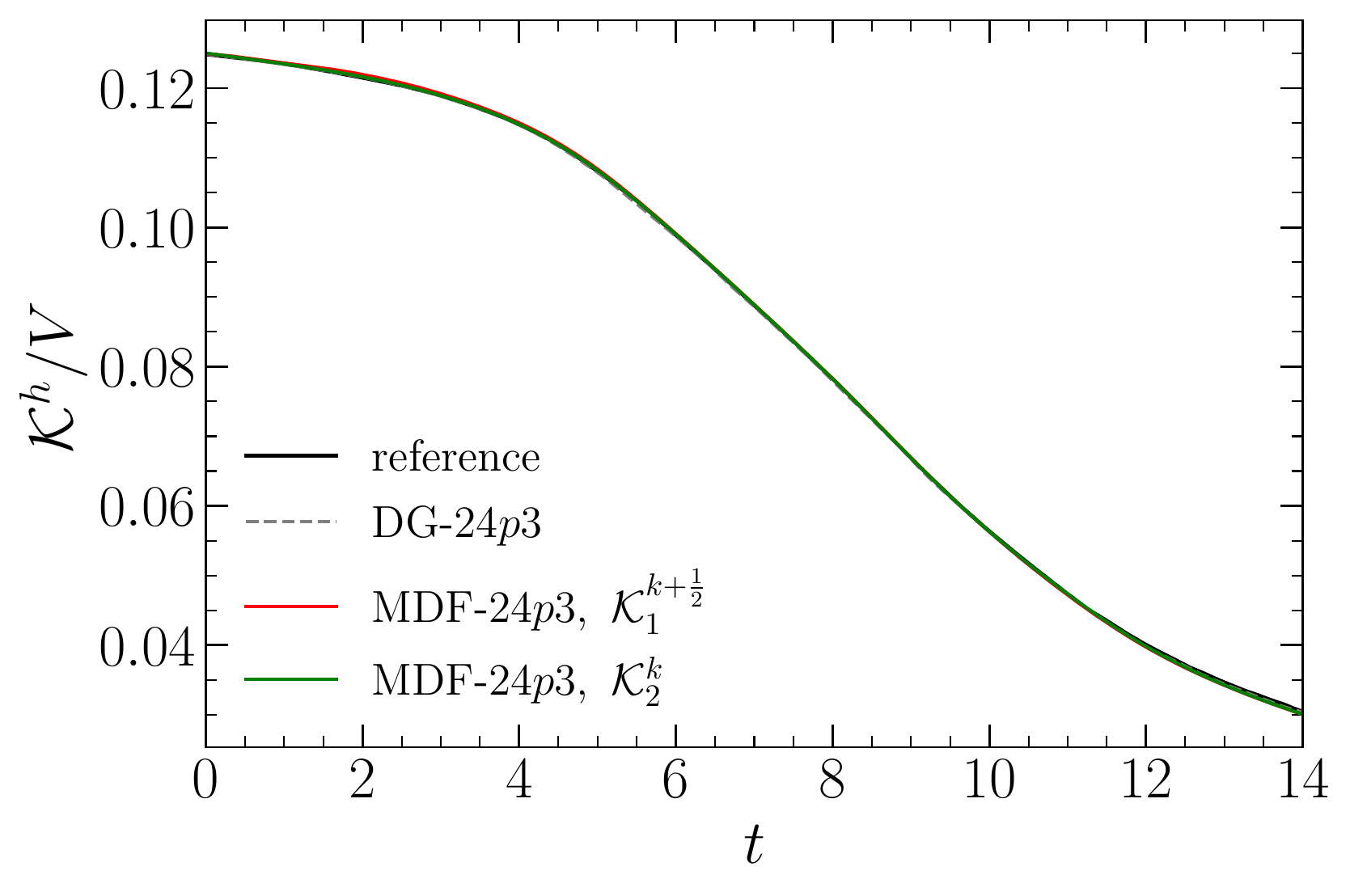}
			\end{minipage}
		}
		\subfloat{
			\begin{minipage}[b]{0.5\textwidth}
				\centering
				\includegraphics[width=0.77\linewidth]{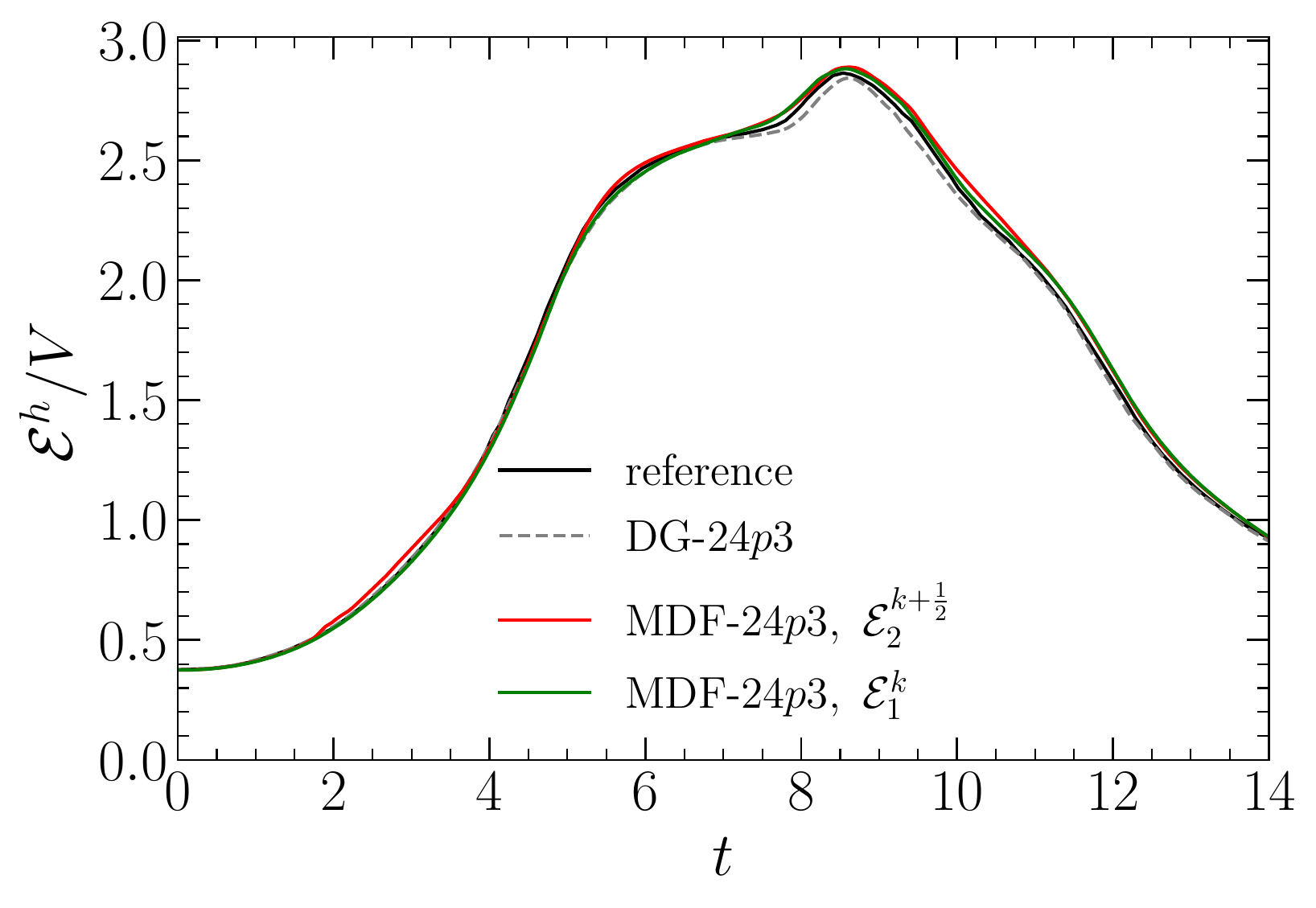}
			\end{minipage}
		}\\
		\subfloat{
			\begin{minipage}[b]{0.5\textwidth}
				\centering
				\includegraphics[width=0.76\linewidth]{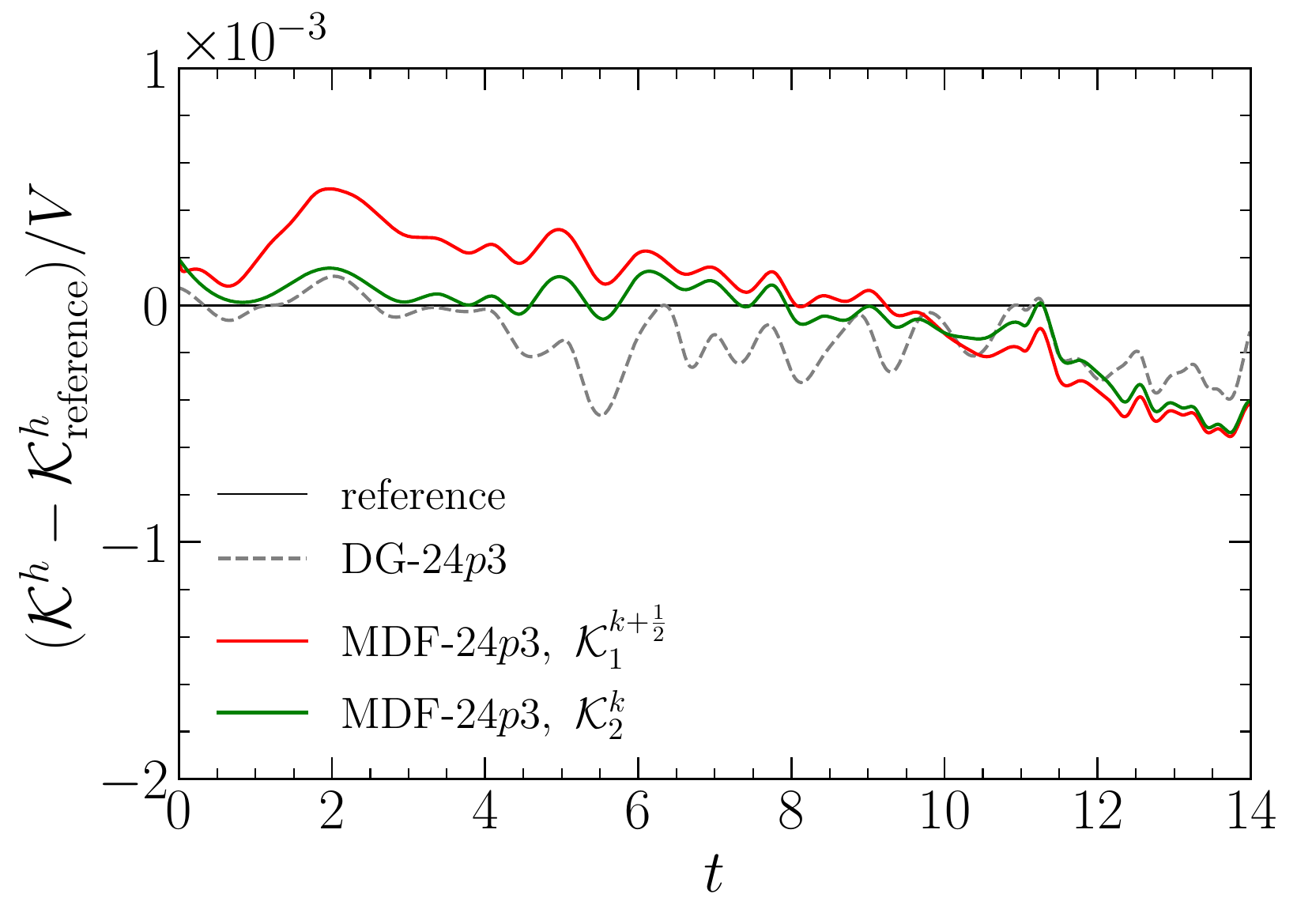}
			\end{minipage}
		}
		\subfloat{
			\begin{minipage}[b]{0.5\textwidth}
				\centering
				\includegraphics[width=0.77\linewidth]{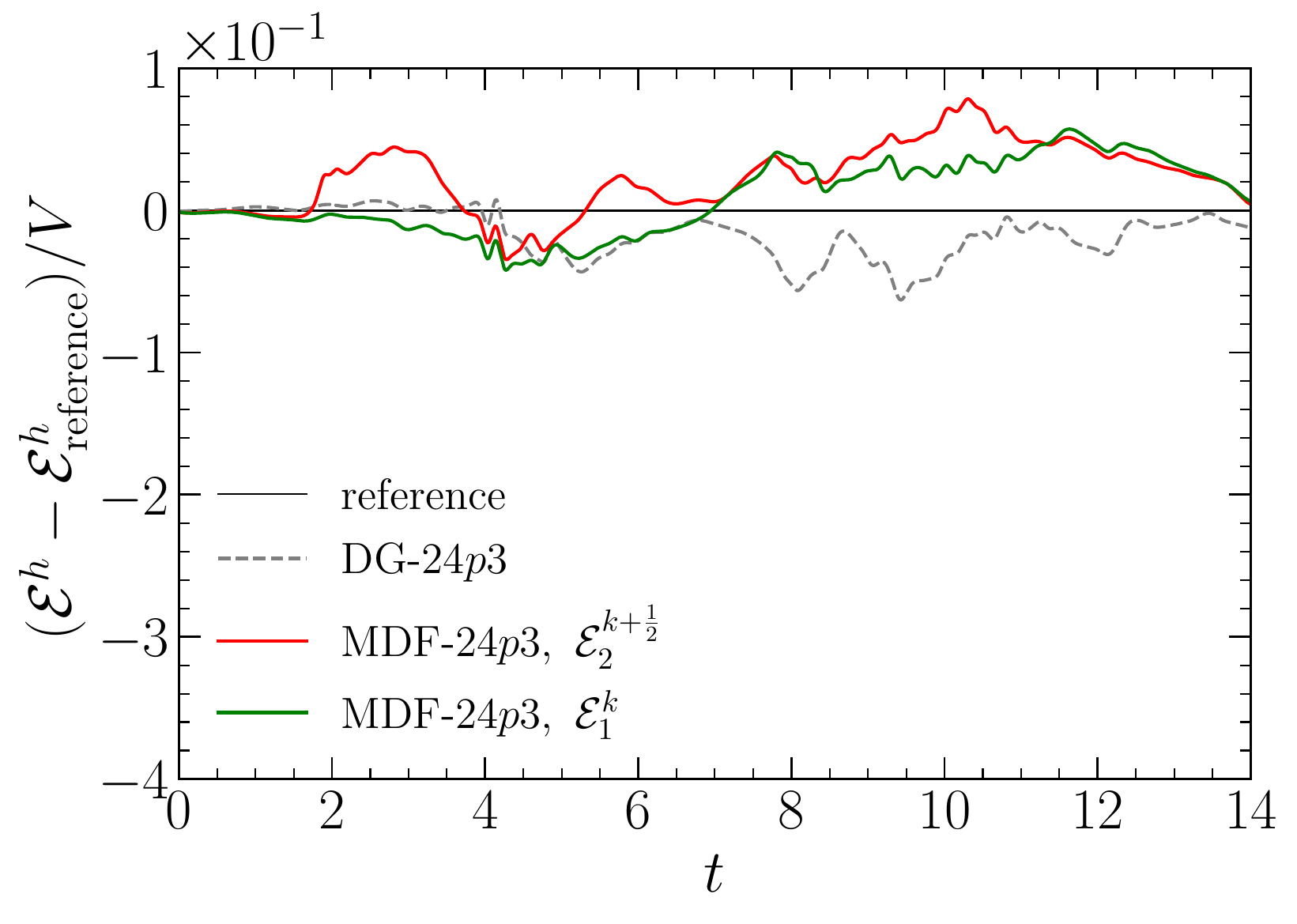}
			\end{minipage}
		}
		\caption{\REVT{A comparison of kinetic energy and enstrophy results of the TGV test. The reference and the DG-$ 24p3 $ results are taken from \cite{chapelier2012inviscid}. The reference method is a $ 128 $th order spectral method. DG stands for a discontinuous Galerkin method. MDF stands for the mimetic dual-field method. $ 24p3 $ represents $ h=1/24 $ and the degree of the polynomial spaces is 3. The time step interval is selected to be $ \varDelta t =1/50 $ for MDF-$24p3 $.}}
		\label{fig: TGV 2}}
\end{figure}

In Fig.~\ref{fig: TGV 1} and Fig.~\ref{fig: TGV 2}, results of total kinetic energy and total enstrophy are presented. These results are compared to benchmarks taken from \cite{chapelier2012inviscid}. 
In Fig.~\ref{fig: TGV 1} we can see that the proposed mimetic dual-field method, compared to a discontinuous Galerkin (DG) method of the same order $ (p=2) $ and in the same mesh ($ 32\times32\times32 $ elements), produces better results in terms of the error to the results produced by a reference, a very high ($ 128 $th) order spectral method. This is mostly clear in the enstrophy results near $ t=9 $ when the total enstrophy reaches its peak; the DG method is not able to capture the peak of the total enstrophy while the mimetic dual-field method captures it well for both of the dual solutions. \MOD{Similar} comparisons are made for more resolved simulations in Fig.~\ref{fig: TGV 2}, where improved results are seen especially near the peak of total enstrophy; the DG method now is able to capture the peak and the mimetic dual-field method captures the shape of the peak better.

\begin{figure}[h!]
	\centering
	\includegraphics[width=0.48\linewidth]{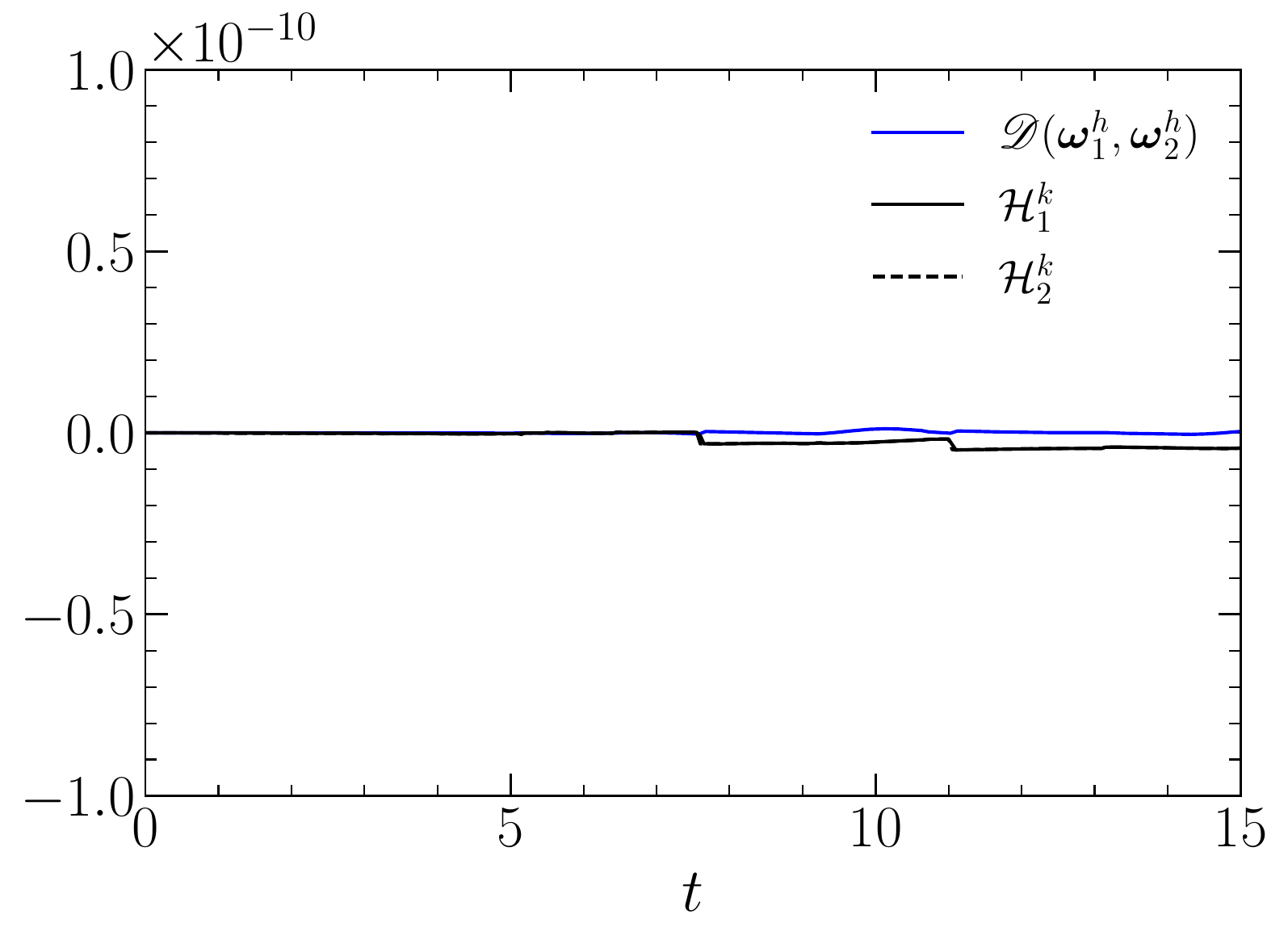}
	\caption{Total helicity and its dissipation rate $ \mathscr{D}(\bw_{1}^{h},\bw_{2}^{h}) $, see \eqref{Eq: discrete helicity disspation H1 H2 1}, versus time of the TGV test for MDF-$ 8p2 $. The time step interval is selected to be $ \varDelta t =1/20 $.}
	\label{fig:h1h2dw1w2}
\end{figure}
In Fig.~\ref{fig:h1h2dw1w2}, some results of the total helicity versus time for the TGV flow is shown. It is seen that, as the flow evolves, the total helicity remains zero (to the machine precision). Such a phenomenon is consistent with the fact that the dissipation rate of helicity, see \eqref{Eq: discrete helicity disspation H1 H2 1}, is constantly zero (to the machine precision)  as shown in the same diagram.

\begin{figure}[h!]
	\centering{
		\subfloat{
			\begin{minipage}[b]{0.5\textwidth}
				\centering
				\includegraphics[width=0.97\linewidth]{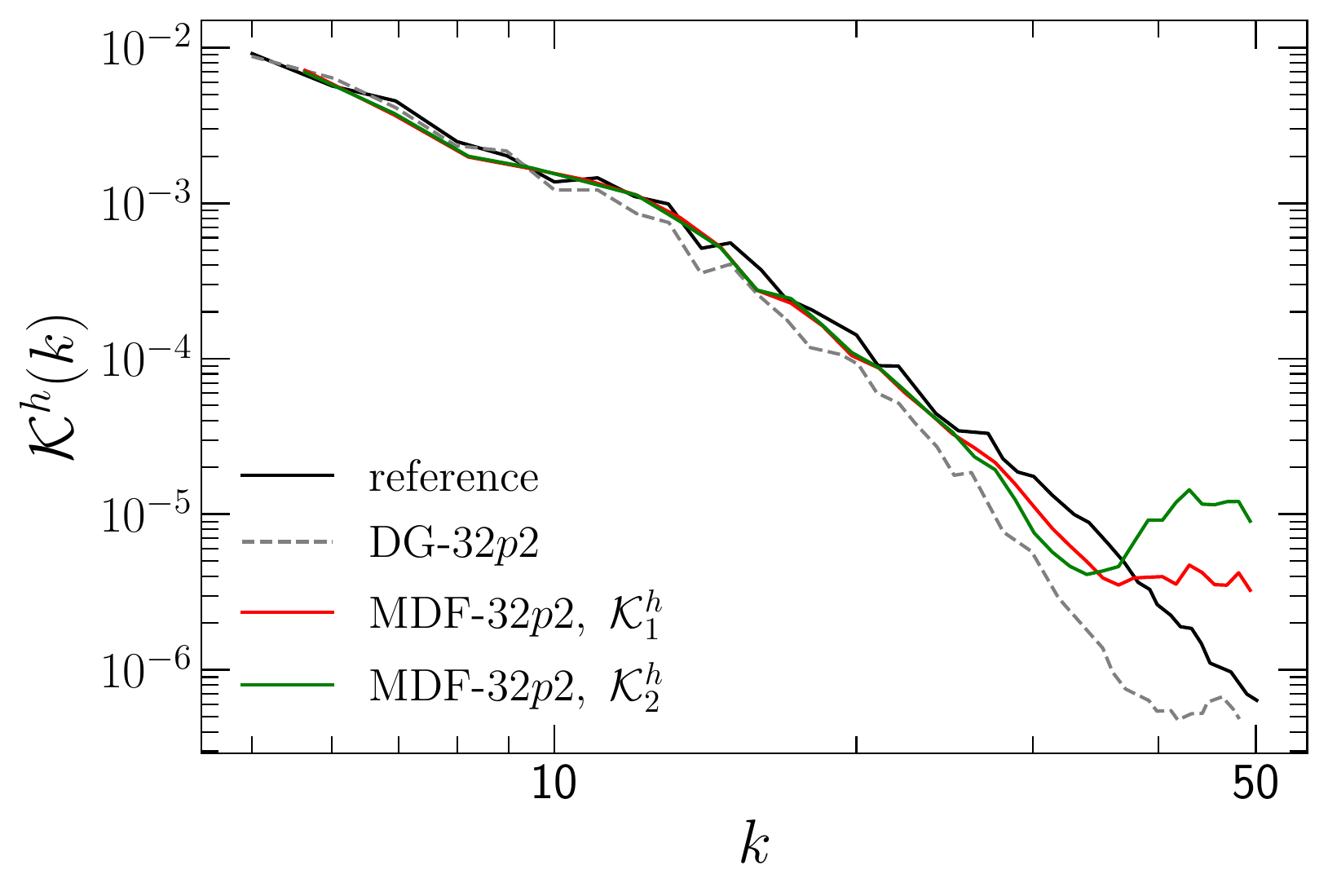}
			\end{minipage}
		}
		\subfloat{
			\begin{minipage}[b]{0.5\textwidth}
				\centering
				\includegraphics[width=0.97\linewidth]{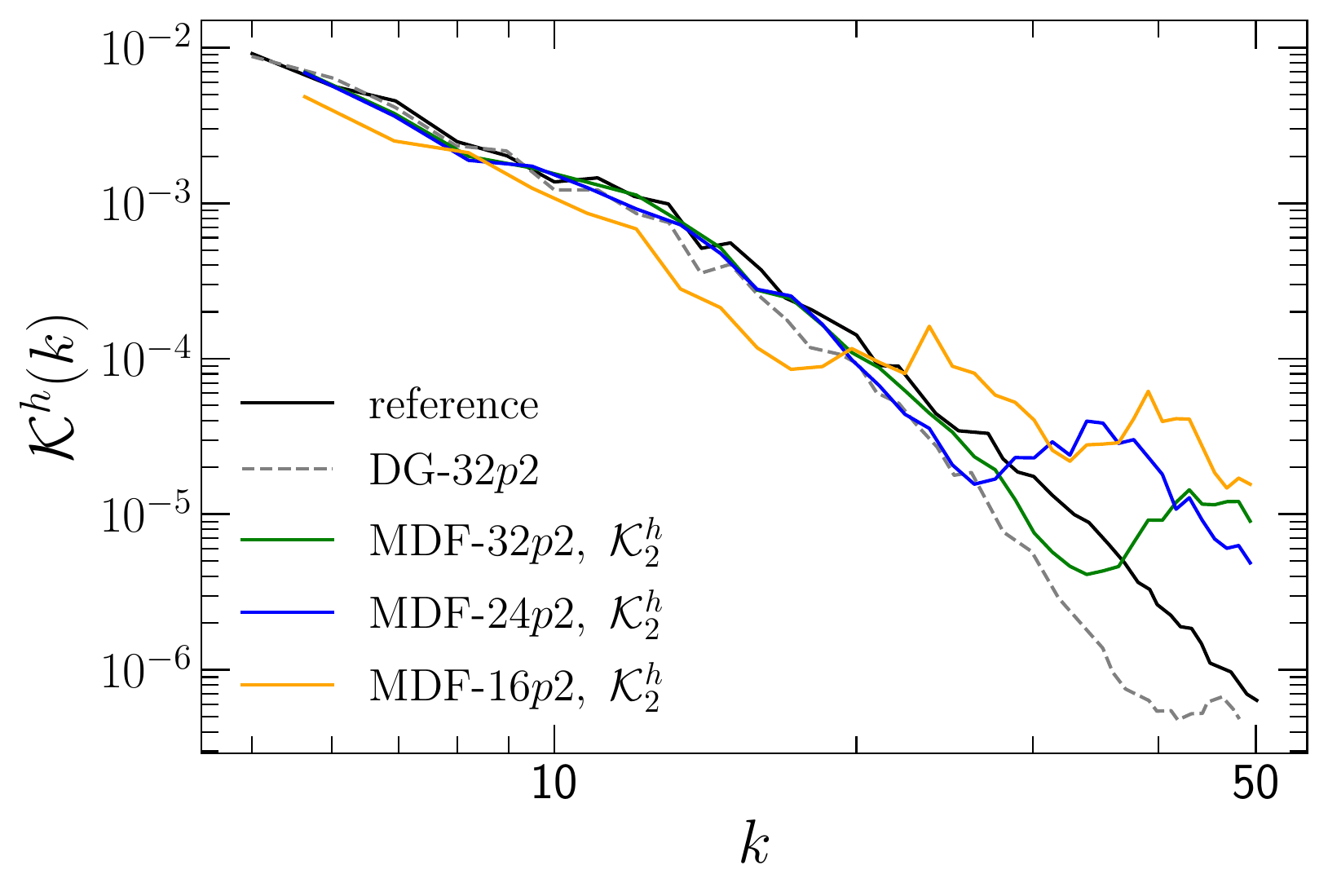}
			\end{minipage}
		}
		\caption{Kinetic energy spectra of the TGV test. The reference and the DG-$ 32p2 $ results are taken from \cite{chapelier2012inviscid}. The reference method is a $ 128 $th order spectral method. DG stands for a discontinuous Galerkin method. MDF stands for the mimetic dual-field method. $ 32p2 $, $ 24p2 $ and $ 16p2 $ represents that the mesh size $ h=1/32$, $ 1/24 $ and $1/16 $, respectively, and the degree of the polynomial spaces is 2. The time step interval is selected to be $ \varDelta t =1/50$,  $1/40 $ and $ 1/30 $ for MDF-$32p2 $, MDF-$24p2 $ and MDF-$16p2 $, \MOD{respectively}.}
		\label{fig:kespectrare500}}
\end{figure}

In Fig.~\ref{fig:kespectrare500}, the results of kinetic energy spectra at $ t=9.1 $ are presented. In the left diagram, it is seen that, in terms of kinetic energy, the mimetic dual-field method has similar accuracy as the DG method for large scales ($ k\leq 10 $). For medium scales ($ 10 < k \leq 35 $), both methods start to deviate from the high order spectral reference results with the proposed dual-field method showing less overdissipation. For small scales ($ k>35 $), both methods show large deviations from the reference results. The interesting aspect is that\MOD{, for small scales,} the DG method and the proposed dual-field method present different behaviors: the DG method over dissipates the energy and the dual-field method accumulates energy. The accumulation of energy at \MOD{small} scales is expected due to the energy conservation properties of the dual-field method. The energy cascade occurs up to the resolved scales and then it is stored (and accumulates at the smaller scales). It is the authors opinion that this can be an advantage of this method since subscale grid methods can specifically target these small scales and introduce the required dissipation that is not resolved. In opposition, the DG method already over dissipates the energy, therefore it is challenging for a dissipation based sub-grid scale model to improve the results for these smaller scales. This is a topic of interest for the authors and will be further researched in the future. A partial support for this claim is the results presented in the right diagram of Fig.~\ref{fig:kespectrare500} where \MOD{it is seen} that the value of $ k $ where energy accumulation starts decreases when a less resolved discretization is employed.

\section{Summary and future work}\label{Sec: Conclusions}
\subsection{Summary}
In this paper, we introduce \MOD{a discretization} which satisfies pointwise mass conservation and, if in the absence of dissipative terms, conserves total kinetic energy and total helicity and, otherwise, properly captures the dissipation rates of total kinetic energy and total helicity for the 3D incompressible Navier-Stokes equations. The discretization is based on a novel dual-field mixed weak formulation where two \MOD{evolution} equations are employed. \REVT{A staggered temporal discretization {linearizes the convective terms} and reduces the size of the discrete systems, which can be regarded as a big advantage of the proposed method in terms of the computational efficiency.} A mimetic spatial discretization enables the validity of the conservation properties and the dissipation rates at the fully discrete level.

\subsection{Future work}

In this paper, $\bu_1$ and $\bu_2$ ($\bw_1$ and $\bw_2$) are approximated in different function spaces, and, therefore, equality between them will not hold unless the flow is fully resolved. 
\REVO{In other words, we are not able to construct a square, time-independent and explicit discrete Hodge operator. Instead, this method implicitly defines a time-dependent discrete Hodge operator. By allowing the time evolution of the discrete Hodge operator we can construct a helicity conserving scheme.}
Based on the promising results reported in this paper, we want to apply the \emph{algebraic dual polynomial spaces}, \cite{JAIN2020}, such that solutions $\bu^{h}_1$ and $\bu^{h}_2$ ($\bw^{h}_1$ and $\bw^{h}_2$) are two representations in a pair of algebraic dual polynomial spaces. As a result, we expect the difference between  $\bu^h_1$ and $\bu^h_2$ ($\bw^{h}_1$ and $\bw^{h}_2$) to be smaller and using
the vorticity from the other subset of equations, see \eqref{Eq: TD1} and \eqref{Eq: TD2}, to be more consistent.

In the kinetic energy spectra of the dual field formulation, Fig.~\ref{fig:kespectrare500}, we see that for high wave numbers the energy decay is insufficient. This is attributed to the fact that the scheme is non-dissipative and the grids are too coarse for energy at the small scales to dissipate. In future work we want to add a sub-grid scale model on the momentum equations for $\bu_1$ and $\bu_2$ of the form $\epsilon \Delta (\bu_1 - \bu_2)$ to the $ \bu_1$ equation, \eqref{Eq: WF d} and $\epsilon \Delta (\bu_2 - \bu_1)$ to the $ \bu_2$ equation, \eqref{Eq: WF a}. If we define
\[ \bar{\bu} = \frac{1}{2} \left (  \bu_1 + \bu_2\right ) \;,\]
this sub-grid scale diffusion cancels from the average, while it only acts on the difference between the two fields
\[  \bu' = \frac{1}{2} \left (  \bu_1 - \bu_2\right ) \;.\]
So the numerical dissipation only acts on the difference of the two fields. This implies that the added diffusion is only active for the large wave numbers, \REVT{where the} difference between $\bu_1$ and $\bu_2$ is significant, while for the small wave numbers where $\bu_1$ and $\bu_2$  are almost the same, no dissipation takes place of $\bu'$. In this sense the dual field formulation can be used as a turbulence model. Future work needs to establish how the parameter $\epsilon$ should be chosen.

Other steps we want to report in the future includes, for example, error analysis, mesh adaptivity based on the local difference between $\bu^{h}_1$ and $\bu^{h}_2$ ($\bw^{h}_1$ and $\bw^{h}_2$) and the extension from periodic boundary conditions to general boundary conditions.

\section*{Acknowledgments}
Yi Zhang is supported by China Scholarship council under grant number 201607720010 and Leo G. Rebholz is supported by US National Science Foundation grant DMS2011490. The authors also would like to thank Dr. Chapelier for sharing reference results (spectral and DG methods) that enabled the comparison of the proposed to existing methods. \MOD{We also thank the reviewers for their valuable comments.}

\bibliographystyle{elsarticle-num}
\bibliography{ref}
\end{document}